\documentclass[reqno,centertags]{amsart}
\usepackage{amsmath,amsthm,amscd,amssymb}
\usepackage{latexsym}
\usepackage{upref}
\newcommand{\bbN}{{\mathbb{N}}}
\newcommand{\bbR}{{\mathbb{R}}}

\newcommand{\bbZ}{{\mathbb{Z}}}
\newcommand{\bbC}{{\mathbb{C}}}
\newcommand{\calA}{{\mathcal A}}
\newcommand{\calB}{{\mathcal B}}

\newcommand{\calD}{{\mathcal D}}
\newcommand{\calE}{{\mathcal E}}

\newcommand{\calH}{{\mathcal H}}
\newcommand{\calI}{{\mathcal I}}

\newcommand{\calM}{{\mathcal M}}
\newcommand{\calN}{{\mathcal N}}

\newcommand{\calU}{{\mathcal U}}

\newcommand{\no}{\nonumber}
\newcommand{\lb}{\label}

\newcommand{\ul}{\underline}
\newcommand{\ol}{\overline}
\newcommand{\ti}{\tilde  }
\newcommand{\wti}{\widetilde  }
\newcommand{\Oh}{O}
\newcommand{\oh}{o}

\newcommand{\loc}{\text{\rm{loc}}}
\newcommand{\rank}{\text{\rm{rank}}}
\newcommand{\ran}{\text{\rm{ran}}}

\newcommand{\supp}{\text{\rm{supp}}}

\newcommand{\bi}{\bibitem}
\newcommand{\hatt}{\widehat}

\renewcommand{\Re}{\text{\rm Re}}
\renewcommand{\Im}{\text{\rm Im}}

\DeclareMathOperator{\tr}{tr}

\DeclareMathOperator*{\slim}{s-lim}
\DeclareMathOperator*{\wlim}{w-lim}
\allowdisplaybreaks
\numberwithin{equation}{section}

\newtheorem{theorem}{Theorem}[section]
\newtheorem{lemma}[theorem]{Lemma}
\newtheorem{corollary}[theorem]{Corollary}

\theoremstyle{definition}
\newtheorem{definition}[theorem]{Definition}
\newtheorem{example}[theorem]{Example}

\theoremstyle{remark}
\newtheorem{remark}[theorem]{Remark}

\begin{document}
\title[Herglotz Functions]{On Matrix-Valued Herglotz
Functions}
\author[Gesztesy and Tsekanovskii]{Fritz Gesztesy
and Eduard
Tsekanovskii}
\address{Department of Mathematics, University of
Missouri, Columbia, MO
65211, USA}
\email{fritz@math.missouri.edu\newline
\indent{\it URL:}
http://www.math.missouri.edu/people/faculty/fgesztesy.html}
\address{Department of Mathematics, University of Missouri,
Columbia, MO
65211, USA}
\email{tsekanov@math.missouri.edu}
\thanks{Research supported by the US National Science
Foundation under
Grant No.~DMS-9623121.}
\date{\today}
\subjclass{Primary 30D50, 30E20, 47A10; Secondary 47A45.}

\begin{abstract}

We provide a comprehensive analysis of matrix-valued
Herglotz functions and illustrate their applications
in the spectral
theory of self-adjoint Hamiltonian systems including
matrix-valued
Schr\"odinger and Dirac-type operators. Special emphasis
is devoted to
appropriate matrix-valued extensions of the well-known
Aronszajn-Donoghue theory concerning
support properties of measures in their
Nevanlinna-Riesz-Herglotz
representation. In
particular, we study a class of linear fractional
transformations
$M_A(z)$ of a given $n \times n$ Herglotz matrix $M(z)$
and prove that
the minimal support of the absolutely continuos part of
the measure
associated to $M_A(z)$ is invariant under these linear
fractional
transformations.

Additional applications discussed in detail include
self-adjoint
finite-rank perturbations of self-adjoint operators,
self-adjoint
extensions of densely defined symmetric linear operators
(especially,
Friedrichs and Krein extensions), model operators for
these two
cases, and
associated realization theorems for certain classes of
Herglotz
matrices.

\end{abstract}

\maketitle
\section{Introduction}\lb{s1}

The spectral analysis of self-adjoint ordinary
differential operators,
or
more generally, that of self-adjoint Hamiltonian systems
(including
matrix-valued Schr\"odinger and Dirac-type operators),
is well-known
to
be intimately connected with the
Nevanlinna-Riesz-Herglotz
representation
of matrix-valued Herglotz functions. The latter
terminology is not
uniformly adopted in the literature and postponing
its somewhat
controversial
origin to the beginning of Section~\ref{s2}, we recall
that $M(z)$
is said
to be a matrix-valued Herglotz function if
$M:\bbC_+\to M_n(\bbC)$ is
analytic and $\Im(M(z))\geq 0$ for $z\in\bbC_+$.
(Here $\bbC_+$
denotes the
open complex upper half-plane, $M_n(\bbC)$,
$n\in\bbN$ the set of
$n \times n$ matrices with entries in $\bbC$, and
$\Re(M)=(M+M^*)/2$,
$\Im(M(z))=(M-M^*)/(2i)$ the real and imaginary parts
of the
matrix $M$).
The Nevanlinna-Riesz-Herglotz representation of $M(z)$
is of the type
\begin{equation}
M(z)=C+Dz+\int_\bbR d\Omega(\lambda)
((\lambda-z)^{-1}-\lambda
(1+\lambda^2)^{-1}), \lb{1.1}
\end{equation}
where
\begin{equation}
C=\Re(M(i)), \quad D=\lim_{\eta\uparrow\infty}(\frac{1}
{i\eta}M(i\eta))\geq 0,
\lb{1.2}
\end{equation}
and $d\Omega(\lambda)$ denotes an $n \times n$
matrix-valued measure
satisfying
\begin{equation}
\int_\bbR (\ul c,d\Omega(\lambda)\ul c)_{\bbC^n}
(1+\lambda^2)^{-1}<\infty
\text{ for all } \ul c\in \bbC^n \lb{1.3}
\end{equation}
(with $(\cdot,\cdot)_{\bbC^n}$ the scalar product
in $\bbC^n$).
The Stieltjes
inversion formula for $\Omega$ then reads
\begin{equation}
\frac{1}{2}\Omega(\{\lambda_1\})+\frac{1}{2}
\Omega(\{\lambda_2\}) +
\Omega((\lambda_1,\lambda_2)) =\pi^{-1}
\lim_{\varepsilon\downarrow 0}
\int_{\lambda_1}^{\lambda_2} d\lambda \,
\Im(M(\lambda+i\varepsilon))
\lb{1.4}
\end{equation}
and its absolutely continuous part $\Omega_{ac}$
(w.r.t.~Lebesgue
measure)
is given by
\begin{equation}
d\Omega_{ac}(\lambda)=\pi^{-1}
\Im(M(\lambda+i\varepsilon))d\lambda
\lb{1.5}
\end{equation}
(cf. Section~\ref{s5} for a detailed exposition of
these facts).
Spectral
analysis of ordinary differential operators (with
matrix-valued
coefficients)
then boils down to an analysis of (matrix-valued) measures
$d\Omega(\lambda)$
in the representation \eqref{1.1} for $M(z)$.
These Herglotz
matrices are
traditionally called Weyl-Titchmarsh M-functions
in honor of
Weyl, who
introduced the concept in the special (scalar)
Sturm-Liouville
case, and
Titchmarsh, who recognized and first employed its
function-theoretic
content. Since different self-adjoint boundary conditions
associated to a
given formally symmetric (matrix-valued) differential
expression
$\tau$
yield self-adjoint realizations of $\tau$ whose
corresponding
$M$-functions
are related via certain linear fractional
transformations (cf.
\cite{KS74}),
we study in
depth transformations of the type
\begin{equation}
M(z) \longrightarrow M_A(z)=(A_{2,1}+A_{2,2}M(z))(A_{1,1}
+A_{1,2}M(z))^{-1},
\lb{1.6}
\end{equation}
where
\begin{align}
&A=\big(A_{p,q}\big)_{1\leq p,q\leq 2} \in{\calA}_{2n},
\no \\
&{\calA}_{2n}=\{A\in M_{2n}(\bbC)\,|\, A^*J_{2n}A
=J_{2n} \},
\quad J_{2n}=\left(\begin{array}{cc}0 & -I_n\\ I_n&
0\end{array}\right)
\lb{1.7}
\end{align}
($I_n$ the identity matrix in $\bbC^n$). $M_A(z)$,
$A\in{\calA}_{2n}$ are
Herglotz matrices whenever $M(z)$ is a Herglotz matrix.
Moreover,
denoting
the measure in the Nevanlinna-Riesz-Herglotz representation
\eqref{1.1}
for $M_A(z)$ by $d\Omega_A(\lambda)$, we provide a
matrix-valued
extension of the well-known
Aronszajn-Donoghue theory relating support properties of
$d\Omega_A(\lambda)$
and $d\Omega(\lambda)$, originally
inspired by Sturm-Liouville boundary value problems. As
one of our
principal new results we prove that the minimal support
of the
absolutely
continuous part $d\Omega_{A,ac}(\lambda)$ of
$d\Omega_A(\lambda)$
is
independent of $A\in{\calA}_{2n}$, which represents
the proper
generalization
of Aronszajn's celebrated result \cite{Ar57} for
Sturm-Liouville
operators.

Concrete applications of our formalism include self-adjoint
finite-rank
perturbations of self-adjoint operators and self-adjoint
extensions
of
densely defined symmetric closed linear operators $H$
with finite
deficiency
indices especially emphasizing Friedrichs and Krein
extensions
in the
special case where $H$ is bounded from below. Moreover,
we describe
in detail
associated model operators and realization theorems
for certain
classes of
Herglotz matrices. These results appear to be of
independent
interest in
operator theory.

Finally we briefly describe the content of each
section. In
Section~\ref{s2} we review basic facts on scalar
Herglotz functions
and their representation theorems. Section~\ref{s3}
reviews the
Aronszajn-Donoghue theory concerning support
properties of
$d\Omega_A(\lambda)$ in the scalar case and
Section~\ref{s4}
describes
a variety of applications of the scalar formalism.
Some of these
applications to self-adjoint extensions of symmetric
operators with
deficiency indices $(1,1)$ (such as
Theorem~\ref{t4.1}\,(iv) and
Theorems~\ref{t4.3} and \ref{t4.5}) appear to be new.
Section~\ref{s5}
provides the necessary background material for
matrix-valued
Herglotz
functions and their representation theorems.
Section~\ref{s6}, the
principal section of this paper, is devoted to a
detailed study of
support properties of $d\Omega_A(\lambda)$ and
Theorem~\ref{t6.6}
contains the invariance result of the minimal
support of
$d\Omega_{A,ac}(\lambda)$ with respect to
$A\in{\calA}_{2n}$. In our
final Section~\ref{s7} we again treat applications
to self-adjoint
finite-rank perturbations and self-adjoint extensions
of symmetric
operators with finite deficiency indices. We pay
particular attention
to Friedrichs and Krein extensions of symmetric operators
bounded from
below and prove a variety of realization theorems for
certain classes
of Herglotz matrices. To the best of our knowledge, most of
the applications discussed in Section~\ref{s7} are new.
For the
convenience of the reader
we collect some examples of scalar Herglotz functions in
Appendix~\ref{A}. Appendix~\ref{B} contains a detailed
discussion of
Krein's formula, relating self-adjoint extensions of a
symmetric
operator, and its application to linear fractional
transformations
of associated Weyl-Titchmarsh matrices.

It was our aim to provide a rather comprehensive account on
matrix-valued Herglotz functions. We hope the enormous
number of
applications of this formalism to the theory of self-adjoint
extensions of symmetric operators,
the spectral theory of ordinary (matrix-valued)
differential and
difference operators, interpolation problems, and
factorizations of matrix and operator functions
\cite{AK97}, \cite{At64}--\cite{BGK79},
\cite{Be68}, \cite{Br71},
\cite{BEME96},
\cite{Bu97}, \cite{Cl90}, \cite{DM91}-- \cite{DMT88},
\cite{DSS90}, \cite{Dy89}, \cite{Ev64}--\cite{Fu76},
\cite{GKS97a},
\cite{GKS97b},
\cite{HKS97}--\cite{HS86}, \cite{Jo87}, \cite{Ko50},
\cite{KR74},
\cite{KS88}--\cite{Kr95}, \cite{Li57}--\cite{Ma97},
\cite{MHR77},
\cite{Ni72}, \cite{Or76}, \cite{Wa74}, \cite{We87},
inverse spectral
theory \cite{AM63}, \cite{AG94}--\cite{AGS96}, \cite{Ga68},
\cite{GH97}, \cite{GHL98}, \cite{GR98}, \cite{JL97},
\cite{Ma92}--\cite{Ma97}, \cite{MO82}, \cite{MPT95},
\cite{Sa92}--\cite{Sa97},
\cite{Su90}, \cite{SF70},
\cite{WK74}, \cite{Yu92}--\cite{Zh95},  and completely
integrable
hierarchies of matrix-valued nonlinear evolution equations
\cite{AK90},
\cite{Ch96}, \cite{Di91}, \cite{Du83}, \cite{DMN76},
\cite{GD77},
\cite{IKK94}, \cite{Ma78}, \cite{Ma88}, \cite{MM79},
justifies this effort.

\section{Basic Facts on Scalar Herglotz Functions} \lb{s2}

This section provides a quick review of scalar Herglotz
functions
and
their representation theorems. These results are well-known,
in fact,
classical by now, and we include them for later reference to
achieve
a certain
degree of completeness, and partly to fix our notation.

\begin{definition}\label{d2.1} Let $\bbC_\pm=\{z\in\bbC
\mid \Im(z)\gtrless 0 \}$.
$m:\mathbb{C_+}\to \mathbb{C}$ is called a Herglotz
function if $m$ is analytic on $\mathbb{C}_+$ and
$m(\mathbb{C}_+)\subseteq \mathbb{C}_+$.
\end{definition}

It is customary to extend $m$ to $\bbC_-$ by reflection,
that is,
one defines
\begin{equation}
m(z)=\overline{m(\overline z)},
\quad z\in\mathbb{C}_-. \lb{2.1a}
\end{equation}
We will adopt this convention in this paper. While
\eqref{2.1a}
defines
an analytic function on $\bbC_-$, $m\big|_{\bbC_-}$ in
\eqref{2.1a}, in
general, does not represent the analytic continuation of
$m\big|_{\bbC_+}$ (cf.~Lemma~\ref{l2.5} for more details in
this connection).

There appears to be considerable disagreement concerning the
proper name
of functions satisfying the conditions in Definition~\ref{d2.1}.
In
addition to the presently used notion of Herglotz functions one
can also find
the names Pick, Nevanlinna, Nevanlinna-Pick, and $R$-functions
(sometimes
depending on the geographical origin of authors and occasionally
whether
the open upper half-plane $\mathbb{C}_+$ or the conformally
equivalent
open unit disk $D$ is involved). Following a tradition in
mathematical
physics, we decided to adopt the notion of Herglotz
functions
in this paper.

If $m(z)$ and $n(z)$ are Herglotz functions, then
$m(z)+n(z)$
and
$m(n(z))$ are also Herglotz. Elementary examples of
Herglotz
functions are
\begin{align}&c+id,\quad c+dz,\quad c\in \mathbb{R},
\quad d>0,
\label{2.1}\\
&z^r ,\quad 0< r <1,\label{2.2}\\
&\ln(z),\label{2.3}
\end{align}
choosing the obvious branches in
\eqref{2.2} and \eqref{2.3},
\begin{align}
&\tan (z), \quad -\cot(z), \label{2.4}\\
&\frac{a_{2,1}+a_{2,2}z}{a_{1,1}+a_{1,2}z},\lb{2.5}\\
&a=\left(\begin{array}{cc}a_{1,1} & a_{1,2}\\ a_{2,1} &
a_{2,2}\end{array}\right) \in M_2(\mathbb{C}),
\quad a^*j_2a=j_2,
\quad
j_2=\left(\begin{array}{cc}0 & -1\\ 1&
0\end{array}\right),\lb{2.6}
\end{align}
with $M_n(\mathbb{C})$ the
set of $n\times n$ matrices with entries in $\mathbb{C}$,
and hence
\begin{equation}
-1/z\lb{2.7}
\end{equation}
as a special case of
\eqref{2.5}. Equations \eqref{2.5} and \eqref{2.6} define
the group
of automorphisms on $\bbC_+$ (or $\bbC_-$). Finally,
we mention a
less elementary example, the digamma
function \cite{AS72}, Ch. 6,
\begin{equation}
\psi (z)=\Gamma '(z)/\Gamma(z),\lb{2.8}
\end{equation}
with $\Gamma(z)$ Euler's gamma function. Further
examples are
described
in detail in Appendix~\ref{A}.

As a consequence,
\begin{align}
&-1/m(z),\quad m(-1/z),\lb{2.9}\\
&\ln(m(z)),\lb{2.10}
\end{align}
and
\begin{equation}
m_a(z)=\frac{a_{2,1}+a_{2,2}m(z)}{a_{1,1}+a_{1,2}m(z)},
\lb{2.11}
\end{equation}
with $a\in M_2(\mathbb{C})$ satisfying
\eqref{2.6}, are all Herglotz functions whenever $m(z)$
is Herglotz.
More generally, and most relevant in the context of spectral
theory for
linear operators, let $H$ be a self-adjoint operator
in a (complex,
separable) Hilbert space
$\mathcal{H}$ with $(\cdot ,\cdot )_{\mathcal{H}}$ the scalar
product on
$\mathcal{H}\times \mathcal{H}$ linear in the second factor. Let
$(H-z)^{-1}$, $z\in \mathbb{C}\backslash\mathbb{R}$
denote the
resolvent of $H$.
Then for all $f\in
\mathcal{H}$,
\begin{equation}
(f,(H-z)^{-1}f)_{\mathcal{H}}\lb{2.12}
\end{equation}
is a scalar Herglotz function (it suffices to appeal
to the
spectral theorem
for $H$ and apply the functional calculus to
$(H-z)^{-1}$), while
$(H-z)^{-1}$ represents a $\calB (\calH)$-valued Herglotz
function
($\calB (\calH)$ the set of bounded linear operators
mapping
$\calH$ to
itself).

The fundamental result on Herglotz functions and their
representations on
Borel transforms, in part due to Fatou, Herglotz, Luzin,
Nevanlinna,
Plessner, Privalov, de la Vall{\'e}e Poussin, Riesz, and
others,
then reads as follows.

\begin{theorem}
\mbox{\rm (\cite{AG93}, Ch.~VI, \cite{AD56},
\cite{Do74}, Chs.~II, IV,
\cite{KK74}, \cite{Ko80}, Ch.~6, \cite{Pr56}, Chs.~II, } \\
\mbox{\rm IV, \cite{RR94}, Ch.~5).} \label{t2.2}
Let $m(z)$ be a Herglotz function. Then

\noindent (i). $m(z)$ has finite normal limits $m(\lambda
\pm i0)=\lim_{\varepsilon
\downarrow 0} m(\lambda \pm i\varepsilon)$ for
a.e.~$\lambda \in \mathbb{R}$.

\noindent (ii).~If $m(z) $ has a zero normal limit on a
subset of
$\mathbb{R}$ having positive Lebesgue measure, then
$m\equiv 0$.

\noindent (iii). There exists a Borel measure $\omega$ on
$\mathbb{R}$
satisfying
\begin{equation}\lb{2.13}
\int_{\mathbb{R}}d\omega (\lambda )(1+\lambda
^2)^{-1}<\infty
\end{equation}
such that the Nevanlinna, respectively, Riesz-Herglotz
representation
\begin{align}
m(z)&=c+dz+\int_{\mathbb{R}}
d\omega (\lambda)((\lambda -z)^{-1}-\lambda
(1+\lambda^2)^{-1}),
\quad z\in\bbC_+, \lb{2.14} \\
c&=Re(m(i)),\quad d=\lim_{\eta \uparrow
\infty}m(i\eta )/(i\eta )
\geq 0 \no
\end{align}
holds.

\noindent (iv). Let $(\lambda _1,\lambda_2)\subset
\mathbb{R}$,
then the
Stieltjes inversion formula for
$\omega$ reads
\begin{equation}\lb{2.15}
\frac{1}{2}\omega\left(\left\{\lambda_1
\right\}\right)+\frac{1}{2}\omega
\left(\left\{\lambda_2\right\}\right) +\omega
((\lambda_1,\lambda_2))=\pi^{-1}\lim_{\varepsilon\downarrow
0}\int^{\lambda_2}_{\lambda_1}d\lambda \, \Im(m(\lambda
+i\varepsilon)).
\end{equation}

\noindent (v). The absolutely continuous $({\it ac})$ part
$\omega_{ac}$ of
$\omega$ with respect to Lebesgue measure $d\lambda$ on
$\mathbb{R}$ is
given by
\begin{equation}\lb{2.16} d\omega_{ac}(\lambda)=\pi^{-1}\Im
(m(\lambda+i0))d\lambda .\end{equation}

\noindent (vi). Any poles and isolated zeros of $m$ are simple
and located
on the real axis, the residues at poles being
negative.
\end{theorem}

It is quite illustrative to compare the various measures
$\omega$ for
the examples in \eqref{2.1}--\eqref{2.5}, \eqref{2.7}, and
\eqref{2.8}
and hence we provide these details and also a few more
sophisticated
examples in Appendix A.

Further properties of Herglotz functions are collected
in the
following
theorem. We denote by
$\omega =\omega_{ac}+\omega_s=\omega_{ac}+\omega_{sc}
+\omega_{pp}$
the
decomposition of $\omega$ into its absolutely continuous
$({\it ac})$,
singularly continuous $({\it sc})$, pure point $({\it pp})$,
and singular
$({\it s})$ parts
with respect to Lebesgue measure on $\mathbb{R}$.

\begin{theorem} [\cite{AD56}, \cite{KK74}, \cite{Si95a},
\cite{Si95b}] \label{t2.3} Let
$m(z)$ be a Herglotz function with representation \eqref{2.14}.
Then

\noindent (i).
\begin{align}
& d=0 \text{ and } \int_{\mathbb{R}}d\omega
(\lambda)(1+|\lambda|^s)^{-1}<\infty \text{ for some }
s\in (0,2)
\no \\
& \text{if and only if }
\int^\infty_1d\eta \, \eta^{-s} \, \Im (m(i\eta ))<\infty.
\lb{2.17}
\end{align}

\noindent (ii). Let $(\lambda_1,\lambda_2)\subset \mathbb{R}$,
$\eta_1>0$.
Then there is a constant
$C(\lambda_1,\lambda_2,\eta_1)>0$ such that
\begin{equation}\lb{2.18}
\eta|m(\lambda+i\eta)|\leq
C(\lambda_1,\lambda_2,\eta_1),\quad (\lambda,\eta)\in
[\lambda_1,\lambda_2]\times (0,\eta_1).\end{equation}

\noindent (iii). \begin{align} \lb{2.19}
&\sup_{\eta
>0}\eta |m(i\eta)|<\infty \text{ if and only if }
m(z)=\int_{\mathbb{R}}d\omega (\lambda)(\lambda-z)^{-1} \no \\
& \hspace*{5cm} \text{ and }
\int_{\mathbb{R}}d\omega(\lambda)<\infty.
\end{align}
In this case,
\begin{equation}\lb{2.20}
\int_{\mathbb{R}}d\omega
(\lambda)=\sup_{\eta >0}\eta |m(i\eta )|=-i\lim_{\eta \uparrow
\infty}\eta m(i\eta).
\end{equation}

\noindent (iv). For all $\lambda \in \mathbb{R}$,
\begin{align}\lb{2.21}
&\lim_{\varepsilon\downarrow 0}\varepsilon
\Re(m(\lambda+i\varepsilon ))=0,\\
&\omega(\{\lambda\})=\lim_{\varepsilon\downarrow 0}\varepsilon
\Im
(m(\lambda+i\varepsilon))=-i\lim_{\varepsilon \downarrow 0}
\varepsilon
m(\lambda+i\varepsilon ).\lb{2.22}
\end{align}

\noindent (v). Let $L>0$ and suppose $0\leq \Im (m(z))\leq L$
for all $z\in
\mathbb{C}_+$. Then $d=0$, $\omega$ is purely absolutely
continuous,
$\omega = \omega_{ac}$, and
\begin{equation} \lb{2.23}
0\leq
\frac{d\omega(\lambda)}{d\lambda}=\pi^{-1}\lim_{\varepsilon
\downarrow 0}\Im
(m(\lambda+i\varepsilon))\leq \pi^{-1}L \text{ for a.e. }
\lambda\in
\mathbb{R}.
\end{equation}

\noindent (vi). Let $p\in (1,\infty)$, $[\lambda_3,\lambda_4]
\subset
(\lambda_1,\lambda_2)$, $[\lambda_1,\lambda_2]\subset
(\lambda_5,\lambda_6)$. If
\begin{equation} \lb{2.24}
\sup_{0<\varepsilon <1}\int^{\lambda_2}_{\lambda_1}d\lambda
|\Im (m(\lambda
+i\varepsilon))|^p < \infty,
\end{equation}
then $\omega=\omega _{ac}$ is purely
absolutely continuous on $(\lambda_1,\lambda_2)$,
$\frac{d\omega_{ac}}{d\lambda}\in L^p((\lambda_1,
\lambda_2);d\lambda)$, and
\begin{equation} \lb{2.25}
\lim_{\varepsilon\downarrow 0}\|\pi^{-1}\Im
(m(\cdot +i\varepsilon ))-\frac{d\omega
_{ac}}{d\lambda}\|_{L^p((\lambda_3,\lambda_4);d\lambda)}=0.
\end{equation}
Conversely, if $\omega$ is purely absolutely continuous on
$(\lambda_5,\lambda_6)$, and if
$\frac{d\omega_{ac}}{d\lambda}\in$ \linebreak
$L^p((\lambda_5,\lambda_6);d\lambda)$,
then \eqref{2.24} holds.

\noindent (vii). Let $(\lambda_1,\lambda_2)
\subset\mathbb{R}$.
Then a local version
of Wiener's theorem reads for $p\in (1,\infty)$,
\begin{align}
&\lim_{\varepsilon\downarrow 0}\varepsilon^{p-1}
\int^{\lambda_2}_{\lambda_1}d\lambda |\Im
(m(\lambda+i\varepsilon ))|^p \no \\
&=\frac{\Gamma (\frac{1}{2})\Gamma (p-\frac{1}{2})}
{\Gamma (p)} \bigg( \frac{1}{2}\omega (\{\lambda_1\})^p
+\frac{1}{2}\omega (\{\lambda_2\})^p+
\sum_{\lambda\in (\lambda_1,\lambda_2)}
\omega (\{\lambda \})^p
\bigg). \lb{2.26}
\end{align}
Moreover, for $0<p<1$,
\begin{equation}\lb{2.27}
\lim_{\varepsilon \downarrow
0}\int^{\lambda_2}_{\lambda_1}d\lambda |\pi^{-1}\Im
(m(\lambda + i\varepsilon))|^p
=\int^{\lambda_2}_{\lambda_1}d\lambda
\left|\frac{d\omega_{ac}(\lambda)}{d\lambda}\right|^p.
\end{equation}
\end{theorem}

It should be stressed that Theorems~\ref{t2.2} and
\ref{t2.3}
record only
the tip of an iceberg of results in this area. In
addition
to the
references already mentioned, the reader will find a
great deal of
interesting results, for instance, in \cite{AD64},
\cite{RJL98}, \cite{RST94}, \cite{DT85},
\cite{Do74}, \cite{Gi84}, \cite{HKS97}--\cite{HSW97},
\cite{No60},
Ch.~III, \cite{Pe94}, \cite{SW86}.

Together with $m(z)$, $\ln (m(z))$ is a Herglotz
function by
\eqref{2.10}. Moreover, since
\begin{equation}\lb{2.28}
0\leq \Im (\ln (m(z))=\arg (m(z))\leq
\pi,\quad z\in \mathbb{C}_+,
\end{equation}
the measure $\hatt{\omega}$ in
the representation \eqref{2.14} of $\ln (m(z))$, that is,
in the
exponential representation of $m(z)$, is purely absolutely
continuous by
Theorem~\ref{t2.3}\,(v), $d\hatt{\omega}(\lambda)=\xi
(\lambda )d\lambda$
for some $0\leq\xi\leq 1$.
These exponential representations have been studied in
detail
by Aronszajn
and Donoghue \cite{AD56},
\cite{AD64} and we record a few of their properties below.

\begin{theorem}
[\cite{AD56}, \cite{AD64}] \label{t2.4}
Suppose $m(z)$ is a Herglotz function with representation
\eqref{2.14}. Then

\noindent (i). There exists a $\xi \in L^\infty
(\mathbb{R})$, $0\leq \xi \leq 1$ a.e., such that
\begin{align}
\ln(m(z))&=k+\int_{\mathbb{R}}d\lambda \,
\xi(\lambda)((\lambda-z)^{-1}-\lambda
(1+\lambda^2)^{-1}), \quad z\in\bbC_+, \lb{2.29} \\
k&=\Re (\ln(m(i))), \no
\end{align}
where
\begin{equation}\lb{2.30}
\xi (\lambda)=\pi^{-1}\lim_{\varepsilon
\downarrow 0}\Im (\ln (m(\lambda+i\varepsilon)))
\text{ a.e.}
\end{equation}

\noindent (ii). Let $\ell_1,\ell_2\in\mathbb{N}$ and
$d=0$ in
\eqref{2.14}. Then
\begin{align}
&\int^0_{-\infty}d\lambda \, \xi (\lambda )|\lambda
|^{\ell_1}(1+\lambda^2)^{-1}+\int^\infty_0d\lambda \,
\xi (\lambda)|\lambda
|^{\ell_2}(1+\lambda^2)^{-1}<\infty \text{ if and
only if}\no \\
& \int^0_{-\infty}d\omega
(\lambda)|\lambda|^{\ell_1}(1+\lambda^2)^{-1}
+\int^\infty_0d\omega
(\lambda)|\lambda|^{\ell_2}(1+\lambda^2)^{-1}<\infty \no \\
& \text{and } \lim_{z\to
i\infty}m(z)=c-\int_{\mathbb{R}}d\omega (\lambda) \, \lambda
(1+\lambda^2)^{-1}>0. \lb{2.31}
\end{align}

\noindent (iii). \begin{align}
& \xi (\lambda)=0 \text{ for } \lambda <0 \text{ if and
only if }
\no \\
& d=0, \quad [0,\infty) \text{ is a support for } \omega \,\,
( \text{i.e., }
\omega ((-\infty,0))=0), \\
& \int^\infty_0 d\omega (1+\lambda)^{-1}<\infty, \text{ and }
c \geq\int^\infty_0d\omega (\lambda)\lambda (1+\lambda^2)^{-1}.
\no
\end{align}
In this case
\begin{equation}\lb{2.33}
\lim_{\lambda\downarrow
-\infty}m(\lambda)=c-\int^\infty_0d\omega (\lambda
')\lambda'(1+{\lambda'}^{2})^{-1}
\end{equation}
and
\begin{equation}\lb{2.34}
c>\int^\infty_0d\omega(\lambda )\lambda
(1+\lambda ^2)^{-1} \text{ if and only if }
\int^\infty_0d\lambda\,
\xi (\lambda )(1+\lambda)^{-1}<\infty.
\end{equation}

\noindent (iv). Let $(\lambda_1,\lambda_2)\subset
\mathbb{R}$
and suppose
$0\leq A\leq \xi (\lambda)\leq B\leq 1$ for a.e.
$\lambda \in (\lambda_1,\lambda_2)$ with $(B-A)<1$. Then
$\omega$ is
purely absolutely continuous in
$(\lambda_1,\lambda_2)$ and $\frac{d\omega}{d\lambda}\in
L^p((\lambda_3,\lambda_4);d\lambda)$ for $[\lambda _3,
\lambda_4]\subset
(\lambda_1,\lambda_2)$ and all $p<(B-A)^{-1}$.

\noindent (v). The measure $\omega$ is purely singular,
$\omega=\omega_s$,
$\omega_{ac}=0$ if and only if $\xi$ equals the
characteristic
function
of a measurable subset $A\subseteq\mathbb{R}$, that is,
$\xi=\chi_A$.
\end{theorem}

As mentioned after Definition~\ref{d2.1}, the definition of
$m\big|_{\bbC_-}$ by means of reflection as in \eqref{2.1},
in general, does not
represent the analytic continuation of $m\big|_{\bbC_+}$.
The
following
result of Greenstein \cite{Gr60} clarifies those
circumstances
under
which $m$ can be analytically
continued from $\bbC_+$ into a subset of $\bbC_-$
through an
interval
$(\lambda_1,\lambda_2)\subseteq \bbR$.

\begin{lemma}  [\cite{Gr60}] \lb{l2.5}
Let $m$ be a Herglotz function with representation \eqref{2.14}
and
$(\lambda_1,\lambda_2)\subseteq\bbR$, $\lambda_1<\lambda_2$.
Then $m$ can
be analytically continued from
$\bbC_+$ into a subset of $\bbC_-$ through the interval
$(\lambda_1,\lambda_2)$
if and only if the
associated measure $\omega$ is purely absolutely continuous on
$(\lambda_1,\lambda_2)$, $\omega {\big|_{(\lambda_1,\lambda_2)}}
=\omega {\big|_{(\lambda_1,\lambda_2),ac}}$, and the density
$\omega'$ of $\omega$ is real-analytic on $(\lambda_1,
\lambda_2)$.
In this case, the
analytic continuation of $m$ into some domain
$\calD_-\subseteq\bbC_-$
is given by
\begin{equation}
m(z)=\overline{m(\overline z)} + 2\pi i \omega'(z),
\quad z\in\calD_-, \lb{2.36}
\end{equation}
where $\omega'(z)$ denotes the complex-analytic extension
of
$\omega'(\lambda)$ for $\lambda\in(\lambda_1,\lambda_2)$.
In
particular, $m$ can be analytically continued through
$(\lambda_1,\lambda_2)$ by reflection, that is,
$m(z)=\overline{m(\overline z)}$ for all $z\in\bbC_-$ if
and only
if $\omega$ has no support in $(\lambda_1,\lambda_2)$.
\end{lemma}

If $m$ can be analytically continued through $(\lambda_1,
\lambda_2)$
into some region $\calD_-\subseteq\bbC_-$, then
$\widetilde m(z):=m(z)-\pi i \omega'(z)$ is
real-analytic on
$(\lambda_1,\lambda_2)$ and hence can be continued
through
$(\lambda_1,\lambda_2)$ by reflection. Similarly,
$\omega'(z)$,
being real-analytic, can be continued through $(\lambda_1,
\lambda_2)$
by reflection. Hence \eqref{2.36} follows from
\begin{equation}
m(z)-\pi i \omega'(z)=\widetilde m(z)
=\overline{{\widetilde m}(\overline z)}
=\overline{m(\overline z)} + \pi i \omega'(z),
\quad z\in\calD_-. \lb{2.37}
\end{equation}

Formula \eqref{2.36} shows that any possible singularity
behavior
of $m\big|_{\bbC_-}$ is determined by that of
$\omega'\big|_{\bbC_-}$. (Note that
$m$, being Herglotz, has no singularities in $\bbC_+$.)
Moreover,
analytic continuations through different intervals on
$\bbR$ may
lead to different $\omega'(z)$ and hence to branch cuts
of $m\big|_{\bbC_-}$.

The following result of Kotani \cite{Ko84}, \cite{Ko87}
is
fundamental in connection with applications of Herglotz
functions
to reflectionless Schr\"odinger and Dirac-type operators
on $\bbR$
(i.e., solitonic, periodic, and certain classes of
quasi-periodic
and almost-periodic operators \cite{Cr89}, \cite{GHL98},
\cite{GKT96},
\cite{GS96b}, \cite{Ko84}, \cite{Ko87}, \cite{KK88},
\cite{Le87}).

\begin{lemma}  [\cite{Ko84}, \cite{Ko87}] \lb{l2.6}
Let $m$ be a Herglotz function and $(\lambda_1,\lambda_2)
\subseteq\bbR$,
$\lambda_1<\lambda_2$. Suppose $\lim_{\varepsilon\to 0}
\Re(m(\lambda+i\varepsilon))=0$ for a.e.
$\lambda\in(\lambda_1,
\lambda_2)$.
Then $m$ can be analytically continued from $\bbC_+$ into
$\bbC_-$ through the interval $(\lambda_1,\lambda_2)$ and
\begin{equation}
m(z)=-\overline{m(\overline z)}. \lb{2.38}
\end{equation}
In addition, $\Im(m(\lambda+i0))>0$, $\Re(m(\lambda+i0))=0$
for all $\lambda\in(\lambda_1,\lambda_2)$.
\end{lemma}

\section{Support Theorems in the Scalar Case} \lb{s3}

This section further reviews the case of scalar Herglotz
functions
and
focuses on support theorems for $\omega$,
$\omega_{ac}$, $\omega_s$, etc., in \eqref{2.14}. In
addition,
we recall
the main results of the Aronszajn-Donoghue theory
relating $m_a(z)$,
$a\in \mathcal{A}_2$ (cf.~\eqref{1.7}) and $m(z)$ as in
\eqref{2.11}.

Let $\mu,\nu$ be Borel measures on $\mathbb{R}$. We
recall
that $S_\mu$
is called a support of $\mu$ if $\mu (\mathbb{R}\backslash
S_\mu)=0$. The topological support $S^{cl}_\mu$ of
$\mu $ is
then the
smallest closed support of
$\mu$. In addition, a support $S_\mu$ of $\mu $ is
called minimal
relative to $\nu$ if for any smaller support
$T_\mu\subseteq
S_\mu$ of $\mu$, $\nu(S_\mu\backslash T_\mu)=0$
(or equivalently,
$\widetilde{T}\subset S_\mu$ with
$\mu (\widetilde{T})=0$ implies $\nu (\widetilde{T})=0$).
Minimal
supports are
unique up to sets of $\mu$ and $\nu$ measure zero and
\begin{equation}S\sim T\quad\text{ if and only if}
\quad\mu
(S\vartriangle
T)=0=\nu (S\vartriangle T)\lb{3.1}\end{equation}
defines an
equivalence class $\mathcal{E}_\nu (\mu)$ of minimal
supports
of $\mu$
relative to $\nu$ (with
$S\vartriangle T=(S\backslash T)\cup (T\backslash S)$
the symmetric
difference of $S$ and $T$).

Two measures, $\mu$ and $\nu$, are called orthogonal,
$\mu \bot
\nu$, if some
of their supports are disjoint.

If $\mu_1,\mu_2$ are absolutely continuous
with respect to $\nu$, $\mu_j\ll \nu$, $j=1,2$, and $\mu_1$ and
$\mu_2$
have a common support minimal relative to $\nu$, then $\mu_1$
and $\mu_2$
are equivalent, $\mu_1\sim\mu_2$.

From now on the reference measure $\nu$ will be chosen to be
Lebesgue
measure on $\mathbb{R}$, ``minimal'' without further
qualifications will
always refer to minimal relative to Lebesgue measure on
$\mathbb{R}$, and
the corresponding equivalence class $\mathcal{E}_\nu (\mu)$
will simply
be denoted by $\mathcal{E}(\mu)$.

For pure point measures $\mu=\mu_{pp}$ we agree to
consider
only the smallest
support (i.e., the countable set of points with
positive $\mu_{pp}$
measure). If the support of a pure point measure
$\mu=\mu_{pp}$
contains no finite accumulation points we call it a
discrete point
measure and denote it by $\mu_d$.

It can be shown that there always exists a minimal support
$S_\mu$ of
$\mu$ such that $\overline{S_\mu}=S^{cl}_\mu$
(cf.~\cite{GP87},
Lemma 5), but in general, a minimal support and the
corresponding topological support $S^{cl}_\mu$ may differ
by a set of
positive Lebesgue measure. Frequently, minimal supports
are called
essential supports in the literature.

\begin{theorem} [\cite{Ar57}, \cite{GS96b}, \cite{Gi89},
\cite{GP87}] \label{t3.1}
Let $m$ be a Herglotz function with representations
\eqref{2.14} and \eqref{2.29}. Then

\noindent (i). \begin{align}
S_{\omega_{ac}}=\{\lambda \in \mathbb{R} \,
|\lim_{\varepsilon
\downarrow 0}\Im (m(\lambda +i\varepsilon )) \text{ exists
finitely
and } 0<\Im (m(\lambda+i0))<\infty\} \lb{3.2}
\end{align}
is a minimal support of $\omega_{ac}$.

\noindent (ii). \begin{equation}
S_{\omega_s}=\{\lambda \in
\mathbb{R} \, |\lim_{\varepsilon \downarrow 0}
\Im (m(\lambda+i\varepsilon
))=+\infty\}\lb{3.3}
\end{equation}
and
\begin{equation}
S_{\omega_{sc}}=\{\lambda \in
S_{\omega_s}|\lim_{\varepsilon \downarrow 0}\varepsilon
 \Im (m(\lambda
+i\varepsilon ))=0\}\lb{3.4}
\end{equation}
are minimal supports of
$\omega_s$ and $\omega_{sc}$, respectively.

\noindent (iii). \begin{equation}\lb{3.5}
S_{\omega_{pp}}=\{\lambda \in
\mathbb{R} \, | \lim_{\varepsilon \downarrow 0}\varepsilon
\Im (m(\lambda+i\varepsilon ))=
-i\lim_{\varepsilon\downarrow 0}
\varepsilon m(\lambda+i\varepsilon) >0\}
\end{equation}
is the smallest support of $\omega_{pp}$.

\noindent (iv). $S_{\omega_{ac}}$, $S_{\omega_{sc}}$, and
$S_{\omega_{pp}}$
are mutually disjoint minimal supports and
\begin{align}
S_\omega &=\{\lambda \in
\mathbb{R} \, |\lim_{\varepsilon
\downarrow 0}\Im (m(\lambda+i\varepsilon ))\leq
+\infty \text{ exists and } 0<\Im (m(\lambda+i0))\leq
+\infty\} \no \\
&=S_{\omega_{ac}} \cup S_{\omega_{s}} \lb{3.6}
\end{align}
is a minimal support for $\omega$.

\noindent (v). \begin{equation}\lb{3.7}
\widetilde{S}_{\omega
_{ac}}=\{\lambda \in \mathbb{R} \, | \, 0<\xi (\lambda )<1\}
\end{equation}
is a
minimal support for $\omega_{ac}$.
\end{theorem}

Of course
\begin{align} \lb{3.8}
\widehat{S}_{\omega_{ac}}=\{\lambda \in \mathbb{R}
\, | \lim_{\varepsilon \downarrow
0}m(\lambda +i\varepsilon) \text{ exists finitely
and } 0<\Im (m(\lambda+i0))<\infty\}
\end{align}
is also a minimal support of
$\omega_{ac}$. Later on we shall use the analog of
\eqref{3.8}
in the matrix-valued context (cf.~Section~\ref{s6}).

The equivalence relation \eqref{3.1} motivates the
introduction of
equivalence classes associated with $\omega$ and its
decompositions
$\omega_{ac}$, $\omega_s$, etc. We will, in particular use
the following
two equivalence classes
\begin{equation}\lb{3.9}
\mathcal{E}(\omega_{\underset{s}{ac}}) := \text{the
equivalence class of minimal supports of }
\omega_{\underset{s}{ac}}
\end{equation}
in Theorem~\ref{t3.2} below.

Next, we turn our attention to $m_a(z)$ in \eqref{2.11}, with
$a\in
M_2(\mathbb{C})$ satisfying \eqref{2.6}. We abbreviate the
identity matrix
in $M_n(\bbC)$ by $I_n$, the unit circle in $\bbC$ by
$S^1=\partial D$,
and introduce the set (cf.~\eqref{2.6}),
\begin{equation} \lb{3.10}
\mathcal{A}_2=\{a\in M_2(\mathbb{C})| \, a^*j_2a=j_2\}.
\end{equation}
We note that
\begin{equation} \lb{3.10b}
|\det(a)|=1, \quad a\in\calA_2,
\end{equation}
and
\begin{equation}
(a_{1,1}/a_{1,2}), (a_{1,1}/a_{2,1}), (a_{2,2}/a_{1,2}),
(a_{2,2}/a_{2,1})
\in\bbR, \quad a\in\calA_2 \lb{2.10c}
\end{equation}
as long as $a_{1,2}\neq 0$, respectively, $a_{2,1}\neq 0$.
Moreover, we
recall (cf. \eqref{2.11}),
\begin{equation}
m_a(z)=\frac{a_{2,1}+a_{2,2}m(z)}{a_{1,1}+a_{1,2}m(z)},
\quad z\in\bbC_+ \lb{3.10a}
\end{equation}
and its general version
\begin{equation}\lb{3.17}
m_a(z)=\frac{(ab^{-1})_{2,1}+(ab^{-1})_{2,2}m_b(z)}
{(ab^{-1})_{1,1}+
(ab^{-1})_{1,2}m_b(z)}, \quad a,b \in\calA_2, \,\,
z\in\bbC_+.
\end{equation}

The corresponding
equivalence classes of minimal supports of
$\omega_{\underset{s}
{ac}}$ and
$\omega_{a,{\underset{s}{ac}}}$ are then denoted by
$\mathcal{E}(\omega_{\underset{s}{ac}})$ and
$\mathcal{E}(\omega_{a,{\underset{s}{ac}}})$.

The celebrated
Aronszajn-Donoghue theory then revolves around the
following
result.

\begin{theorem} [\cite{Ar57}, \cite{Do65}, see also
\cite{GP87}, \cite{SW86}] \label{t3.2}
Let $m(z)$ and $m_a(z)$,
$a\in \mathcal{A}_2$ be Herglotz functions related by
\eqref{3.10a},
with corresponding measures $\omega$ and $\omega_a$,
respectively. Then

\noindent (i). For all $a\in \mathcal{A}_2$,
\begin{equation}\lb{3.11}
\mathcal{E}(\omega_{a,ac})=\mathcal{E}(\omega_{ac}),
\end{equation}
that is, $\mathcal{E}(\omega_{a,ac})$ is independent of
$a\in \mathcal{A}_2$
(and hence denoted by $\mathcal{E}_{ac}$ below) and
$\omega_{a,ac}\sim
\omega_{ac}$ for all $a\in\calA_2$.

\noindent (ii). Suppose $\omega_b$ is a discrete point
measure,
$\omega_b=\omega_{b,d}$, for some $b\in \mathcal{A}_2$.
Then
$\omega_a=\omega_{a,d}$ is a discrete point measure for
all $a\in
\mathcal{A}_2$.

\noindent (iii). Define
\begin{equation}
S=\{\lambda \in \mathbb{R} \, | \,
\text{there is no } a\in \mathcal{A}_2 \text{ for which }
\Im (m_a(\lambda+i0)) \text{ exists and equals }
0\}.\lb{3.12}
\end{equation}
Then $S\in \mathcal{E}_{ac}$.

\noindent (iv). Suppose $a_{1,2}\neq 0$ (i.e., $a\in
\mathcal{A}_2\backslash \{\gamma I_2\}$, $\gamma\in S^1$). If
$\omega_{a,s} (\mathbb{R})>0$ or
$\omega_s(\mathbb{R})>0$, then
$\mathcal{E}(\omega_{a,s})\neq \mathcal{E}(\omega_s)$ and
there exist
$S_{a,s}\in\calE(\omega_{a,s})$, $S_s\in\calE(\omega_s)$
such that
$S_{a,s}\cap S_s =\emptyset$ (i.e., $\omega_{a,s}\bot\,
\omega_s$).
In particular,
\begin{equation} \lb{3.13}
\widetilde{S}_{\omega_{a,s}}=\{\lambda \in
\mathbb{R} \, |\lim_{\varepsilon \downarrow 0}m(\lambda
+i\varepsilon)=-a_{1,1}/a_{1,2}\}
\end{equation}
is a minimal support for
$\omega_{a,s}$ and the smallest support of $\omega_{a,pp}$
equals
\begin{equation}\lb{3.14}
S_{\omega_{a,pp}}=\{\lambda
\in \mathbb{R} \, |\lim_{\varepsilon \downarrow 0}
m(\lambda+i\varepsilon
)=-a_{1,1}/a_{1,2},\, \int_{\mathbb{R}}d\omega
(\lambda')(\lambda'-\lambda)^{-2}<\infty\}.
\end{equation}
Moreover,
\begin{equation}\lb{3.15}
\omega_a(\{\lambda
\})=|a_{1,2})|^{-2}\left(d+\int_{\mathbb{R}}d\omega
(\lambda')(\lambda'-\lambda)^{-2}\right)^{-1},\quad \lambda \in
\mathbb{R}.
\end{equation}

\noindent (v). Suppose $\omega_b$ is a discrete point measure
for some
(and hence for all) $b\in\calA_2$.
Assume that $\supp(\omega)=\{\lambda_{I_2,n}\}_{n\in\calI}$ and
$\supp(\omega_a)=\{\lambda_{a,n}\}_{n\in\calI}$ for some
$a\in\calA_2$
with $a_{1,2}\neq 0$ are given, where $\calI$ is either
$\bbN$, $\bbZ$, or
a finite non-empty index set. Suppose in addition that one
of the following
conditions hold: (1) $\omega(\bbR)$ is known, or (2) $m(z_0)$
is known for
some $z_0\in\bbC_+$, or (3) $\lim_{z\to
i\infty}(m(z)-m^0 (z))=0$, where
$m^0 (z)$ is a known Herglotz function. Then the system of
measures
$\{\omega_b\}_{b\in\calA_2}$ and hence the system of Herglotz
functions
$\{m_b(z)\}_{b\in\calA_2}$ is uniquely determined.
\end{theorem}

\noindent{\it Sketch of Proof.} (i), (iii), and (iv) follow
from \eqref{3.10} and \eqref{3.10a} which imply
\begin{equation}
\Im (m_a(z))=\frac{\Im
(m(z))}{|a_{1,1}+a_{1,2}m(z)|^2},\lb{3.16}
\end{equation}
from
Theorem~\ref{t2.2}\,(i), (ii), and from Theorem~\ref{t3.1}.
Note that
$a_{1,1}=a_{1,2}=0$ cannot occur in \eqref{3.16} since
this would
contradict \eqref{3.10}. (ii) follows from \eqref{3.17}
and the fact that
$\omega_a=\omega_{a,d}$ if and only if $m_a(z)$ is
meromorphic on
$\mathbb{C}$. In order to prove (v) we define
\begin{equation}
F(z)=m(z)+\frac{a_{1,1}}{a_{1,2}}, \quad z\in\bbC_+.
\lb{3.17a}
\end{equation}
Then $F$ is a meromorphic Herglotz function with simple
zeros at
$\{\lambda_{a,n}\}_{n\in\calI}$ and simple poles at
$\{\lambda_{I_2,n}\}_{n\in\calI}$. In particular, its
zeros and
poles
necessarily interlace and the exponential Herglotz
representation
\eqref{2.29} for $F$ then yields
\begin{equation}
F(z)=\exp\big ( k+\int_{\mathbb{R}}d\lambda \, \xi(\lambda)
((\lambda-z)^{-1}-\lambda(1+\lambda^2)^{-1}) \big ),
\lb{3.17b}
\end{equation}
with $\xi$ a piecewise constant function. Analyzing
\eqref{2.30}
shows
that
\begin{equation}
\xi(\lambda)=\chi_{\{\lambda\in\bbR\,|\,F(\lambda)<0\}}
(\lambda),
\lb{3.17c}
\end{equation}
where $\chi_{\calM}$ denotes the characteristic function
of a set
$\calM\subseteq\bbR$ and hence $\xi$ is uniquely determined
by $\supp(\omega)$
and $\supp(\omega_a)$. Thus $F(z)$ is uniquely determined
except for
the constant $k\in\bbR$ (which cannot be determined from
$\supp(\omega)$
and $\supp(\omega_a)$). Either one of the conditions (1)--(3)
then will
determine $k$ and hence $F(z)$, $z\in\bbC_+$. Thus $m(z)$,
and hence by
\eqref{3.10a} $m_b(z)$ for all $b\in\calA_2$, are uniquely
determined,
which in turn determine $\omega_b$ for all $b\in\calA_2$.
\hfill$\square$ \\

For connections between Theorem~\ref{t3.2} \,(iv) and Hankel
operators see \cite{MPT95}, Sect. III.10.

The relationship between $\Im (\ln (m_a(z)))$ (respectively,
$\xi _a(\lambda ))$ and
$\Im (\ln (m(z)))$ (respectively, $\xi (\lambda ))$, analogous
to \eqref{3.10a},
in general, is quite involved. The special case $a=j_2$,
that is,
\begin{equation}\lb{3.18}
m_{j_2}(z)=-1/m(z),
\end{equation}
however,
is particularly simple and leads to
\begin{equation} \lb{3.19}
\xi_{j_2}(\lambda )=1-\xi (\lambda )
\text{ for a.e. } \lambda \in\mathbb{R}.
\end{equation}
We also state the following elementary result.

\begin{lemma} \lb{l3.3}
Suppose $a,b\in\calA_2$ and that $m_a(z)$ is a nonconstant
Herglotz
function. Then $m_a(z)=m_b(z)$ for all $z\in\bbC_+$ if and
only if
$a=\gamma b$ for some $\gamma\in S^1$.
\end{lemma}

\begin{proof}
First we note that \eqref{2.5} and \eqref{2.6} determine a
subgroup of the group of M\"obius transformations
(characterized
by leaving $\bbC_{\pm}$ invariant and normalized by
$|\det (a)|=1$).
Hence $m_a(z)=m_b(z)$ if and only if $a=\gamma b$ for some
$\gamma \in
\bbC\backslash\{0\}$. The normalization
$|\det (a)|=|\det (b)|=1$ then
yields $\gamma\in S^1$.
\end{proof}

\section{Further Applications of Scalar Herglotz
Functions} \lb{s4}

For additional applications of the Aronszajn-Donoghue theory
described in
Theorem~\ref{t3.2} we now consider self-adjoint rank-one
perturbations
of self-adjoint operators, Friedrichs and Krein extensions of
densely defined symmetric operators bounded
from below with deficiency indices $(1,1)$, and Sturm-Liouville
operators on a half-line.

Some of the following results are well-known (the material
mainly being taken from \cite{DM91}, \cite{DM95}, \cite{DMT88},
\cite{Do65}, \cite{Kr47a}, \cite{Si95b}, \cite{SW86}, and
\cite{Sk79}) and hence this section is partially expository
in character. However, we do supply
simplified proofs of various
results below and prove several realization results for
different classes of Herglotz functions which appear to be
new. Additional material connecting rank-one perturbations
with Hankel operators, respectively, Krein's spectral shift
function can be found in \cite{MPT95}, Sect.~III.10,
respectively, \cite{Po97}.\\

We start with self-adjoint rank-one perturbations of
self-adjoint
operators (following \cite{Do65} and \cite{SW86}). Let
$\calH$ be a
separable complex Hilbert space with scalar product
$(\cdot,\cdot)_{\calH}$, $H_0$ a self-adjoint operator in
$\calH$ (which may or may not be bounded) with simple spectrum.
Suppose
$f_1\in\calH$, $\|f_1\|_{\calH}=1$ is a cyclic vector for
$H_0$ (i.e.,
$\calH=\overline{\text{linspan}\{(H_0-z)^{-1}f_1\in\calH\,|\,
z\in\bbC\backslash\bbR\}}$, or equivalently,
$\calH$$=\overline{\text{linspan}
\{E_0(\lambda)f_1\in\calH\,| }$ $\overline{ \lambda\in\bbR\}}$,
$E_0(\cdot)$ the
family of orthogonal spectral projections of $H_0$) and define
\begin{equation}
H_{\alpha}=H_0 +\alpha P_1, \quad P_1=(f_1,\cdot)_{\calH}f_1,
\quad \alpha\in\bbR, \lb{4.1}
\end{equation}
with $\calD (H_{\alpha})=\calD (H_0)$, $\alpha\in\bbR$
($\calD (\cdot)$ abbreviating
the domain of a linear operator). Denote by $E_{\alpha}(\cdot)$
the family
of orthogonal spectral projections of $H_{\alpha}$ and define
\begin{equation}
d\omega_{\alpha}(\lambda)=d\|E_{\alpha}(\lambda)f_1\|_{\calH}^2,
\quad
\int_{\bbR} d\omega_{\alpha}(\lambda)=\|f_1\|_{\calH}^2=1.
\lb{4.2}
\end{equation}
By the spectral theorem for self-adjoint operators (cf., e.g.,
\cite{Na68}, Ch.~VI), $H_{\alpha}$ in
$\calH$ is unitarily equivalent to ${\hatt H}_{\alpha}$ in
$\hatt \calH_\alpha=L^2(\bbR;d\omega_{\alpha})$, where
\begin{align}
&({\hatt H}_{\alpha} \hat g) (\lambda)=\lambda \hat g(\lambda),
\quad
\hat g\in\calD ({\hatt H}_{\alpha})=
L^2(\bbR;(1+\lambda^2)d\omega_{\alpha}), \lb{4.3} \\
&H_{\alpha}=U_{\alpha} {\hatt H}_{\alpha} U_{\alpha}^{-1},
\quad
\calH=U_{\alpha} L^2(\bbR;d\omega_{\alpha}), \lb{4.4}
\end{align}
with $U_{\alpha}$ unitary,
\begin{equation}
U_\alpha:\widehat \calH_\alpha=L^2(\bbR;d\omega_\alpha)
\to\calH,
\quad
\hat g\to(U_\alpha \hat g)=\slim_{N\to\infty}
\int_{-N}^N
d(E_\alpha(\lambda)f_1)\hat g(\lambda). \lb{4.4a}
\end{equation}
Moreover,
\begin{equation}
f_1=U_{\alpha}\hat f_1, \quad \hat f_1 (\lambda)=1,
\quad \lambda\in\bbR.
\lb{4.5}
\end{equation}
The family of spectral projections
$\widehat E_\alpha(\lambda)$,
$\lambda\in\bbR$ of $\widehat H_\alpha$ is then given by
\begin{align}
&(\widehat E_\alpha (\lambda) \hat g)(\mu)=\theta(\lambda -
\mu)
\hat g(\mu) \text{ for $\omega_\alpha$--a.e. } \mu\in\bbR,
\,\, \hat g\in
L^2(\bbR;d\omega_\alpha), \\ \lb{4.5a}
& \hspace*{6.7cm} \theta(x)=\begin{cases} 1,&x\geq 0 \\
0,&x < 0. \end{cases} \no
\end{align}
Introducing the Herglotz function
\begin{equation}
m_{\alpha}(z)=(f_1,(H_{\alpha}-z)^{-1}f_1)_{\calH}=\int_{\bbR}
\frac{d\omega_{\alpha}}{\lambda-z}, \quad z\in\bbC_+, \lb{4.6}
\end{equation}
one verifies
\begin{equation}
m_{\beta}(z)=\frac{m_{\alpha}(z)}{1+(\beta-\alpha) m_{\alpha}(z)},
\quad \alpha, \beta \in\bbR. \lb{4.7}
\end{equation}
A comparison of \eqref{4.7} and \eqref{3.10a} suggests an
introduction of
\begin{equation}
a(\alpha,\beta)=\left(\begin{array}{cc} 1 & \beta-\alpha
\\ 0 &
1 \end{array}\right) \in\calA_2, \quad \alpha,\beta\in\bbR.
\lb{4.8}
\end{equation}
Moreover, since $\omega_{\alpha}(\bbR)=1$, Theorem~\ref{t3.2}
applies (with
$a_{1,1}(\alpha,\beta)=a_{2,2}(\alpha,\beta)=1$,
$a_{1,2}(\alpha,\beta)
=\beta-\alpha$, $a_{2,1}(\alpha,\beta)=0$).

If $f_1$ is not a cyclic vector for $H_0$, then as discussed
in \cite{Do65},
$\calH$ (not necessarily assumed to be separable at
this point)
decomposes into two orthogonal subspaces $\calH^1$ and
$\calH^{1,\bot}$,
\begin{equation}
\calH=\calH^1\oplus\calH^{1,\bot}, \lb{4.8aa}
\end{equation}
with $\calH^1$ separable,
each of which is a reducing subspace for all $H_{\alpha}$,
$\alpha\in\bbR$. One then has $\calH^1=
\overline{\text{linspan }
\{(H_0 -z)^{-1}f_1\in\calH\,|\, z\in\bbC\backslash\bbR \}}$
and
\begin{equation}
H_{\alpha}=H_{\beta} \text{ on } \calD(H_0)\cap\calH^{1,\bot}
\text{ for all } \alpha, \beta\in\bbR. \lb{4.8a}
\end{equation}
In particular,
\begin{align}
& H_0=H_0^1\oplus H_0^{1,\bot}, \quad
H_{\alpha}=H_{\alpha}^1\oplus H_0^{1,\bot}, \quad
\alpha\in\bbR,
\lb{4.8b} \\
& f_1=f_1^1\oplus 0, \lb{4.8c}
\end{align}
where
\begin{equation}
H_0^1=H \big |_{\calD(H_0)\cap\calH^1}, \quad H_0^{1,\bot}=
H_0\big |_{\calD(H_0)\cap\calH^{1,\bot}}, \lb{4.28d}
\end{equation}
implying
\begin{equation}
(f_1,(H_{\alpha} -z)^{-1}f_1)_{\calH}=
(f_1^1,(H_{\alpha}^1 -z)^{-1}f_1^1)_{\calH^1}=m_{\alpha}^1 (z),
\quad \alpha\in\bbR. \lb{4.8e}
\end{equation}
Thus, $\alpha$-dependent spectral properties of $H_{\alpha}$ in
$\calH$ are effectively
reduced to those of $H_{\alpha}^1$ in $\calH^1$, where
$H_{\alpha}^1$
are self-adjoint operators with simple spectra and cyclic vector
$f_1^1\in\calH^1$. \\

Introducing the following set of Herglotz functions (we will
choose the usual
symbol $\calN$ for these sets in honor of R.~Nevanlinna)
\begin{equation}
\calN_1 =\{m:\bbC_+\to\bbC_+ \text{ analytic }|\,m(z)=
\smallint_{\bbR}
d\omega(\lambda)(\lambda -z)^{-1}, \,\,\,
 \smallint_{\bbR} d\omega(\lambda)<\infty \}, \lb{4.9a}
\end{equation}
we now turn to a realization theorem for Herglotz functions
of the type
\eqref{4.6}.

\begin{theorem} \mbox{\rm } \lb{t4.1a}

\noindent (i). Any $m\in\calN_1$ with associated measure
$\omega$ can be
realized in the form
\begin{equation}
m(z)=(f_1,(H-z)^{-1}f_1)_{\calH}, \quad z\in\bbC_+, \lb{4.9b}
\end{equation}
where $H$ denotes a self-adjoint operator in some separable
complex
Hilbert space $\calH$, $f_1\in\calH$, and
\begin{equation}
\int_{\bbR}d\omega(\lambda)=\|f_1\|^2_{\calH}. \lb{4.9c}
\end{equation}

\noindent (ii). Suppose $m_\ell\in\calN_1$ with corresponding
measures
$\omega_\ell$, $\ell=1,2$, and $m_1\neq m_2$. Then $m_1$ and
$m_2$ can be
realized as
\begin{equation}
m_\ell(z)=(f_1,(H_\ell -z)^{-1}f_1)_{\calH}, \quad \ell=1,2,
\,\, z\in\bbC_+, \lb{4.9d}
\end{equation}
where $H_\ell$, $\ell=1,2$ are self-adjoint rank-one
perturbations
of one and
the same self-adjoint operator $H_0$ in some complex Hilbert
space
$\calH$ (which may be chosen separable) with $f_1\in\calH$,
that is,
\begin{equation}
H_\ell=H_0+\alpha_\ell P_1, \quad P_1=(f_1,\cdot)_{\calH}f_1
\lb{4.9e}
\end{equation}
for some $\alpha_\ell\in\bbR$, $\ell=1,2$, if and only if the
following
conditions hold:
\begin{equation}
\int_{\bbR}d\omega_1 (\lambda)=\int_{\bbR}d\omega_2 (\lambda)
=\|f_1\|^2_{\calH}, \lb{4.9f}
\end{equation}
and for all $z\in\bbC_+$,
\begin{equation}
m_2(z)=\frac{m_1(z)}{1+
\|f_1\|^{-2}_{\calH}(\alpha_2 -\alpha_1)m_1(z)}. \lb{4.9g}
\end{equation}
\end{theorem}

\begin{proof}
Define the self-adjoint operator $H_0$ of multiplication by
$\lambda$ in
$\calH=L^2(\bbR;d\omega)$ by
\begin{equation}
(H_0g)(\lambda)=\lambda g(\lambda), \quad g\in\calD(H_0)=
L^2(\bbR;(1+\lambda^2)d\omega), \lb{4.9h}
\end{equation}
where $\omega$ denotes the measure in the Herglotz
representation
of
$m(z)$, and consider $f_1=1\in\calH$. One infers
\begin{equation}
(f_1,(H_0 -z)^{-1}f_1)_{\calH}=\int_{\bbR}d\omega(\lambda)
(\lambda-z)^{-1}
=m(z). \lb{4.9i}
\end{equation}
Since $\wlim_{z\to i\infty} (-z)(H_0-z)^{-1}=I_{\calH}$,
the identity in $\calH$, and $|iy(\lambda-iy)^{-1}|\leq 1$,
Lebesgue's
dominated convergence theorem yields \eqref{4.9c} and hence
part (i).
The necessity of condition \eqref{4.9g} in part (ii) was
proved by
Donoghue \cite{Do65} (assuming $\|f_1\|_{\calH}=1$). Indeed,
applying the last part in the argument proving (i) to
$m_1(z)$
and
$m_2(z)$ immediately proves \eqref{4.9f}. Identifying
$\alpha_1=
\alpha$, $\alpha_2=\beta$, $H_1=H_\alpha$, $H_2=H_\beta$,
$\|f_1\|^{-2}_{\calH}m_1(z)=m_\alpha (z)$, and
$\|f_1\|^{-2}_{\calH}m_2(z)=m_\beta (z)$, \eqref{4.9g} is
easily seen
to be equivalent to \eqref{4.7}. Conversely, assume
\eqref{4.9f}
and
\eqref{4.9g}. By part (i), we may realize $m_1(z)$ as
\begin{equation}
\|f_1\|^{-2}_{\calH}m_1(z)=\|f_1\|^{-2}_{\calH}
(f_1,(H_1 -z)^{-1}f_1)_{\calH}. \lb{4.9j}
\end{equation}
By \eqref{4.8e} we may assume that $\calH$ is separable and
$H_1$ has
simple spectrum and hence identify it with $H_\alpha$ in
\eqref{4.1}.
Define $H_\beta$ as in \eqref{4.1} for
$\beta\in\bbR\backslash\{\alpha\}$
and consider
\begin{equation}
m_\beta(z)=\|f_1\|^{-2}_{\calH}
(f_1,(H_\beta -z)^{-1}f_1)_{\calH}.
\lb{4.9k}
\end{equation}
By \eqref{4.7} one obtains ($m_\alpha (z)=
\|f_1\|^{-2}_{\calH} m_1(z)$)
\begin{equation}
m_\beta (z)=\frac{\|f_1\|^{-2}_{\calH} m_1(z)}{1+ (\beta
-\alpha)\|f_1\|^{-2}_{\calH} m_1 (z)}. \lb{4.9l}
\end{equation}
A comparison of \eqref{4.9l} and \eqref{4.7} then yields
$\|f_1\|^2_{\calH} m_\beta (z)=m_2(z)$ for
$(\alpha_2 -\alpha_1)=
(\beta -\alpha)$, completing the proof.
\end{proof}

Of course we could have normalized $f_1$, $\|f_1\|_{\calH}=1$,
and then
added the constraint $\smallint_{\bbR} d\omega(\lambda)=1$ to
\eqref{4.9a}. By \eqref{4.8e}, \eqref{4.9b} can be realized in
nonseparable Hilbert spaces. \\

Next we turn to a characterization of Friedrichs and Krein
extensions
of densely defined operators bounded from below with
deficiency indices $(1,1)$ (following \cite{DM91}, \cite{DM95},
\cite{DMT88},
\cite{Do65}, \cite{Kr47a}, \cite{KO78}, \cite{Sk79}, and
\cite{Ts92}).

We start by describing a canonical representation of densely
defined
closed symmetric operators with deficiency indices $(1,1)$ as
discussed
by \cite{Do65}. Let $\calH$ be a separable
complex Hilbert space, $H$ a closed densely defined symmetric
operator
with domain $\calD (H)$ and deficiency indices $(1,1)$. Choose
$u_{\pm}\in\ker(H^{*}\mp i)$ with $\|u_{\pm}\|_{\calH}=1$ and
denote by
$H_{\alpha}$, $\alpha\in [0,\pi)$ all self-adjoint extensions
of
$H$, that is,
\begin{align}
&H_{\alpha}(g+u_++e^{2i\alpha}u_-)=Hg+iu_+ -ie^{2i\alpha}u_-,
\no \\
&\calD (H_{\alpha})=\{(g+ u_+ +e^{2i\alpha}u_-)\in\calD
(H^{*})\,|\, g\in\calD (H), \, u_{\pm}\in\ker(H^*\mp i) \} \lb{4.9}
\end{align}
by von Neumann's formula for self-adjoint extensions of $H$.
Let
$E_{\alpha}(\cdot)$ be the family of spectral projections
of $H_{\alpha}$
and suppose $H_{\alpha}$ has simple spectrum for some (and
hence for all)
$\alpha\in [0,\pi)$ (i.e., $u_+$ is a cyclic vector for
$H_{\alpha}$
for all $\alpha\in [0,\pi)$). Define
\begin{equation}
d\nu_{\alpha}(\lambda)=d\|E_{\alpha}(\lambda)u_+\|_{\calH}^2,
\quad
\int_{\bbR}d\nu_{\alpha}(\lambda)=\|u_+\|_{\calH}^2=1, \quad
\alpha\in [0,\pi), \lb{4.10}
\end{equation}
then $H_{\alpha}$ is unitarily equivalent to multiplication
by $\lambda$
in $L^2(\bbR;d\nu_{\alpha})$ and $u_+$ can be mapped into the
function
identically $1$. However, it is more convenient to define
\begin{equation}
d\omega_{\alpha}(\lambda)=(1+\lambda^2)d\nu_{\alpha}(\lambda),
\lb{4.11}
\end{equation}
such that
\begin{equation}
\int_{\bbR}\frac{d\omega_{\alpha}(\lambda)}{1+\lambda^2}=1,
\quad
\int_{\bbR} d\omega_{\alpha}(\lambda)=\infty, \quad
\alpha\in [0,\pi) \lb{4.12}
\end{equation}
(by \eqref{4.10} and the fact that $u_+\notin\calD (H_{\alpha})$).
Thus,
$H_{\alpha}$ is unitarily equivalent to ${\hatt H}_{\alpha}$ in
$\hatt \calH_\alpha=L^2(\bbR;d\omega_{\alpha})$, where
\begin{align}
 ({\hatt H}_{\alpha}\hat g)(\lambda)&=\lambda \hat g(\lambda),
\quad
 \hat g\in\calD ({\hatt H}_{\alpha})=
L^2(\bbR;(1+\lambda^2)d\omega_{\alpha}), \lb{4.13} \\
& H_{\alpha}=U_{\alpha}{\hatt H}_{\alpha}U_{\alpha}^{-1},
\quad
\calH =U_{\alpha}L^2(\bbR;d\omega_{\alpha}), \lb{4.14}
\end{align}
with $U_{\alpha}$ unitary,
\begin{align}
& U_{\alpha}: \hatt \calH_{\alpha}=L^2(\bbR;d\omega_{\alpha}) \to
\calH, \no \\
& \hat g \to U_{\alpha}\hat g =\slim_{N\to\infty} \int_{-N}^N
 d(E_{\alpha}(\lambda)u_+)(\lambda-i)\hat g(\lambda).
\lb{4.14a}
\end{align}
Moreover,
\begin{equation}
u_{+}=U_{\alpha}\hat u_{+}, \quad \hat u_{+}(\lambda)=
(\lambda - i)^{-1}, \lb{4.15}
\end{equation}
and
\begin{equation}
({\hatt H}(\alpha)\hat g)(\lambda)=\lambda\hat g(\lambda),
\quad
\hat g\in\calD ({\hatt H}(\alpha))=\{\hat h \in\calD
({\hatt H}_{\alpha})\,|\,
\smallint_{\bbR}d\omega_{\alpha}(\lambda)\hat h(\lambda)=0\},
\lb{4.16}
\end{equation}
where
\begin{equation}
H=U_{\alpha}{\hatt H}(\alpha)U_{\alpha}^{-1}. \lb{4.17}
\end{equation}
Thus $\hatt H(\alpha)$ in $L^2(\bbR;d\omega_{\alpha})$ is a
canonical
representation for a densely defined closed symmetric
operator $H$ with deficiency indices $(1,1)$ in a separable
Hilbert space
$\calH$ with cyclic deficiency vector $u_+\in\ker(H^*-i)$.
We shall
prove in Theorem~\ref{t4.1b} below that $\hatt H(\alpha)$ in
$L^2(\bbR;d\omega_{\alpha})$ is actually a model for all
such operators.
Moreover, since
\begin{align}
& ((H-\overline z)g,U_{\alpha}(\cdot -z)^{-1})_{\calH}=
\int_{\bbR}d\omega_{\alpha}(\lambda)
(\lambda -z)\ol{(U_\alpha^{-1} g)(\lambda)}(\lambda -z)^{-1}=0,
\lb{4.17c} \\
& \hspace*{7.6cm} g\in\calD(H), \,\, z\in\bbC\backslash\bbR \no
\end{align}
by \eqref{4.16}, one infers that
$U_{\alpha}(\cdot -z)^{-1}\in\calD(H^{*})$. Since $\calD(H)$
is dense in
$\calH$, one concludes
\begin{equation}
\ker(\hatt H(\alpha)^{*} -z)=\{c(\cdot -z)^{-1}\,|\, c\in\bbC\},
\quad z\in\bbC\backslash\bbR, \lb{4.17a}
\end{equation}
where
\begin{equation}
H^{*}=U_{\alpha}\hatt H (\alpha)^{*}U_{\alpha}^{-1}. \lb{4.17b}
\end{equation}

If $u_+$ is not cyclic for $H_{\alpha}$ then, as shown in
\cite{Do65},
$\calH$ (not necessarily assumed to be separable at this point)
decomposes into two orthogonal subspaces $\calH^0$ and
$\calH^{0,\bot}$,
\begin{equation}
\calH=\calH^0\oplus\calH^{0,\bot}, \lb{4.20aa}
\end{equation}
with $\calH^0$
separable, each of which is
a reducing subspace for all $H_{\alpha}$, $\alpha\in [0,\pi)$
and
\begin{align}
& \calH^0=\overline{\text{linspan}\{(H_\alpha-z)^{-1}
u_{+}\in\calH
\,|\,z\in\bbC\backslash\bbR \}} \text{ is independent of }
\alpha\in [0,\pi), \lb{4.19aa} \\
& (H_{\alpha} -z)^{-1}=(H_{\beta} -z)^{-1} \text{ on }
\calH^{0,\bot}
\text{ for all } \alpha,\beta\in [0,\pi),
\,\, z\in\bbC_+. \lb{4.19a}
\end{align}
In particular, the part $H^{0,\bot}$ of $H$ in
$\calH^{0,\bot}$ is
then self-adjoint,
\begin{align}
& H=H^0\oplus H^{0,\bot}, \quad H_{\alpha}=
H_{\alpha}^0\oplus H^{0,\bot},
\quad \alpha\in [0,\pi), \lb{4.19b} \\
& \ran(H^{0,\bot} -z) = \calH^{0,\bot},
\quad z\in\bbC\backslash\bbR,
\lb{4.19c} \\
& u_+ = u_+^0 \oplus 0, \lb{4.19d}
\end{align}
with $H^0$ a densely defined closed symmetric operator
in $\calH^0$
and deficiency indices $(1,1)$. One then computes
\begin{align}
& z\|u_+\|^2_{\calH} + (1+z^2)(u_+,(H_{\alpha} -z)^{-1}
u_+)_{\calH} \no \\
& =z\|u_+^{0}\|^2_{\calH^0} + (1+z^2)(u_+^0,
(H_{\alpha}^0 -z)^{-1}u_+^0)_{\calH^0}, \quad
\alpha\in [0,\pi)  \lb{4.19e}
\end{align}
and hence $\alpha$-dependent spectral properties of
$H_{\alpha}$
in $\calH$ are effectively
reduced to those of $H_{\alpha}^0$ in $\calH^0$, where
$H_{\alpha}^0$
are self-adjoint operators with simple spectra and cyclic
vector
$u_+^0\in\ker((H^0)^* -i)$.

Next we show the model character of
$(\hatt \calH_\alpha, \hatt H(\alpha), \hatt H_\alpha)$
following the
approach outlined by Donoghue \cite{Do65}.

\begin{theorem} [\cite{Do65}] \lb{t4.1b}
Let $H$ be a densely defined closed symmetric operator with
deficiency
indices $(1,1)$ and normalized deficiency vectors
$u_\pm\in\ker(H^* \mp i)$, $\|u_\pm\|_\calH$ $=1$ in some
separable
complex Hilbert space $\calH$. Let
$H_\alpha$ be a self-adjoint extension of $H$ with simple
spectrum (i.e.,
$u_+$ is a cyclic vector for $H_\alpha$). Then the pair
$(H,H_\alpha)$
in $\calH$ is unitarily equivalent to the pair
$(\hatt H(\alpha), \hatt H_\alpha)$ in $\hatt \calH$
defined in
\eqref{4.16} and \eqref{4.13} with unitary operator
$U_\alpha$ defined in
\eqref{4.14a} (cf. \eqref{4.17} and \eqref{4.14}).
Conversely, given a
measure $d\widetilde \omega$ satisfying
\begin{equation}
\int_{\bbR}\frac{d\widetilde\omega(\lambda)}{1+\lambda^2}=1,
\quad
\int_{\bbR} d\widetilde\omega(\lambda)=\infty, \lb{4.19f}
\end{equation}
define the
self-adjoint operator $\widetilde H$ of multiplication by
$\lambda$
in $\widetilde \calH=L^2(\bbR;d\widetilde \omega)$,
\begin{equation}
(\widetilde H g)(\lambda)=\lambda g(\lambda),
\quad g\in\calD(\widetilde H)
= L^2(\bbR;(1+\lambda^2)d\widetilde \omega), \lb{4.47}
\end{equation}
and the linear operator $H$ in $\widetilde \calH$,
\begin{equation}
\calD(H)=\{g\in\calD(\widetilde H)\,|\, \smallint_{\bbR}
d\widetilde \omega(\lambda)g(\lambda)=0 \}, \quad H=
\widetilde H \big|_{\calD(H)}. \lb{4.51}
\end{equation}
Then $H$ is a densely defined closed symmetric operator in
$\widetilde \calH$ with
deficiency indices $(1,1)$ and deficiency spaces
\begin{equation}
\ker(H^{*} \mp i)=\{c(\lambda \mp i)^{-1} \,|\, c\in\bbC \}.
\lb{4.56}
\end{equation}
\end{theorem}
\begin{proof}
The first part of the theorem (with the exception of the
explicit
expression for the unitary operator $U_\alpha$ in
\eqref{4.14a}) is
due to
Donoghue \cite{Do65} and we essentially sketched the major
steps in
\eqref{4.10}--\eqref{4.19e} above. For the sake of
completeness we
add two more details. First, in connection with proving
the unitary
equivalence stated in \eqref{4.14}, one observes that
$U_\alpha(\hatt H_\alpha -z)^{-1} \hat u_+=(H_\alpha -z)^{-1}
u_+$.
Using the first resolvent identity for $\hatt H_\alpha$ and
$H_\alpha$ then
yields $U_\alpha(\hatt H_\alpha-z)^{-1}
((\hatt H_\alpha -z')^{-1}\hat u_+)
=(H_\alpha-z)^{-1}((H_\alpha-z')^{-1}u_+)$. Since
$z'\in\bbC\backslash\bbR$
is arbitrary, one obtains \eqref{4.14} from the fact that
$u_+$ is
cyclic for $H_\alpha$ . Secondly, in connection with
the domain of $\hatt H(\alpha)$ in \eqref{4.16} one makes
use of the
well-known fact that $\hat h\in\widehat\calH_\alpha$ belongs
to $\calD (\hatt H(\alpha))$ if and only if
$\hat h\in\calD (\hatt H_\alpha)$ and $\hat h$ is orthogonal
to
$\ker({\hatt H}^* -i)$
in the topology of the graph of ${\hatt H}^*$, that is,
\begin{equation}
({\hatt H}^* {\hat h},{\hatt H}^*u_+)_{\hat \calH_\alpha}
+({\hat h},u_+)_{\hat \calH_\alpha}=0 \text{ or }
i(\hatt H_\alpha {\hat h}, u_+)_{\hat \calH_\alpha}+
(\hat h,u_+)_{\hat \calH_\alpha}=0. \lb{4.47aa}
\end{equation}
This is easily seen to be equivalent to $\smallint_\bbR
d\omega_\alpha
(\lambda)\hat h(\lambda)=0$ in \eqref{4.16}.

Since the second part of Theorem~\ref{t4.1b} is stated
but not
explicitly proved in \cite{Do65}, we now sketch such a proof.

Define $\widetilde \calH_{2r}=L^2(\bbR;
(1+\lambda^2)^rd\widetilde \omega)$,
$r\in\bbR$,
$\widetilde \calH_0=\widetilde \calH$ and consider the
isometric
isomorphism (unitary operator)
$R$ from $\widetilde \calH_2$ onto $\widetilde \calH_{-2}$,
\begin{align}
& R:\widetilde \calH_2\to\widetilde \calH_{-2}, \quad
\widetilde f\to (1+\lambda^2)\widetilde f, \lb{4.48} \\
& (\widetilde f,\widetilde g)_{\tilde \calH_2}=
(\widetilde f,R\widetilde g)_{\tilde \calH}=
(R\widetilde f,\widetilde g)_{\tilde \calH}=
(R\widetilde f,R\widetilde g)_{\tilde \calH_{-2}},
\quad \widetilde f,\widetilde g\in\widetilde \calH_2,
\lb{4.49} \\
& (\widetilde u,\widetilde v)_{\tilde \calH_{-2}}=
(\widetilde u,R^{-1}\widetilde v)_{\tilde \calH}=
(R^{-1}\widetilde u,\widetilde v)_{\tilde \calH}=
(R^{-1}\widetilde u,R^{-1}\widetilde v)_{\tilde \calH_2}, \quad
\widetilde u,\widetilde v\in\widetilde \calH_{-2}. \lb{4.50}
\end{align}
We note that $\bbC\subset\widetilde \calH_{-2}$. Since
$\widetilde g\in\widetilde \calH_2$ implies
$\widetilde g\in L^1(\bbR;d \widetilde \omega)$
using $|\widetilde g(\lambda)|=(1+\lambda^2)^{-1/2}
(1+\lambda^2)^{1/2}|\widetilde g(\lambda)|$
and Cauchy's inequality, $\calD(H)$ is well-defined.
Moreover, as a
restriction of the self-adjoint operator $\widetilde H$, $H$
is clearly
symmetric. One infers from
\eqref{4.48} and \eqref{4.51} that
\begin{equation}
\calD(\widetilde H)=\widetilde \calH_2=
\calD(H)\oplus_{\tilde \calH_2}R^{-1}\bbC,
\lb{4.52}
\end{equation}
where, in obvious notation, $\oplus_{\tilde \calH_2}$
denotes the direct
orthogonal sum in $\widetilde \calH_2$. Next, to prove that
$\calD(H)$ is dense
in $\widetilde \calH$, suppose there exists a
$\widetilde g\in\widetilde \calH$
such that $\widetilde g\bot \calD(H)$. Then
\begin{equation}
0=(\widetilde f,\widetilde g)_{\tilde \calH}=
(\widetilde f,R^{-1}\widetilde g)_{\tilde \calH_2}
\text{ for all } \widetilde f\in\calD(H)
\lb{4.53}
\end{equation}
and hence $R^{-1}\widetilde g\in R^{-1}\bbC$, that is,
$\widetilde g=c\in\bbC$ a.e.~by
\eqref{4.52}, and consequently, $\widetilde g\in
\widetilde \calH$ if
and only if
$c=\widetilde g=0$. Next, $H$ is a closed operator, either by
\eqref{4.52} or directly by its definition \eqref{4.51}
($\lim_{n\to\infty}\|\widetilde f_n-\widetilde g\|_{\tilde
\calH}=0$,
$\lim_{n\to\infty}\|H\widetilde f_n-\widetilde g\|_{\tilde
\calH}=0$ for
$\{\widetilde f_n\}_{n\in\bbN}\subset\calD(H)$,
$\widetilde f,\widetilde g\in\widetilde \calH$ imply
$\widetilde f\in\widetilde \calH_2$
and $\widetilde g=\widetilde H\widetilde f$ by passing
to appropriate
subsequences of
$\{\widetilde f_n\}_{n\in\bbN}$ and $\{H\widetilde
f_n\}_{n\in\bbN}$,
and $\smallint_{\bbR}
d \widetilde \omega(\lambda)\widetilde f(\lambda)=
(R^{-1}1,\widetilde f)_{\tilde \calH_2}
=\lim_{n\to\infty}(R^{-1}1,\widetilde f_n)_{\tilde
\calH_2}=0$ then yields
$\widetilde f\in\calD(H)$).
Since $\widetilde H$ is self-adjoint, $\ran(\widetilde H -z)=
\widetilde \calH$
for all $z\in\bbC\backslash\bbR$, and
$(\widetilde H \pm i):\widetilde \calH_2\to\widetilde
\calH$ is unitary,
\begin{equation}
((\widetilde H \pm i)\widetilde f,
(\widetilde H \pm i)\widetilde g))_{\tilde \calH}=
\int_{\bbR}(1+\lambda^2)d \omega(\lambda)
\overline{ \widetilde f(\lambda)}\widetilde g(\lambda)
=(\widetilde f,\widetilde g)_{\tilde \calH_2}, \quad
\widetilde f,\widetilde g\in\widetilde \calH_2. \lb{4.54}
\end{equation}
Thus, \eqref{4.52} and \eqref{4.54} yield
\begin{align}
\widetilde \calH&=(\widetilde H\pm i)\widetilde \calH_2 =
(\widetilde H\pm i)(\calD(H)\oplus_{\tilde \calH_2}R^{-1}\bbC)
 =(H\pm i)\calD(H)\oplus_{\tilde \calH}
\big \{\tfrac{\lambda\pm i}{1+\lambda^2}c
\,\big |\, c\in\bbC \big \} \no \\
& =\ran(H\pm i)\oplus_{\tilde \calH}\{c(\lambda \mp i)^{-1} \,|\,
c\in\bbC \} \lb{4.55}
\end{align}
and hence \eqref{4.56}.
\end{proof}

\medskip

If $H_\alpha$ and $H_\beta$ are two distinct self-adjoint
extensions
of the symmetric operator $H$ with deficiency indices $(1,1)$
considered in Theorem~\ref{t4.1b}, then, in
contrast to the case of deficiency indices $(n,n)$ to be
studied in detail
in Section~\ref{s7}, $\calD (H_\alpha)$ and $\calD
(H_\beta)$ have
a trivial intersection, that is,
\begin{equation}
\calD (H_\alpha) \cap \calD (H_\beta) = \calD (H)
\text{ for all }
\alpha, \beta \in [0,\pi), \alpha\neq\beta. \lb{4.56a}
\end{equation}

Introducing the Herglotz function
\begin{align}
m_{\alpha}(z)&=\int_{\bbR}d\omega_{\alpha}(\lambda)
((\lambda-z)^{-1}-
\lambda(1+\lambda^2)^{-1}) \lb{4.18} \\
&=z+(1+z^2)(u_+,(H_{\alpha}-z)^{-1}u_+)_{\calH} \lb{4.18a}
\end{align}
(the last equality being a simple consequence of
$\int_{\bbR}d\omega_{\alpha}(\lambda)
(1+\lambda^2)^{-1}=1$) one verifies
\begin{equation}
m_{\beta}(z)=\frac{-\sin(\beta-\alpha) +
\cos(\beta-\alpha) m_{\alpha}(z)}{\cos(\beta-\alpha) +
\sin(\beta-\alpha) m_{\alpha}(z)}, \quad
\alpha,\beta\in [0,\pi). \lb{4.19}
\end{equation}
A comparison of \eqref{4.19} and \eqref{3.10a} suggests
invoking
\begin{equation}
a(\alpha,\beta)=\left(\begin{array}{cc} \cos(\beta-\alpha)
& \sin(\beta-\alpha) \\ -\sin(\beta-\alpha) &
\cos(\beta-\alpha) \end{array}\right) \in\calA_2, \quad
\alpha,\beta\in [0,\pi). \lb{4.20}
\end{equation}
Moreover, since $m_{\gamma}(i)=i$ for all $\gamma\in [0,\pi)$,
Theorem~\ref{t3.2} applies (with $a_{1,1}(\alpha,\beta)=
a_{2,2}(\alpha,\beta)=
\cos(\beta-\alpha)$, $a_{1,2}(\alpha,\beta)=-a_{2,1}(\alpha,\beta)
=\sin(\beta-\alpha)$).\\

Next, assuming that $H$ is nonnegative, $H\geq 0$, we intend
to characterize
the Friedrichs and Krein extensions, $H_F$ and $H_K$, of $H$.
In order
to apply Krein's results \cite{Kr47a} (see also \cite{AT82},
\cite{Ts80},
\cite{Ts81}) in a slightly different
form (see, e.g., \cite{Sk79}, Sect.~4 for an efficient
summary of
Krein's results most relevant in our context) we state

\begin{theorem} \mbox{\rm } \lb{t4.1c}

\noindent (i). $H_{\alpha}=H_F$ for some $\alpha\in [0,\pi)$
if and only if
$\int_R^{\infty} d\|E_{\alpha}(\lambda)u_+\|_{\calH}^2 \lambda
=\infty$, or
equivalently, if and only if  $\int_R^{\infty}
d\omega_{\alpha}(\lambda) \lambda^{-1}=\infty$ for all $R>0$.

\noindent (ii). $H_{\beta}=H_K$ for some $\beta\in [0,\pi)$
if and only if
$\int^R_{0} d\|E_{\beta}(\lambda)u_+\|_{\calH}^2
\lambda^{-1}=\infty$, or
equivalently, if and only if $\int^R_{0}
d\omega_{\beta}(\lambda) \lambda^{-1} =\infty$ for all $R>0$.

\noindent (iii). $H_{\gamma}=H_F=H_K$ for some
$\gamma\in [0,\pi)$ if and only if
$\int_R^{\infty}d\|E_{\gamma}(\lambda)u_+\|_{\calH}^2
\lambda =$ $
\int^R_{0} d\|E_{\gamma}(\lambda)u_+\|_{\calH}^2
\lambda^{-1} =\infty$, or
equivalently, if and only if $\int_R^{\infty}
d\omega_{\gamma}(\lambda) \lambda^{-1} =\infty =$ $
\int^R_{0} d\omega_{\gamma}(\lambda) \lambda^{-1}$
for all $R>0$.
\end{theorem}
\begin{proof}
In order to reduce the above statements (i)--(iii) to
those in Krein
\cite{Kr47a} (as summarized in Skau \cite{Sk79}),
it suffices
to argue as follows. From  $(\mu +1)^{-1}=(\mu -i)^{-1}-
(1+i)(\mu +1)^{-1}(\mu -i)^{-1}$ one infers
\begin{align}
\|E_{\alpha}(\lambda)U_{\alpha}(\cdot +1)^{-1}\|_{\calH}^2 =
&\|E_{\alpha}(\lambda)u_+\|_{\calH}^2 +
2\|E_{\alpha}(\lambda)(H_{\alpha} +1)^{-1}u_+\|_{\calH}^2
\no \\
& -2\|E_{\alpha}(\lambda)(H_{\alpha} +1)^{-1/2}
u_+\|_{\calH}^2 \lb{4.20a}
\end{align}
and since
\begin{equation}
\int_{\calB}d\|E_{\alpha}(\lambda)(H_{\alpha} +1)^{-r}
u_+\|_{\calH}^2 =
\int_{\calB}\frac{d\|E_{\alpha}(\lambda)u_+\|_{\calH}^2}
{(\lambda +1)^{2r}},
\quad r\geq 0, \lb{4.20b}
\end{equation}
one concludes that
\begin{align}
&\int_{\calB} d\|E_{\alpha}(\lambda)U_{\alpha}(\cdot +1)^{-1}
\|_{\calH}^2
\text{ is finite (infinite)} \no \\
&\text{if and only if } \int_{\calB}
d\|E_{\alpha}(\lambda)u_+\|_{\calH}^2
\text{ is finite (infinite).} \lb{4.20c}
\end{align}
Here $\calB$ denotes any Borel subset of $[0,\infty)$.
\end{proof}

We also recall that
\begin{equation}
\inf \sigma (H_F)=\inf \{(g,Hg)_{\calH}\in\bbR\,|\,g\in\calD (H),
\|g\|_{\calH}=1\} \geq 0, \lb{4.21}
\end{equation}
whereas
\begin{equation}
\inf \sigma (H_K)=0 \lb{4.22}
\end{equation}
(here $\sigma(\cdot)$ abbreviates the spectrum of a
linear operator).
Moreover, all nonnegative self-adjoint extensions
$\widetilde H$ of
$H$ satisfy
\begin{equation}
0 \leq (H_F-\mu)^{-1}\leq (\widetilde H-\mu)^{-1}
\leq (H_K-\mu)^{-1},
\quad \mu\in (-\infty,0). \lb{4.22a}
\end{equation}
and hence $H$ has a unique nonnegative self-adjoint extension
if and
only if $H_K=H_F$.

Theorem~\ref{t4.1c} then yields the following result.

\begin{theorem} \mbox{\rm (\cite{DM91}, \cite{DM95},
\cite{DMT88},
\cite{Do65}, \cite{KO78}). } \lb{t4.1}

\noindent (i). $H_{\alpha}=H_F$ for some $\alpha\in [0,\pi)$
if and only if  $\lim_{\lambda\downarrow
-\infty} m_{\alpha}(\lambda)=-\infty$.

\noindent (ii). $H_{\beta}=H_K$ for some $\beta\in [0,\pi)$
if and only if  $\lim_{\lambda\uparrow 0}m_{\beta}(\lambda)
=\infty$.

\noindent (iii). $H_{\gamma}=H_F=H_K$ for some
$\gamma\in [0,\pi)$ if and only if $\lim_{\lambda\downarrow
-\infty}m_{\gamma}(\lambda)=-\infty$ and $\lim_{\lambda\uparrow
0}m_{\gamma}(\lambda)=\infty$.

\noindent (iv). Suppose $\alpha_F\in [0,\pi)$ corresponds to
$H_{\alpha_F} =H_F$, $\beta_K\in [0,\pi)$ to $H_{\beta_K}
=H_K$,
and $\gamma\in [0,\pi)$. Then
\begin{align}
&\lim_{\lambda\downarrow -\infty} m_{\gamma}(\lambda)=
-\cot(\gamma - \alpha_F)=-\int_{\bbR} d\omega_{\gamma}(\lambda)
\lambda(1+\lambda^2)^{-1}, \quad \gamma\neq\alpha_F,
\lb{4.22b} \\
&\lim_{\lambda\uparrow 0} m_{\gamma}(\lambda)
=-\cot(\gamma-\alpha_K)
=\int_{\bbR} d\omega_{\gamma}(\lambda)(\lambda^{-1}-
\lambda(1+\lambda^2)^{-1}), \quad \gamma\neq\alpha_K,
\lb{4.22c} \\
&\int_{\bbR}d\omega_{\gamma}(\lambda)\lambda^{-1}
=\cot(\gamma-\alpha_F)
- \cot(\gamma - \alpha_K), \quad \gamma\neq\alpha_F, \,\,
\gamma\neq\alpha_K. \lb{4.22d}
\end{align}
\end{theorem}

\begin{proof}
If $H_{\delta}\geq 0$ for some
$\delta\in [0,\pi)$ one infers
\begin{equation}
m_{\delta}(\lambda)=\int_0^{\infty} d\omega_{\delta}(\lambda')
((\lambda'+|\lambda|)^{-1}-\lambda'(1+{\lambda'}^2)^{-1}),
\quad \lambda <0. \lb{4.23}
\end{equation}
Next, suppose $H_{\alpha}=H_F$. Then, since
$[\lambda'(1+{\lambda'}^2)^{-1}
-(\lambda'+|\lambda|)^{-1}]$ is monotone increasing as
$\lambda\downarrow
-\infty$, $\lim_{\lambda\downarrow -\infty}m_{\alpha}(\lambda)=
-\int_0^{\infty} d\omega_{\alpha}(\lambda')\lambda'
(1+{\lambda'}^2)^{-1}
=-\infty$ by the monotone convergence theorem and
Theorem~\ref{t4.1b}\,(i).
Conversely, suppose $\lim_{\lambda\downarrow -
\infty}m_{\alpha}(\lambda)
=-\infty$, then necessarily $\int_0^{\infty}
d\omega_{\alpha}(\lambda)
(1+\lambda)^{-1}=\infty$ and hence $H_{\alpha}=H_F$
again by
Theorem~\ref{t4.1b}\,(i).
This proves (i). Items (ii) and (iii) then follow
analogously from
Theorem~\ref{t4.1b}\,(ii)  and (iii) above. Equation
\eqref{4.22b} is a direct consequence of
\eqref{4.18}, \eqref{4.19}, (i)--(iii), and the fact
that all operators $H_{\alpha}$,
$\alpha\in [0,\pi)$ are bounded from below (and hence
$m_{\alpha}(z)$
are real-valued for $z\in (-\infty,-c(\alpha)]$ and
analytic in
$\bbC\backslash (-\infty,-c(\alpha)]$ for $c(\alpha)>0$
sufficiently large). Equation \eqref{4.22c} is proved in
the same manner
observing that $\sigma(H_{\alpha}) \cap (-\infty,0)$,
$\alpha\in [0,\pi)$ consists of at most one eigenvalue.
Finally,
\eqref{4.22d} is just the difference of \eqref{4.22c} and
\eqref{4.22b}.
\end{proof}

The following represents an elementary example illustrating
these concepts.

\begin{example} \lb{e4.2}
Let $r\in (-1,1)$ and consider the measure
\begin{equation}
d\mu_r(\lambda)=\begin{cases} (2/\pi)\sin((r+1)\pi/2)
\lambda^r d\lambda,
\quad \lambda\geq 0 \\
0, \quad \lambda <0. \end{cases} \lb{4.24}
\end{equation}
Then \eqref{4.12} is easily verified and
\begin{equation}
\int_R^{\infty}\frac{d\mu_r(\lambda)}{\lambda}=\begin{cases}
\infty
\quad \text{if } 0\leq r <1 \\
< \infty \quad \text{if } -1 <r < 0, \end{cases} \quad
\int_0^R\frac{d\mu_r(\lambda)}{\lambda}=\begin{cases} \infty
\quad \text{if } -1< r \leq 0 \\
< \infty \quad \text{if }  0 < r < 1 \end{cases} \lb{4.25}
\end{equation}
for all $R>0$. Define the closed symmetric operator
${\hatt H}(r)\geq 0$
in $L^2((0,\infty);d\mu_r)$ with deficiency indices
$(1,1)$ by
\begin{align}
({\hatt H}(r)\hat g)(\lambda)&=\lambda \hat g(\lambda),
\lb{4.26} \\
\hat g\in\calD ({\hatt H}(r))&=\{\hat h \in L^2((0,\infty);
(1+\lambda^2)d\mu_r)\,|\,
\smallint_0^{\infty} d\mu_r(\lambda)\hat h(\lambda)=0\} \no
\end{align}
and the self-adjoint (maximally defined multiplication)
operator
\begin{equation}
({\hatt H}_r \hat g)(\lambda)=\lambda \hat g(\lambda),
\quad
\hat g\in\calD ({\hatt H}_r) = L^2((0,\infty);
(1+\lambda^2)d\mu_r).
 \lb{4.27}
\end{equation}
Then ${\hatt H}_r$ represents the Friedrichs extension
${\hatt H}(r)_F$
of ${\hatt H}(r)$ for $0\leq r < 1$ and the Krein extension
${\hatt H}(r)_K$
of ${\hatt H}(r)$ for $-1 < r \leq 0$. In particular,
${\hatt H}(r)_F=
{\hatt H}(r)_K$ if and only if $r=0$.
\end{example}

Next, we turn to a realization theorem for Herglotz
functions of the type
\eqref{4.18a}. It will be convenient to introduce the
following sets of Herglotz functions,
\begin{align}
\calN_0 &=\{m:\bbC_+\to\bbC_+ \text{ analytic }|\,m(z)=
\smallint_{\bbR}
d\omega(\lambda)((\lambda -z)^{-1} -
\lambda(1+\lambda^2)^{-1}), \no \\
& \hspace*{.63cm}  \smallint_{\bbR} d\omega(\lambda)
=\infty, \,\,
\smallint_{\bbR} d\omega(\lambda)
(1+\lambda^2)^{-1} <\infty \}, \lb{4.28} \\
\calN_{0,F} &=\{m\in\calN_0\,|\, \supp(\omega)
\subseteq[0,\infty),\,\,
\smallint_R^{\infty} d\omega(\lambda)\lambda^{-1}
=\infty \text{ for some }R>0 \},
\lb{4.29} \\
\calN_{0,K} &=\{m\in\calN_0\,|\, \supp(\omega)
\subseteq[0,\infty),\,\,
\smallint_0^{R} d\omega(\lambda)\lambda^{-1}
=\infty \text{ for some }R>0 \},
\lb{4.30} \\
\calN_{0,F,K} &=\{m\in\calN_0\,|\, \supp(\omega)
\subseteq[0,\infty),\,\,
\smallint_R^{\infty} d\omega(\lambda)\lambda^{-1}
=\smallint_0^{R} d\omega(\lambda)\lambda^{-1}
=\infty \no \\
& \hspace*{.53cm}\text{ for some }R>0 \}
=\calN_{0,F}\cap\calN_{0,K}. \lb{4.31}
\end{align}
The sets $\calN_{0,F}$, $\calN_{0,K}$, and
$\calN_{0,F,K}$ are of course
independent of $R>0$.

\begin{theorem}  \mbox{\rm } \lb{t4.3}

\noindent (i). Any \, $\widetilde m\in\calN_0$ can
be realized in
the form
\begin{equation}
\widetilde m(z)=z\|u_+\|^2_{\tilde\calH}+(1+z^2)
(u_+,(\widetilde H -z)^{-1}u_+)_{\tilde\calH},
\quad z\in\bbC_+, \lb{4.32}
\end{equation}
where $\widetilde H$ denotes the self-adjoint extension
of some densely
defined closed symmetric operator $H$ with deficiency
indices $(1,1)$
and deficiency vector $u_+\in\ker(H^{*}-i)$ in
some separable complex Hilbert space $\widetilde \calH$.

\noindent (ii). Any \, $\widetilde m_{F(\text{resp.~}K)}\in
\calN_{0,F(\text{resp.~}K)}$ can be realized in the form
\begin{equation}
\widetilde m_{F(\text{resp.~}K)}(z)=z\|u_+\|^2_{\tilde\calH}+
(1+z^2)(u_+,(\widetilde H_{F(\text{resp.~}K)} -z)^{-1}
u_+)_{\tilde\calH}, \quad z\in\bbC_+, \lb{4.33}
\end{equation}
where $\widetilde H_{F(\text{resp.~}K)}\geq 0$ denotes
the Friedrichs
(respectively, Krein) extension of some densely
defined closed operator $H\geq 0$ with deficiency
indices $(1,1)$
and deficiency vector $u_+\in\ker(H^{*}-i)$ in
some separable complex Hilbert space $\widetilde \calH$.

\noindent (iii). Any \, $\widetilde m_{F,K}\in
\calN_{0,F,K}$ can be realized in the form
\begin{equation}
\widetilde m_{F,K}(z)=z\|u_+\|^2_{\tilde\calH}+(1+z^2)
(u_+,(\widetilde H_{F,K} -z)^{-1}u_+)_{\tilde\calH},
\quad z\in\bbC_+,
\lb{4.34}
\end{equation}
where $\widetilde H_{F,K}\geq 0$ denotes the unique
nonnegative self-adjoint
extension of some densely
defined closed operator $H\geq 0$ with deficiency
indices $(1,1)$
and deficiency vector $u_+\in\ker(H^{*}-i)$ in
some separable complex Hilbert space $\widetilde \calH$.

\noindent In each case (i)--(iii) one has
\begin{equation}
\int_{\bbR}d\widetilde \omega(\lambda)(1+\lambda^2)^{-1}=
\|u_+\|^2_{\tilde\calH}, \lb{4.46}
\end{equation}
where $\widetilde \omega$ denotes the measure in the
Herglotz representation of $\widetilde m(z)$.
\end{theorem}

\begin{proof}
We use the notation established in Theorem~\ref{t4.1b}.
Define
\begin{equation}
u_+(\lambda)=(\lambda -i)^{-1},  \lb{4.57}
\end{equation}
then $\|u_+\|^2_{\tilde\calH}=\smallint_{\bbR}
d \widetilde \omega(\lambda)
(1+\lambda^2)^{-1}$ and
\begin{align}
&z\|u_+\|^2_{\tilde\calH}+
(1+z^2)(u_+,(\widetilde H -z)^{-1}u_+)_{\tilde\calH}
\no \\
&=\int_{\bbR} d \widetilde \omega(\lambda)
(z(1+\lambda^2)^{-1} +
(1+z^2)(\lambda -z)^{-1}(1+\lambda^2)^{-1}) \no \\
&=\int_{\bbR} d \widetilde \omega(\lambda)
((\lambda -z)^{-1} -
\lambda(1+\lambda^2)^{-1}) =\widetilde m(z) \lb{4.58}
\end{align}
proves \eqref{4.32} and hence part (i). Parts (ii)
and (iii) then
follow in the same manner from Theorems~\ref{t4.1b}
and \ref{t4.1c}.
\end{proof}

Of course we could have normalized $u_+$,
$\|u_+\|_{\tilde\calH}=1$, and
then added the constraint $\smallint_{\bbR} d\omega(\lambda)
(1+\lambda^2)^{-1}=1$ to \eqref{4.28}--\eqref{4.31}. By
\eqref{4.19e}, \eqref{4.32}--\eqref{4.34} can be realized in
nonseparable Hilbert spaces.

\begin{theorem} \lb{t4.5}
Suppose $m_\ell\in\calN_0$ with corresponding measures
$\omega_\ell$ in
the Herglotz representation of $m_\ell$, $\ell=1,2$,
and $m_1\neq m_2$.
Then $m_1$ and $m_2$ can be realized as
\begin{equation}
m_\ell(z)=z\|u_+\|_{\calH}^2 +(1+z^2)(u_+,(H_\ell -z)^{-1}
u_+)_{\calH},
\quad \ell=1,2, \,\, z\in\bbC_+, \lb{4.69}
\end{equation}
where $H_\ell$, $\ell=1,2$ are distinct self-adjoint
extensions of one
and the same
densely defined closed symmetric operator $H$ with
deficiency indices
$(1,1)$ and deficiency vector $u_+\in\ker(H^*-i)$ in
some complex
Hilbert space $\calH$ (which may be chosen to be separable)
if and only if
the following conditions hold:
\begin{equation}
\int_{\bbR} d\omega_1(\lambda)(1+\lambda^2)^{-1}=
\int_{\bbR} d\omega_2(\lambda)(1+\lambda^2)^{-1}=
\|u_+\|^2_{\calH},
\lb{4.70}
\end{equation}
and for all $z\in\bbC_+$,
\begin{equation}
m_2(z)=\frac{-\|u_+\|^2_{\calH}+h\,m_1(z)}
 {h+\|u_+\|^{-2}_{\calH}m_1(z)}
 \text{ for some }h\in\bbR. \lb{4.71}
\end{equation}
\end{theorem}

\begin{proof}
The necessity of condition \eqref{4.71} has been
proved by Donoghue
\cite{Do65} (he assumed $\|u_+\|_{\calH}=1$). Indeed,
assuming
\eqref{4.69}, the fact
\begin{equation}
m_\ell(i)=i\|u_+\|^2_{\calH}=i\int_{\bbR}d\omega_\ell(\lambda)
(1+\lambda^2)^{-1}, \quad \ell=1,2 \lb{4.72}
\end{equation}
yields \eqref{4.70}. Identifying $h=\cot(\beta-\alpha)$,
$H_1=H_{\alpha}$, $H_2=H_{\beta}$,
$\|u_+\|^{-2}_{\calH}m_1(z)=m_{\alpha}(z)$,
 and $\|u_+\|^{-2}_{\calH}m_2(z)$ $=m_{\beta}(z)$,
\eqref{4.71} is seen to be equivalent to \eqref{4.19}.
(Here we may,
without loss of generality, assume that $H_\ell$,
$\ell=1,2$ have simple
spectra since otherwise one can apply the reduction
\eqref{4.19e}.)
Conversely, assume \eqref{4.70} and \eqref{4.71}. By
Theorem~\ref{t4.3}\,(i), we may realize $m_1(z)$ as
\begin{equation}
\|u_+\|_{\calH}^{-2}m_1(z)=z +(1+z^2)\|u_+\|_{\calH}^{-2}
(u_+,(H_1 -z)^{-1}u_+)_{\calH}. \lb{4.73}
\end{equation}
Again by \eqref{4.19e} we may assume that $H_1$ has simple
spectrum and identify it with $H_{\alpha}$,
$\alpha\in [0,\pi)$ in \eqref{4.9}. If $H_\beta$,
$\beta\in [0,\pi)\backslash\{\alpha\}$ is any other
self-adjoint extension of $H$ defined as in \eqref{4.9}
(the actual
normalization of $u_{\pm}$ being immaterial in this context),
introduce
\begin{equation}
m_\beta(z)=z+(1+z^2)\|u_+\|^{-2}_{\calH}(u_+,
(H_\beta -z)^{-1}u_+)_{\calH}. \lb{4.74}
\end{equation}
By \eqref{4.19} one obtains ($m_{\alpha}(z)=
\|u_+\|^{-2}_{\calH}m_1(z)$)
\begin{equation}
m_\beta(z)=
\frac{-1+\cot(\beta-\alpha)\|u_+\|^{-2}_{\calH}m_1(z)}
{\cot(\beta-\alpha)+\|u_+\|^{-2}_{\calH}m_1(z)}. \lb{4.75}
\end{equation}
A comparison of \eqref{4.71} and \eqref{4.75} then yields
$\|u_+\|^2_{\calH}m_\beta(z)=m_2(z)$ for
$h=\cot(\beta-\alpha)$,
completing the proof.
\end{proof}

\begin{remark} \lb{rm4.7}
For simplicity we studied Friedrichs $H_F$ and Krein $H_K$
extensions of a
densely defined closed operator $H\geq 0$ with deficiency
indices $(1,1)$ in Theorems~\ref{t4.1} and
\ref{t4.3}\,(ii),(iii).
In other words, we studied the special case where $H$
admitted at
least one self-adjoint extension with the spectral
gap $(-\infty,0)$
(in general, there is a one-parameter family of such
self-adjoint
extensions with $H_F$ and $H_K$ as extreme points, cf.,
\eqref{4.22a}).
There is no difficulty in extending all our results to
the case of
symmetric operators and their self-adjoint extensions with
arbitrary gaps
$(\lambda_1,\lambda_2)$, $-\infty < \lambda_1 < \lambda_2
< \infty$ in
their spectrum. In
fact, assuming $H$ to be densely defined and closed in
some complex
Hilbert space $\calH$, the condition $H\geq 0$ is now
replaced by
\begin{equation}
\|(H-\tfrac{\lambda_1 + \lambda_2}{2})f\|_{\calH}
\geq\tfrac{\lambda_2-\lambda_1}{2}\|f\|_{\calH}, \quad
f\in\calD(H).
\lb{4.75a}
\end{equation}
In this situation it was proved by Krein \cite{Kr47a} that
$H$ always admits
self-adjoint extensions $\widetilde H$ with the same spectral
gap
$(\lambda_1,\lambda_2)$. In particular, there always exist
two extremal
self-adjoint extensions $H_{F_{\lambda_1}}$ and
$H_{K_{\lambda_2}}$ of
$H$ with the same gap $(\lambda_1,\lambda_2)$ such that
\begin{equation}
(H_{F_{\lambda_1}}-\mu)^{-1}\leq (\widetilde H-\mu)^{-1}
\leq (H_{K_{\lambda_2}}-\mu)^{-1},
\quad \mu\in (\lambda_1,\lambda_2) \lb{4.22ab}
\end{equation}
for any self-adjoint extension $\widetilde H$ of $H$ with
spectral gap
$(\lambda_1,\lambda_2)$. Given the results in \cite{AT82},
\cite{DM91},
\cite{DM95}, \cite{DMT88}, and \cite{Kr47a},
Theorem~\ref{t4.1} immediately extends to general
gaps $(\lambda_1,\lambda_2)$ upon replacing $H_F$ by
$H_{F_{\lambda_1}}$,
$\lim_{\lambda\downarrow -\infty} m_\alpha(\lambda)
=-\infty$ by
$\lim_{\lambda\downarrow \lambda_1, \, \lambda\in (\lambda_1,
\lambda_2)}$
$m_\alpha(\lambda)=-\infty$, $H_K$ by $H_{K_{\lambda_2}}$,
and
$\lim_{\lambda\uparrow 0}m_\beta(\lambda)=\infty$ by
$\lim_{\lambda\uparrow \lambda_2, \, \lambda\in (\lambda_1,
\lambda_2)}
m_\beta(\lambda)$ $=\infty$, etc. Analogous remarks apply to
Theorem~\ref{t4.3}\,(ii),(iii), replacing the condition
$\supp(\omega)\subseteq [0,\infty)$ by $\supp(\omega)
\subseteq\bbR\backslash
(\lambda_1,\lambda_2)$ in \eqref{4.29}--\eqref{4.31}.
\end{remark}

Next we briefly turn to Schr\"odinger operators on a
half-line. Let
$q\in L^1 ([0,R])$ for all $R>0$, $q$
real-valued, and introduce the fundamental system
$\phi_{\alpha}(z,x)$,
$\theta_\alpha(z,x)$, $z\in\bbC$ of solutions of ($\prime$
denotes $d/dx$)
\begin{equation}
-\psi''(z,x)+(q(x)-z)\psi(z,x)=0, \quad x > 0, \lb{4.76}
\end{equation}
satisfying
\begin{equation}
\phi_\alpha (z,0_+)=-\theta'_\alpha (z,0_+)=-\sin(\alpha),
\,\,
\phi'_\alpha (z,0_+)=\theta_\alpha (z,0_+)=\cos(\alpha),
\,\,
\alpha\in [0,\pi). \lb{4.77}
\end{equation}
Next, pick a fixed $z_0\in\bbC_+$ and a solution
$f_0(z_0,\cdot)\in L^2([0,\infty);dx)$ of \eqref{4.76} and
let
$\psi_\alpha(z,x)$ be the unique solution of \eqref{4.76}
satisfying
\begin{align}
& \psi_\alpha (z,\cdot)\in L^2([0,\infty);dx), \quad
\sin(\alpha)\psi'_\alpha (z,0_+)+\cos(\alpha)
\psi_\alpha (z,0_+)=1,
\no \\
& \underset{x\to\infty}{\lim} W(f_0(z_0,x),
\psi_\alpha (z,x))=0, \quad
z\in\bbC_+, \lb{4.78}
\end{align}
the latter condition being superfluous, i.e., automatically
fulfilled, if
$-\tfrac{d^2}{dx^2}+q$ is in the limit point case at
$\infty$. (Here
$W(f(x),g(x))=f(x)g'(x)-f'(x)g(x)$ denotes the Wronskian
of $f$ and $g$.)
Existence and
uniqueness of $\psi_\alpha (z,x)$ is a consequence of
Weyl's theory
(see, e.g., the discussion in Appendix~A of \cite{GS96a}).
Then
$\psi_\alpha (z,x)$ is of the form
\begin{equation}
\psi_\alpha (z,x)=\theta_\alpha(z,x)+m_\alpha (z)
\phi_\alpha (z,x), \lb{4.79}
\end{equation}
with $m_\alpha (z)$ the Weyl-Titchmarsh $m$-function
corresponding to the
operator $H_\alpha$ in $L^2([0,\infty);dx)$ defined by
\begin{align}
&(H_\alpha g)(x)=-g''(x)+q(x)g(x), \quad x>0, \lb{4.79a} \\
&\calD (H_\alpha)=\{g\in L^2([0,\infty);dx)\,
|\,g,g'\in AC([0,R])
\text{ for all $R>0$ }; \no \\
& \hspace*{1.7cm} -g''+qg\in L^2([0,\infty);dx); \,\,
\underset{x\to\infty}{\lim} W(f_0(z_0,x),g(x))=0; \no \\
& \hspace*{1.7cm} \sin(\alpha)g'(0_+)+\cos(\alpha)g(0_+)
=0 \}, \quad
\alpha\in [0,\pi). \no
\end{align}
(Here $AC([a,b])$ denotes the set of absolutely continuous
functions on
$[a,b]$.) Then $m_\alpha (z)$ is a
Herglotz function with representation
\begin{align}
m_\alpha (z) &=c_\alpha+\int_{\bbR}
d\omega_\alpha (\lambda)((\lambda -z)^{-1}-\lambda
(1+\lambda^2)^{-1}),
&\alpha\in [0,\pi), \lb{4.80a} \\
& = \cot(\alpha)+\int_{\bbR} d\omega_\alpha (\lambda)
(\lambda -z)^{-1},
&\alpha\in (0,\pi), \lb{4.80b}
\end{align}
where
\begin{equation}
\int_{\bbR} d\omega_\alpha (\lambda)(1+|\lambda|)^{-1}
\begin{cases}
<\infty, \,\, &\alpha\in (0,\pi), \\
=\infty, \,\, &\alpha=0.
\end{cases} \lb{4.81}
\end{equation}
Moreover, one verifies
\begin{equation}
m_{\beta}(z)=\frac{-\sin(\beta-\alpha) +
\cos(\beta-\alpha) m_{\alpha}(z)}{\cos(\beta-\alpha) +
\sin(\beta-\alpha) m_{\alpha}(z)}, \quad
\alpha,\beta\in [0,\pi) \lb{4.82}
\end{equation}
and hence the corresponding matrix $a(\alpha,\beta)$ is of
the type
\begin{equation}
a(\alpha,\beta)=\left(\begin{array}{cc} \cos(\beta-\alpha)
& \sin(\beta-\alpha) \\ -\sin(\beta-\alpha) &
\cos(\beta-\alpha) \end{array}\right) \in\calA_2, \quad
\alpha,\beta\in [0,\pi). \lb{4.83}
\end{equation}
The asymptotic behavior of $m_\alpha (z)$ is given by
\begin{equation}
m_\alpha (z) \underset{z\to i\infty}{=}
\begin{cases} \cot(\alpha) + \frac{i}{\sin^2(\alpha)}z^{-1/2}
-\frac{\cos(\alpha)}{\sin^3(\alpha)}z^{-1}+\oh (z^{-1}), \,
&\alpha \in (0,\pi), \\
iz^{1/2} + \oh (1), \, &\alpha=0. \end{cases} \lb{4.84}
\end{equation}
Thus, Theorem~\ref{t3.2} applies (with $a_{1,1}(\alpha,
\beta)=
a_{2,2}(\alpha,\beta)=\cos(\beta-\alpha)$,
$a_{1,2}(\alpha,\beta)=-a_{2,1}(\alpha,\beta)=\sin(\beta-
\alpha)$).

Theorem~\ref{t3.2}\,(v), in particular, represents an
alternative
(abstract) approach to Borg-type uniqueness theorems
\cite{Bo46},
\cite{Bo52} (see also \cite{GS96a}, \cite{Le49}, \cite{Le87},
\cite{Ma73} and the
references therein) to the effect that two sets of spectra
(varying
the boundary condition at one end point but keeping it
fixed at the
other) uniquely determine $q(x)$. Its elegant proof using
the exponential
Herglotz representation for $F(z)$ is due to Donoghue
\cite{Do65}.

For simplicity we only discussed the case of a half-line
$[0,\infty)$.
However, the case of a finite interval $[0,R_0]$ for some
$R_0>0$ is
completely analogous, replacing the first and third
condition on
$\psi_\alpha (z,x)$ in \eqref{4.78} by the boundary
condition
$\sin(\gamma)\psi'_\alpha (z,R_0)+\cos(\gamma)\psi_\alpha
(z,R_0)=0$
for some fixed $\gamma\in [0,\pi)$.

It is possible to characterize the set of Herglotz
functions leading
to Weyl-Titchmarsh $m$-functions for $-\tfrac{d^2}{dx^2}
+q$ in
$L^2([0,R_0];dx)$
or $L^2([0,\infty);dx)$ with real-valued $q$ satisfying
$q\in L^1([0,R_0];dx)$ or $q\in L^1([0,R);dx)$ for all
$R>0$, respectively. These
realization theorems, however, are far less elementary,
being
based on
the Gelfand-Levitan formalism of inverse spectral theory
(see, e.g.,
\cite{GL55}, \cite{Le87}, \cite{LG64}, \cite{Ma86}). We
omit further
details at this point.

These considerations extend to singular coefficients $q$
at $x=0$
replacing $q\in L^1([0,R])$ for all $R>0$ by
$q\in L^1_{\loc} ((0,\infty))$. A careful investigation
of the Weyl
limit point/limit circle theory (see, e.g., \cite{CL85},
Ch.~9)
then shows that the fundamental system $\phi_\alpha(z,x)$,
$\theta_\alpha(z,x)$ of \eqref{4.76} can be replaced by
$\phi (z,x)$, $\theta (z,x)$ satisfying \eqref{4.76} with
the
following properties:

\vspace{1.2mm}

\noindent (i). For all $x>0$, $\phi(z,x)$, $\theta (z,x)$
are
entire with respect to $z\in\bbC$ and real-valued for all
$z\in\bbR$.

\vspace{.8mm}

\noindent (ii). $W(\theta (z,x),\phi (z,x))=1$,
$z\in\bbC$, $x>0$.

\vspace{.8mm}

\noindent (iii). $\lim_{x\downarrow 0} W(\phi (z,x),\ol
{\phi (z,x)})=
\lim_{x\downarrow 0} W(\theta (z,x),\ol {\theta (z,x)})
=0$, \\
\hspace*{.73cm} $\lim_{x\downarrow 0} W(\theta (z,x),\ol
{\phi(z,x)})
=1$, $z\in\bbC$.

\vspace{1.2mm}

Introducing
\begin{equation}
\psi (z,x)=\theta (z,x)+m(z)\phi (z,x), \quad z
\in\bbC\backslash\bbR
\lb{4.87}
\end{equation}
satisfying
\begin{equation}
\psi (z,\cdot)\in L^2 ([0,\infty);dx), \quad z
\in\bbC\backslash\bbR
\lb{4.88}
\end{equation}
then yields
\begin{align}
&\Im(m(z))=\Im(z) \int_0^\infty dx |\psi (z,x)|^2, \quad
z\in\bbC\backslash\bbR \lb{4.90} \\
& \text{if and only if } \lim_{x\uparrow \infty} W(\psi (z,x),
\ol {\psi (z,x)})=0, \quad z\in\bbC\backslash\bbR. \lb{4.91}
\end{align}
In particular, $m(z)$ is a Herglotz function if \eqref{4.91}
is
satisfied. For associated self-adjoint boundary conditions
in the
singular case see, for instance, \cite{BG85}, \cite{JR76},
Ch.~III,
and \cite{Re43}.

\section{Basic Facts on Matrix-Valued Herglotz Functions}
\lb{s5}

The main purpose of this section is to carry over some of
the scalar
results of Section~\ref{s2} to matrix-valued Herglotz
functions.

In the following we denote by $M_n(\mathbb{C})$, $n\in
\mathbb{N}$ the
set of $n\times n$ matrices with complex-valued entries,
denote by
$I_n\in M_n(\mathbb{C})$ the identity matrix, by $A^*$
the adjoint
(complex conjugate transpose) of
$A\in M_n(\mathbb{C})$, and by $(\cdot ,
\cdot)_{\mathbb{C}^n}$
the scalar
product in $\mathbb{C}^n$ associated with the standard
Euclidean metric
on $\mathbb{C}^n$ (antilinear in the first and linear in
the second
factor). We recall that a matrix $A\in M_n(\mathbb{C})$
is called
nonnegative (respectively, nonpositive), $A\geq 0$
(respectively,
$A\leq 0$) if $(\ul x,A\ul x)_{\mathbb{C}^n}\geq 0$
(respectively,
$(\ul x,A\ul x)_{\mathbb{C}^n}\leq 0$) for all $\ul x\in
\mathbb{C}^n$. Similarly, $A$ is called positive (positive
definite, or
strictly positive), $A>0$, if
$(\ul x,A\ul x)_{\mathbb{C}^n}>0$ for all $\ul x\in
\mathbb{C}^n\backslash \{0\}$. A
principal submatrix of $A$ is obtained by deleting $k$
rows and
columns, $0\leq k\leq n-1$, which pairwise intersect at
diagonal
elements. Principal minors are determinants of principal
submatrices.
The rank, range, and kernel of $A$ are denoted by $\rank (A)$,
$\ran(A)$, and $\ker(A)$, respectively.

We start with an elementary result on nonnegative matrices
which will be
useful at various places later on.

\begin{lemma} [\cite{HJ94}, Ch.~7] \label{l5.1}
Let $A=(A_{j,k})_{1\leq
j,k\leq n}\in M_n(\mathbb{C})$ and assume $A\geq 0$. Then

\noindent (i). $A\geq 0$ if and only if all principal
minors of $A$ are nonnegative. In particular, all diagonal
elements of
$A$ are nonnegative,
\begin{equation} \lb{5.1}
A_{j,j}\geq 0,\quad 1\leq j\leq n.
\end{equation}

\noindent (ii). For all $1\leq j,k\leq n$,
\begin{equation} \lb{5.2}
|A_{j,k}|\leq A_{j,j}^{1/2}A^{1/2}_{k,k}\leq
\frac{1}{2}(A_{j,j}+A_{k,k}),
\end{equation}
in particular, if $A_{\ell,\ell}=0$
then the $\ell$th row and column of $A$ are zero.

\noindent (iii). Let $\ul x\in \mathbb{C}^n$ and
$(\ul x,A\ul x)_{\mathbb{C}^n}=0$. Then
$A\ul x=0$.

\noindent (iv). Suppose $\rank (A)=r<n$. Then $A$ has an
$r\times r$ positive
definite principal submatrix.
\end{lemma}

Next we briefly turn to (self-adjoint) matrix-valued
measures. The ones
to be used below will be of the type
\begin{equation}
\Sigma(\mathcal{M})=\int_{\mathcal{M}}d\Omega
(\lambda )(1+\lambda ^2)^{-1}, \quad
\Sigma =(\Sigma_{j,k})_{1\leq j,k\leq n},\quad \Omega
=(\Omega_{j,k})_{1\leq j,k\leq n},\lb{5.4}
\end{equation}
where
$\Sigma_{j,k}$, $1\leq j,k\leq n$ are complex (and hence
finite) Borel
measures on $\mathbb{R}$ and
$\Omega_{j,k}$, $1\leq j,k\leq n$ are complex-valued set
functions
defined on the bounded Borel subsets of $\mathbb{R}$ with
the properties,
\begin{enumerate}
\item[(i).] $\Omega (X)=(\Omega_{j,k}(X))_{1\leq j,k\leq
n}\subset M_n(\mathbb{C})$ is nonnegative,
$\Omega (X)\geq 0$, for all bounded Borel sets $X\subset
\mathbb{R}$, and
$\Omega (\phi )=0$.
\item[(ii).] $\Omega _{j,k}(\cup_{\ell\in \mathbb{N}}X_{\ell})
=\Sigma_{\ell\in
\mathbb{N}}\Omega_{j,k}(X_{\ell})$, $1\leq j,k\leq n$ for
each sequence of
disjoint Borel sets $\{X_{\ell}\}_{\ell\in \mathbb{N}}
\subset \mathbb{R}$ with
$\cup_{\ell\in \mathbb{N}}X_{\ell}$ bounded.\end{enumerate}

Clearly, each diagonal element $\Sigma_{j,j}$, $1\leq j
\leq n$ defines a
positive (finite) Borel measure on $\mathbb{R}$. In
addition, we denote by
\begin{equation}\lb{5.5}
\sigma^{tr}=\tr_{\mathbb{C}^n}(\Sigma)=\Sigma_{1,1}+\cdots
+\Sigma_{n,n}\end{equation} the (scalar) trace measure of
$\Sigma$ and
note that \begin{equation}\sigma^{tr}(X)=0\quad\text{if
and only
if}\quad\Sigma (X)=0\lb{5.6}\end{equation} for all Borel
sets
$X\subseteq \mathbb{R}$ since by \eqref{5.2} each
$\Sigma_{j,k}$ is
absolutely continuous with respect to $\Sigma_{j,j}+
\Sigma_{k,k}$ and
hence with respect to $\sigma^{tr}$,
\begin{equation} \lb{5.7}
\Sigma_{j,k}\ll\Sigma_{j,j}+\Sigma_{k,k}\ll
\sigma^{tr},\quad 1\leq j,k\leq n.
\end{equation}
Below we will use the standard Lebesgue decomposition
of matrix-valued
measures with respect to Lebesgue measure on $\bbR$, in
particular, we will
use the fact that $\Omega=\Omega_{ac}$ is purely absolutely
continuous
with respect to Lebesgue measure $d\lambda$ if
and only if $d\Omega(\lambda)=P(\lambda)d\lambda$ for
some nonnegative
locally integrable matrix $P$ on $\bbR$.

Matrix-valued Herglotz
functions are now defined in analogy to
Definition~\ref{d2.1} as follows.

\begin{definition}\label{d5.2}
$M:\mathbb{C_+}\rightarrow M_n(\mathbb{C})$ is called a
matrix-valued Herglotz function (in short, a Herglotz
matrix) if
 $M$ is analytic on $\mathbb{C}_+$ and $\Im (M(z))\geq 0$
for all $z\in
\mathbb{C}_+$.
\end{definition}
As in the scalar case one usually extends $M$ to
$\bbC_-$ by
reflection, that is, by defining
\begin{equation}
M(z)=M(\overline z)^*, \quad z\in\mathbb{C}_-.
\end{equation}
Hence $M$ is analytic on $\bbC\backslash\bbR$ but
$M\big|_{\bbC_-}$
and $M\big|_{\bbC_+}$, in general, are not analytic
continuations of
each other (cf.~Lemma~\ref{l5.5}). Here we follow the
standard notation
\begin{equation}\lb{5.8}
\Im (M)=\frac{1}{2i}(M-M^*),\quad \Re
(M)=\frac{1}{2}(M+M^*).
\end{equation}

In contrast to the scalar case, we cannot in general expect
strict
inequality in $\Im(M(z))\geq 0$. However the kernel of
$\Im(M(z))$
has extremely simple properties. The following result
and its
elementary proof were communicated to us by Dirk Buschmann:

\begin{lemma} \lb{l5.2a}
Let $M(z)\in M_n(\bbC)$ be a matrix-valued Herglotz function.
Then
the kernel $\ker(\Im(M(z)))$ is independent of
$z\in\bbC_+$, in
particular, the rank $r$ of $M(z)$ is constant on $\bbC_+$.
Consequently, upon choosing an orthogonal basis in
$\ker(M(z))$ and
$\ker(M(z))^\bot$, $M(z)$ takes on the form
\begin{equation}
M(z)=\left(\begin{array}{cc} 0 & 0 \\ 0 &
M_r(z) \end{array}\right), \lb{5.8a}
\end{equation}
where $M_r(z)$ is an $r\times r$ matrix-valued Herglotz
matrix
satisfying
\begin{equation}
\Im(M_r(z))>0, \quad  r=n-\dim_\bbC (\ker(\Im(M(z))),
\quad z\in\bbC_+.
\lb{5.8b}
\end{equation}
\end{lemma}
\begin{proof}
Denote $N_z=\ker(\Im(M(z)))$, $z\in\bbC_+$. Pick a
$z_0\in\bbC_+$ and
suppose $0\neq\ul x_0\in N_{z_0}$. Consider the scalar
Herglotz function
$m(z)=(\ul x_0,M(z)\ul x_0)_{\bbC^n}$. Then $m(z_0)\in\bbR$
shows that
the Herglotz function $m(z)-m(z_0)$ has a zero at
$z=z_0\in\bbC_+$ and
hence vanishes identically. Thus $m(z)$ equals a real
constant for all
$z\in\bbC_+$ and hence
\begin{equation}
0=(\ul x_0,\Im(M(z))\ul x_0)_{\bbC^n}=
\|(\Im(M(z)))^{1/2}\ul x_0\|^2_{\bbC^n}, \quad z\in\bbC_+
\lb{5.8c}
\end{equation}
yields
\begin{equation}
\ul x_0 \in\ker((\Im(M(z)))^{1/2})=\ker(\Im(M(z))),
\quad z\in\bbC_+
\lb{5.8d}
\end{equation}
since $\Im(M(z))\geq 0$. In particular, $N_z$, and hence
$r=\rank(\Im(M(z)))$$=n-\dim_\bbC (N_z)$ are
independent of $z\in\bbC_+$. Finally, suppose
$\ker(\Im(M(z_0)))=\{0\}$.
Then $\ker(\Im(M(z_1)))\neq \{0\}$ for some
$z_1\in\bbC\backslash\{z_0\}$
would contradict the fact that $\dim_\bbC(\ker(\Im(M(z))))$
is constant
for $z\in\bbC_+$. Thus $\ker(\Im(M(z)))=\{0\}$ for all
$z\in\bbC_+$
thereby completing the proof.
\end{proof}

The following result, the analog of
Theorem~\ref{t2.2}, is well-known to experts in the theory
of self-adjoint
extensions of symmetric operators and especially, in the
spectral theory of
matrix-valued Schr\"odinger operators, even-order
Hamiltonian systems, and
higher-order ordinary differential and difference
operators. For
relevant material we refer the reader, for instance, to
\cite{AG94}--\cite{AGS96}, \cite{At64}, Ch.~9, \cite{BT92},
\cite{Be68}, Sect.~VI.5, \cite{Br71}, Sect.~I.4,
\cite{BEME96}, \cite{Bu97}, \cite{Ca76},
\cite{Cl90}, \cite{DM91}--\cite{DM95}, \cite{DS88},
Sects.~XIII.5--XIII.7,
\cite{Ev64}, \cite{Fu76}, \cite{GKS97a}, \cite{GKS97b},
\cite{HS93},
\cite{HS84}, \cite{HS86}, \cite{Jo87}, \cite{Ko50},
\cite{KS88}--\cite{Kr89b}, \cite{Kr95},
Chs.~7, 8,
\cite{KO78}, \cite{MHR77}, \cite{Na68}, Ch.~VI,
\cite{Ni72},
\cite{Sa92}--\cite{Sa97}, \cite{Ts92}, \cite{We87},
Sects.~8--10.

However, since proofs are not always readily available
in the literature,
we briefly sketch some pertinent arguments which
essentially reduce the
matrix case to the scalar situation described in
Theorem~\ref{t2.2}.

\begin{theorem}\label{t5.3} Let
$M(z)\in M_n(\mathbb{C})$ be a
matrix-valued Herglotz function. Then

\noindent (i). Each diagonal element $M_{j,j}(z)$,
$1\leq
j\leq n$ of $M(z)$ is a (scalar) Herglotz function.

\noindent (ii). $M(z)$ has finite normal limits
$M(\lambda \pm
i0)=\lim_{\varepsilon \downarrow 0} M(\lambda
\pm i\varepsilon)$ for
a.e.~$\lambda \in \mathbb{R}$.

\noindent (iii). If each diagonal element $M_{j,j}(z)$,
$1\leq j\leq n$ of
$M(z)$ has a zero normal limit on a fixed subset of
$\mathbb{R}$ having
positive Lebesgue measure, then $M(z)=C_0$, where
$C_0=C_0^*$ is a
constant self-adjoint $n\times n$ matrix with vanishing
diagonal elements.

\noindent (iv). There exists a matrix-valued measure $\Omega$
on the bounded
Borel subsets of $\mathbb{R}$ satisfying
\begin{equation}\lb{5.9}
\int_{\mathbb{R}}(\ul c,d\Omega(\lambda)\ul c)_{\bbC^n}
(1+\lambda^2)^{-1}<\infty \text{ for all } \ul c\in\bbC^n
\end{equation}
such that the Nevanlinna, respectively, Riesz-Herglotz
representation
\begin{align} \lb{5.10}
M(z)&=C+Dz+\int_{\mathbb{R}}d
\Omega (\lambda )((\lambda-z)^{-1}-\lambda
(1+\lambda ^2)^{-1}),
\quad z\in\bbC_+, \\
C&=\Re(M(i)),\quad D=\lim_{\eta\uparrow \infty}
\big(\frac{1}{i\eta}M(i\eta
)\big)\geq 0 \no
\end{align}
holds.

\noindent (v). The Stieltjes inversion formula for
$\Omega $ reads
\begin{equation}\lb{5.11}
\frac{1}{2}\Omega
(\{\lambda_1\})+\frac{1}{2}\Omega (\{\lambda_2\})+\Omega
((\lambda_1,\lambda_2))=\pi^{-1}\lim_{\varepsilon\downarrow
0}\int^{\lambda_2}_{\lambda_1}d\lambda \,
\Im (M(\lambda+i\varepsilon)).
\end{equation}

\noindent (vi). The absolutely continuous part
$\Omega _{ac}$ of $\Omega$ is
given by
\begin{equation}\lb{5.12}
d\Omega_{ac}(\lambda)=\pi^{-1}\Im
(M(\lambda+i0))d\lambda .
\end{equation}

\noindent (vii). Any poles of $M$ are simple and located
on the
real axis, the residues at poles being nonpositive
matrices (of rank
$r\in \{1,\ldots ,n\})$.
\end{theorem}

\begin{proof} (i). Since for all $\ul x\in \mathbb{C}^n$,
\begin{equation}
(\ul x,M(z)\ul x)_{\mathbb{C}^n}\quad\text{is a (scalar)
Herglotz
function,}\lb{5.13}
\end{equation}
the choice $\ul x=\ul x_j=(x_{j,1},\ldots
,x_{j,n})^t\in\mathbb{C}^n$, $x_{j,\ell}=\delta_{j,\ell}$
in \eqref{5.13}
proves (i). (Here ``t'' denotes the transpose operation.)

(ii). Consider $\ul x_j=(x_{j,1},\ldots ,x_{j,n})^t\in
\mathbb{C}^n$,
$x_{j,\ell}=\delta _{j,\ell}$ and apply the polarization
identity to
$(\ul x_j,M(z)\ul x_k)$, $j\neq k$ to obtain
\begin{equation}
\begin{split} \lb{5.14}
M_{j,k}(z)&=\frac{1}{4}\big(((\ul x_j+
\ul x_k),M(z)(\ul x_j+\ul x_k))_{\mathbb{C}^n}-
((\ul x_j-\ul x_k),M(z)
(\ul x_j-\ul x_k))_{\mathbb{C}^n}\\ &\quad +
i((\ul x_j-i\ul x_k),M(z)(\ul x_j
-i\ul x_k))_{\mathbb{C}^n}
-i((\ul x_j+i\ul x_k),M(z)(\ul x_j
+i\ul x_k))_{\mathbb{C}^n}\big).
\end{split}
\end{equation}
Combining \eqref{5.13}, \eqref{5.14}, and
Theorem~\ref{t2.2}(i) then proves (ii).

(iv),(v). By \eqref{5.13} and
Theorem~\ref{t2.2}(iii),(iv) one infers
for all $\ul x\in\mathbb{C}^n$ the representation
\begin{equation}\lb{5.15}
(\ul x,M(z)\ul x)_{\mathbb{C}^n}=c_{\ul x}+d_{\ul x}z+
\int_{\mathbb{R}}d\omega_{\ul x}(\lambda)
((\lambda-z)^{-1}-\lambda (1+\lambda ^2)^{-1}),
\end{equation}
with
\begin{align}\lb{5.16}
&\int_{\mathbb{R}}d\omega_{\ul x}
(\lambda )(1+\lambda ^2)^{-1}<\infty,\\
& c_{\ul x}=\Re((\ul x,M(i)\ul x)_{\mathbb{C}^n}),\quad
d_{\ul x}=\lim_{\eta\uparrow
\infty}(\ul x,M(i\eta )\ul x)_{\bbC^n}/(i\eta )\geq 0. \no
\end{align}
In addition, for
$(\lambda_1,\lambda_2)\subset \mathbb{R}$,
\begin{equation}\lb{5.17}
\frac{1}{2}\omega_{\ul x}(\{\lambda
_1\})+\frac{1}{2}\omega_{\ul x}(\{\lambda_2\})+\omega
_{\ul x}((\lambda_1,\lambda_2)) = \pi^{-1}\lim_{\varepsilon
\downarrow
0}\int^{\lambda_2}_{\lambda_1}d\lambda \,
\Im ((\ul x,M(\lambda +i\varepsilon)\ul x)_{\mathbb{C}^n}).
\end{equation}
The polarization identity for
$(\ul x,M(z)\ul y)_{\mathbb{C}^n}$ then yields for all
$\ul x,\ul y\in \mathbb{C}^n$,
\begin{align} \lb{5.18}
(\ul x,M(z)\ul y)_{\mathbb{C}^n}= & C(\ul x,\ul y)
+D(\ul x,\ul y)z \no \\
& +\int_{\mathbb{R}}d\Omega (\ul x,\ul y)(\lambda)
((\lambda-z)^{-1}-\lambda (1+\lambda^2)^{-1}),
\end{align}
where, in obvious notation,
\begin{align}
C(\ul x,\ul y)&=\frac{1}{4}(c_{\ul x+\ul y}-c_{\ul x-\ul y}
+ic_{\ul x-i\ul y}-ic_{\ul x+i\ul y}), \lb{5.19}\\
D(\ul x,\ul y)&=\frac{1}{4}(d_{\ul x+\ul y}-d_{\ul x-\ul y}
+id_{\ul x-i\ul y}-id_{\ul x+i\ul y}),\lb{5.20}
\end{align}
\begin{align}
\frac{1}{2} \Omega (& \ul x,\ul y)  (\{\lambda_1\})
+\frac{1}{2}\Omega
(\ul x,\ul y)(\{\lambda_2\})+
\Omega (\ul x,\ul y)((\lambda _1,\lambda_2)) \lb{5.21} \\
&=\frac{1}{4}\bigg(\frac{1}{2}\omega_{\ul x
+\ul y}(\{\lambda_1\})
-\frac{1}{2}\omega_{\ul x-\ul y}(\{\lambda_1\})
+\frac{i}{2}\omega_{\ul x-i\ul y}(\{\lambda_1\})
-\frac{i}{2}\omega_{\ul x+i\ul y}(\{\lambda_1\}) \no \\
&\quad +\frac{1}{2}\omega_{\ul x+\ul y}(\{\lambda_2\})
-\frac{1}{2}\omega_{\ul x-\ul y}(\{\lambda_2\})+\frac{i}{2}
\omega_{\ul x-i\ul y}(\{\lambda_2\})
-\frac{i}{2}\omega_{\ul x+i\ul y}(\{\lambda_2\}) \no \\
&\quad +\omega_{\ul x+\ul y}((\lambda_1,\lambda_2))
-\omega_{\ul x-\ul y}((\lambda_1,\lambda_2))
+i\omega_{\ul x-i\ul y}((\lambda_1,\lambda_2))
-i\omega_{\ul x+i\ul y}((\lambda_1,\lambda_2)) \bigg).
\no
\end{align}
Since $C(\ul x,\ul y)$ and $D(\ul x,\ul y)$ are
symmetric sesquilinear forms and
$D(\ul x,\ul x)\geq 0$ for all $\ul x\in \mathbb{C}^n$,
one infers
\begin{equation}\lb{5.22}
C(\ul x,\ul y)=(\ul x,C\ul y)_{\mathbb{C}^n},\quad
D(\ul x,\ul y)=(\ul x,D\ul y)_{\mathbb{C}^n}
\end{equation}
for some
\begin{equation}\lb{5.23}
C=C^*\in M_n(\mathbb{C}),\quad 0\leq D\in
M_n(\mathbb{C}).
\end{equation}
Similarly, using the obvious fact that
$\Im ((\ul x,M(z)\ul x)_{\mathbb{C}^n})
=(\ul x,\Im (M(z))\ul x)_{\mathbb{C}^n}$, $\ul x\in
\mathbb{C}^n$,
\eqref{5.17} then becomes
\begin{equation}\lb{5.24}
(5.22)=(\ul x,\pi ^{-1}\lim_{\varepsilon
\downarrow 0}\int^{\lambda_2}_{\lambda_1}d\lambda \, \Im
(M(\lambda+i\varepsilon ))\ul y)_{\mathbb{C}^n}.
\end{equation}
Arbitrariness of
$\ul x,\ul y\in \mathbb{C}^n$ then yields the representation
\eqref{5.10} for
$M(z)$ and the Stieltjes inversion formula \eqref{5.11}. That
$C=\Re (M(i))$ and $D=\lim_{\eta \uparrow
\infty}M(i\eta )/(i\eta)$ is clear from the
corresponding properties in \eqref{5.16}.

(iii). Let $X_0\subseteq\mathbb{R}$ be the fixed subset in
(iii). Then by
hypothesis, $D=0$, $\Omega =0$ using
\eqref{5.10}, \eqref{5.11} and $\Im (M(\lambda +i0))=0$ for
$\lambda
\in X_0$ by Lemma~\ref{l5.1}(ii). Thus $M(z)=C$ is constant
with vanishing diagonal elements.

(vi). Studying $(\ul x, \Im (M(\lambda
+i\varepsilon ))\ul x)_{\mathbb{C}^n}$,
$\ul x\in \mathbb{C}^n$, one can follow the argument in
\cite{Si95b}, Theorem~1.6(iv) step by step.

(vii). First-order poles with nonpositive residues at isolated
singularities of $M(z)$ on the real axis follows from
polarization,
\eqref{5.13}, and Theorem~\ref{t2.2}(vi).\end{proof}

In the scalar case described in Theorem~\ref{t2.2}, isolated
zeros of
$m(z)$ are necessarily simple and located on
$\mathbb{R}$. This can of course be inferred from the fact
that $-1/m(z)$
is a Herglotz function whenever $m(z)$ is one (cf.~(2.10))
and hence
isolated poles of $1/m(z)$ are necessarily simple. This
reformulation
concerning isolated simple real zeros of $m(z)$ extends to
     the matrix
case since we will show later on (cf.~Theorem~\ref{t6.4}(i))
that if $M(z)$ is
invertible on
$\mathbb{C}_+$, then $-M(z)^{-1}$ is a Herglotz matrix
whenever
$M(z)$ is
one. Hence isolated poles of $M(z)^{-1}$ on
$\mathbb{R}$ are necessarily simple.

It should be remarked at this point that Theorem~\ref{t5.3}(iv)
as well
as Theorem~\ref{t5.4}(iii) below, are well-known to extend to
infinite-dimensional situations under appropriate hypotheses
on $M(z)$.
We will return to this circle of ideas elsewhere.

Due to \eqref{5.7}, Theorem~\ref{t2.3}(i)--(vi) and
Theorem~\ref{t2.4}(i) extend to the present matrix-valued
context with
only minor modifications. For later reference we summarize
a few of
these extensions below.

\begin{theorem}\label{t5.4}
Let $M(z)\in M_n (\mathbb{C})$ be a
matrix-valued Herglotz function with representation
\eqref{5.10}.
Then

\noindent (i). For all $\lambda \in \mathbb{R}$,
\begin{align}\lb{5.25}
& \lim_{\varepsilon \downarrow 0} \varepsilon \Re
(M(\lambda +i\varepsilon ))=0,\\
& \Omega (\{\lambda \})=\lim_{\varepsilon \downarrow 0}
\varepsilon \Im (M(\lambda
+i\varepsilon ))=-i\lim_{\varepsilon
\downarrow 0}\varepsilon M(\lambda +i\varepsilon ).\lb{5.26}
\end{align}

\noindent (ii). Let $L>0$ and suppose
$0\leq \Im (M(z))\leq LI_n$ for
all
$z\in \mathbb{C}_+$. Then $D=0$, $\Omega$ is purely
absolutely
continuous, $\Omega =\Omega _{ac}$, and
\begin{equation}\lb{5.27}
0\leq \frac{d\Omega (\lambda
)}{d\lambda}=\pi^{-1}\lim_{\varepsilon \downarrow 0}
\Im (M(\lambda+i\varepsilon
))\leq \pi^{-1}LI_n \text{ for a.e. }\lambda \in
\mathbb{R}.
\end{equation}

\noindent (iii). Assume $M(z)$ is invertible for all
$z\in \mathbb{C}_+$.
Then there exist $\Xi_{j,k}\in L^\infty (\mathbb{R})$,
$1\leq j,k\leq n$,
$0\leq \Xi \leq I_n$ a.e., such that
\begin{align}\lb{5.28}
\ln(M(z)) &=  K+\int_{\mathbb{R}}d\lambda \, \Xi (\lambda
)((\lambda-z)^{-1}-\lambda (1+\lambda^2)^{-1}),
\quad z\in\bbC_+, \\
 K &=\Re(\ln
(M(i))), \no
\end{align}
where
\begin{equation}
\Xi (\lambda )=\pi^{-1}\lim_{\varepsilon \downarrow 0}\Im
(\ln (M(\lambda +i\varepsilon ))) \text{ for a.e. }
\lambda\in\bbR.
\end{equation}
\end{theorem}

\noindent{\it Proof of (iii).} We briefly sketch an
approach by
Carey \cite{Ca76} (designed for the infinite-dimensional
context).
Define $\ln (z)$ with a cut along $(-\infty,0]$ such that
$\ln (\lambda )$ is real-valued for $\lambda >0$, that is,
$0<\arg (\ln
(z))<\pi$ for all $z\in \mathbb{C}_+$.
Since by hypothesis $0\notin\sigma (M(z))$, $z\in
\mathbb{C}_+$, one can define $\ln (M(z))$ for $z\in
\mathbb{C}_+$ by

\begin{equation}
\ln (M(z))=\int_{-\infty}^0 d\lambda
((\lambda -M(z))^{-1} - \lambda(1+\lambda^2)^{-1}). \lb{5.30}
\end{equation}
Next,
introducing $\ln (z;\eta)=\ln (z+i\eta )$ for $\eta > 0$,
$\ln (\cdot
;\eta)$ is analytic on
$\mathbb{C}_+$. Denoting by $W(A)$ the numerical range of $A\in
M_n(\mathbb{C})$ (i.e., $W(A)=\{(\ul x,A\ul x)\,|\,\ul x\in
\mathbb{C}^n,\|\ul x\|_{\mathbb{C}^n}=1\})$, a theorem by Kato
\cite{Ka65} relating $W(A)$ and $W(f(A))$ for analytic
functions
$f$ on closed domains conformally equivalent to $\ol D$ (the
closure of
the open unit disk $D\in\bbC$),
applied to $\ln (z;\eta )$ for $z\in \mathbb{C}_+$, yields
\begin{equation}\lb{5.32}
W(\ln (M(z);\eta ))\subset
\{\zeta\in\mathbb{C}_+ \, | \,0\leq \Im(\zeta) \leq \pi \},
\quad z\in \mathbb{C}_+,
\end{equation}
that is,
\begin{equation}\lb{5.33}
0\leq \Im (\ln (M(z);\eta))\leq \pi I_n , \quad
z\in \mathbb{C}_+.
\end{equation}
Continuity of $\ln (A;\eta )$ with
respect to $\eta$, $\lim_{\eta\downarrow 0}\ln (A;\eta )
=\ln (A)$,
for $A\in
M_n(\mathbb{C})$ nonsingular, then yields
\begin{equation} \lb{5.34}
0\leq \Im (\ln (M(z)))\leq\pi I_n,\quad
z\in \mathbb{C}_+\end{equation} and one can apply part (ii)
(as in
Theorem~\ref{t2.3}).
\hfill$\square$ \\

Finally we state the matrix analogs of Lemmas~\ref{l2.5} and
\ref{l2.6}, the proofs of which we omit since they are
essentially
identical to the scalar case.

\begin{lemma} \lb{l5.5}
Let $M$ be a Herglotz matrix with representation \eqref{5.10}
and
$(\lambda_1,\lambda_2)\subseteq\bbR$, $\lambda_1<\lambda_2$.
Then $M$ can
be analytically continued from
$\bbC_+$ into a subset of $\bbC_-$ through the interval
$(\lambda_1,\lambda_2)$
if and only if the
associated measure $\Omega$ is purely absolutely continuous on
$(\lambda_1,\lambda_2)$, $\omega {\big|_{(\lambda_1,\lambda_2)}}
=\Omega {\big|_{(\lambda_1,\lambda_2),ac}}$, and the density
$\Omega'\geq 0$ of $\Omega$ is real-analytic on
$(\lambda_1,\lambda_2)$.
In this case, the
analytic continuation of $M$ into some domain
$\calD_-\subseteq\bbC_-$
is given by
\begin{equation}
M(z)=M(\overline z)^* + 2\pi i \Omega'(z),
\quad z\in\calD_-, \lb{5.35}
\end{equation}
where $\Omega'(z)$ denotes the complex-analytic extension of
$\Omega'(\lambda)$ for $\lambda\in(\lambda_1,\lambda_2)$. In
particular, $M$ can be analytically continued through
$(\lambda_1,\lambda_2)$ by reflection, that is,
$M(z)=M(\overline z)^*$ for all $z\in\bbC_-$ if and only
if $\Omega$ has no support in $(\lambda_1,\lambda_2)$.
\end{lemma}

\begin{lemma}  [\cite{KS88}] \lb{l.5.6}
Let $M$ be a Herglotz matrix and $(\lambda_1,\lambda_2)
\subseteq\bbR$,
$\lambda_1<\lambda_2$. Suppose $\lim_{\varepsilon\to 0}
\Re(M(\lambda+i\varepsilon))=0$ for a.e. $\lambda
\in(\lambda_1,\lambda_2)$.
Then $M$ can be analytically continued from $\bbC_+$ into
$\bbC_-$ through the interval $(\lambda_1,\lambda_2)$ and
\begin{equation}
M(z)=-M(\overline z)^*. \lb{5.36}
\end{equation}
In addition, $\Im(M(\lambda+i0))>0$, $\Re(M(\lambda+i0))=0$
for all $\lambda\in(\lambda_1,\lambda_2)$.
\end{lemma}

\section{Support Theorems in the Matrix Case} \lb{s6}

The principal aim of this section is to prove a support
theorem for
$\Omega_{ac}$ in connection with the matrix analog of the
Aronszajn-Donoghue theory (cf.~Theorem~\ref{t3.2}(i),(ii)).

Supports $S_\Omega$, topological supports $S_\Omega^{cl}$,
and minimal
supports (with respect to Lebesgue measure on $\mathbb{R}$)
of matrix-valued
measures such as $\Omega$ in \eqref{5.4}, \eqref{5.4} are
defined as
in the beginning of Section~\ref{s3}. Because of
\eqref{5.6}, in
discussing supports of the matrix measure $\Omega$, we will
occasionally
replace $\Omega$ by the (scalar) trace measure
$\omega^{tr}=\tr_{\mathbb{C}^n}(\Omega)$. For pure point
measures,
$\Omega=\Omega_{pp}$, we again consider the smallest
support. If a
pure point measure $\Omega=\Omega_{pp}$ contains no
finite accumulation
points
in its support we call it a discrete point measure and
denote it by
$\Omega_d$.

In order to capture spectral multiplicities in the
matrix-valued case in
connection with applications to differential and
difference operators we
introduce the sets ($1\leq r\leq n$)
\begin{align}
S_{\Omega_{ac},r}&=\{\lambda \in
\mathbb{R} \, | \lim_{\varepsilon \downarrow 0}M(\lambda +
i \varepsilon)
\text{ exists finitely}, \, \rank(\Im (M(\lambda +i0)))=r\},
\lb{6.1} \\
S_{\Omega_{ac}}&=\bigcup_{r=1}^n S_{\Omega_{ac},r},\lb{6.2}\\
S_{\Omega_{pp},r}&=\{\lambda \in \mathbb{R} \, | \,
\rank (\lim_{\varepsilon
\downarrow 0}\varepsilon M(\lambda+i\varepsilon ))=r\},
\quad 1\leq r\leq
n,\lb{6.3}\\
S_{\Omega_{pp}}&=\bigcup_{r=1}^n S_{\Omega_{pp},r},\lb{6.4}\\
S_{\Omega_{s}}&=\{\lambda \in \mathbb{R} \, | \lim_{\varepsilon
\downarrow 0}\Im
(\tr_{\bbC^n}(M(\lambda +i\varepsilon)))=+\infty \},\lb{6.5}\\
S_{\Omega_{sc}}&=\{\lambda \in S_{\Omega_{s}} \, |
\lim_{\varepsilon \downarrow 0}\varepsilon
\tr_{\bbC^n}(M(\lambda +i\varepsilon))=0\},\lb{6.5a}\\
S_\Omega&=S_{\Omega_{ac}}\cup S_{\Omega_{s}}. \lb{6.6}
\end{align}
(Here existence of matrix limits are of course understood for
each individual matrix element.)
Thus, $S_{\Omega_{ac},r}$,
$S_{\Omega_{pp},r'}$, $S_{\Omega_{sc}}$ are all disjoint
for any $1\leq
r,r'\leq n$.

As in \eqref{3.9} we define the equivalence classes
$\mathcal{E}(\Omega_{ac})$ and
$\mathcal{E}_r(\Omega_{ac})$ of
$S_{\Omega_{ac}}$ and $S_{\Omega_{ac},r}$, $1\leq r\leq n$
with respect
to the equivalence relation \eqref{3.1} (with
$\nu$ representing Lebesgue measure on $\mathbb{R}$ and
$\mu=\Omega_{ac}$).

The following result is analogous to Theorem~\ref{t3.1} in
the scalar case and can be reduced to it by studying the
trace measure
$\omega^{tr}$ of $\Omega$.

\begin{theorem}\label{t6.1} Let $M$ be a matrix-valued
Herglotz function
with representations \eqref{5.10} and
\eqref{5.28}. Then

\noindent (i). $S_{\Omega_{ac}}$ is a minimal support of
$\Omega_{ac}$.

\noindent (ii). $S_{\Omega_{sc}}$ is a minimal support
of $\Omega_{sc}$.

\noindent (iii). $S_{\Omega_{pp}}$ is the smallest support
of $\Omega_{pp}$.

\noindent (iv). $S_\Omega$ is a minimal support of $\Omega$.

\noindent (v). If in addition $M(z)$ is invertible for
all $z\in\bbC_+$,
then
\begin{equation}\lb{6.7}
\hspace*{+.4cm} \widetilde{S}_{ac}=\{\lambda
\in S_{\Omega_{ac}}
\,| \, \ln(M(\lambda+i0)) \text{ exists finitely and }
0<\tr(\Xi (\lambda)) < n\}
\end{equation}
is a minimal support of
$\Omega_{ac}$.
\end{theorem}

\noindent{\it Proof of (v).} By definition,
$\widetilde S_{ac}\backslash
S_{\Omega_{ac}} =\emptyset$. Next, suppose
$\tr(\Xi(\lambda))$ equals $0$ or $n$.
Then one concludes from $0\leq\Xi(\lambda)\leq I_n$ for
a.e.~$\lambda\in\bbR$ (cf.~Theorem~5.4\,(iii)) that
$\Im(\ln(M(\lambda+i0)))=0$ or $\Im(\ln(-M(\lambda+i0)))=0$,
that is,
$\ln(M(\lambda+i0))$ or $\ln(-M(\lambda+i0))$ is
self-adjoint. Taking
exponentials, $M(\lambda+i0)$ is self-adjoint and hence
\begin{equation}
\Im(M(\lambda+i0))=0 \text{ for }
\lambda\in \{\nu\in\bbR \, |
\tr(\Xi(\nu))\in \{0,n\} \}. \lb{6.7a}
\end{equation}
Thus, abbreviating Lebesgue measure on $\bbR$ by
$|\cdot|$, one infers
\begin{align}
& |S_{\Omega_{ac}}\backslash \widetilde S_{ac}| =
|\{\lambda\in S_{\Omega_{ac}} \, | \,
\text{either } \Im(\ln(M(\lambda+i0))_{j,j})
\text{ does not exist finitely } \no \\
& \hspace*{2.35cm} \text{for
some } 1\leq j\leq n, \text{ or }
\tr(\Xi(\lambda))\in\{0,n\}\}| \no \\
& =|\{\lambda\in S_{\Omega_{ac}} \, | \,
\Im(\ln(M(\lambda+i0))_{j,j})
\text{ does not exist finitely for some }
1\leq j\leq n\}| \no \\
& =0 \lb{6.7b}
\end{align}
by \eqref{6.7a}, the fact that $\lambda\in
S_{\Omega_{ac}}$ implies
$\Im(M(\lambda+i0)) >0$, and Theorem~\ref{t5.3}\,(ii)
(applied to
$\ln(M(z))$). Thus, $|\widetilde S_{ac}
\triangle S_{\Omega_{ac}}|=0$ and since
$\Omega_{ac}$ is absolutely continuous with respect to
Lebesgue measure
$|\cdot|$, also $\Omega_{ac}(\widetilde S_{ac}\triangle
S_{\Omega_{ac}})=0$.
Consequently, $\widetilde S_{ac}$ and $S_{\Omega_{ac}}$
are equivalent
minimal
supports for $\Omega_{ac}$, $\widetilde S_{ac}
\sim S_{\Omega_{ac}}$.
\hfill$\square$ \\

In order to prove the analog of Theorem~\ref{t3.2}(i) in the
matrix-valued case, that is, the stability of the minimal
support
$S_{\Omega_{ac}}$ with respect to linear fractional
transformations
(generalizing \eqref{2.10} to the matrix case as in
\eqref{6.19}), we
need to introduce a bit of preparatory material.

Define
\begin{equation}\lb{6.8}
J_{2n}=\left(\begin{array}{cc}0 & -I_n\\ I_n &
0\end{array}\right),
\end{equation}
\begin{equation}
\mathcal{A}_{2n}=\{ A\in
M_{2n}(\mathbb{C}) | \, A^*J_{2n}A=J_{2n}\}.\lb{6.9}
\end{equation}
Representing $A\in M_{2n}(\mathbb{C})$ by
\begin{equation}\lb{6.10}
A=\left(\begin{array}{cc}A_{1,1} &
A_{1,2}\\ A_{2,1} & A_{2,2}\end{array}\right),
\quad A_{p,q}\in
M_n(\mathbb{C}),\quad 1\leq p,q\leq 2,
\end{equation}
the condition
$A^*J_{2n}A=J_{2n}$ in \eqref{6.9} explicitly reads
\begin{align}\lb{6.11}
& A^*_{1,1}A_{2,1}=A^*_{2,1}A_{1,1},\quad
A^*_{2,2}A_{1,2}=A^*_{1,2}A_{2,2}, \no \\
& A^*_{2,2}A_{1,1}-A^*_{1,2}A_{2,1}=I_n=A^*_{1,1}A_{2,2}
-A^*_{2,1}A_{1,2},
\end{align}
or equivalently,
\begin{equation}\lb{6.12}
\left(\begin{array}{cc}A^*_{2,2}
&-A^*_{1,2}\\-A^*_{2,1}&
A^*_{1,1}\\\end{array}\right)\left(\begin{array}{cc}A_{1,1}
& A_{1,2}\\
A_{2,1} & A_{2,2}\end{array}\right)=I_{2n}.
\end{equation}
Since left
inverses in $M_{2n}(\mathbb{C})$ are also right
inverses, \eqref{6.12}
implies
\begin{equation}\lb{6.13}
\left(\begin{array}{cc} A_{1,1} &
A_{1,2}\\ A_{2,1} & A_{2,2}\end{array}\right)
\left(\begin{array}{cc}
A^*_{2,2} & -A^*_{1,2}\\ -A^*_{2,1} &
A^*_{1,1}\end{array}\right)=I_{2n},
\end{equation}
that is,
\begin{align} \lb{6.14}
& A_{1,1}A^*_{1,2}=A_{1,2}
A^*_{1,1},\quad A_{2,2}A^*_{2,1}=A_{2,1}A^*_{2,2},
\no \\
& A_{2,2}A^*_{1,1}-A_{2,1}A^*_{1,2}=I_{n}
=A_{1,1}A^*_{2,2}-A_{1,2}
A^*_{2,1},
\end{align}
or equivalently,
\begin{equation}\lb{6.15}
AJ_{2n}A^*=J_{2n}.
\end{equation}
In particular,
\begin{equation}\lb{6.16}
A\in \mathcal{A}_{2n}\text{ if and
only if } A^{-1}\in \mathcal{A}_{2n}.
\end{equation}
Next, let
$A=(A_{p,q})_{1\leq p,q\leq 2}\in \mathcal{A}_{2n}$ and
suppose $M\in
M_n(\mathbb{C})$ is chosen such that $\ker (A_{1,1}
+A_{1,2}M)=\{0\}$,
that is $(A_{1,1}+A_{1,2}M)$ is invertible in
$\mathbb{C}^n$. Define
(cf., e.g.,
\cite{KS74})
\begin{equation}\lb{6.17}
M_A(M)=(A_{2,1}+A_{2,2}M)(A_{1,1}+A_{1,2}M)^{-1}
\end{equation}
to observe
\begin{align}
M_{I_{2n}} (M)&=M, \no\\
M_A(M_B(M))&=M_{AB}(M), \lb{6.18} \\
M_A(M)&=M_{AB^{-1}}(M_B(M)), \no
\end{align}
whenever $M_A(M)$ and $M_B(M)$ exist.

We are particularly interested in the case where $M$
in \eqref{6.17}
equals an $n\times n$ Herglotz matrix $M(z)$. In this
case the existence
of $(A_{1,1}+A_{1,2}M(z_0))^{-1}$ for some $z_0\in
\mathbb{C}_+$ and
analyticity of $M(z)$ on
$\mathbb{C}_+$ proves that $(A_{1,1}+A_{1,2}M(z))^{-1}$
is meromorphic on
$\mathbb{C}_+$. However, since later on we are interested
in analyticity
of $M_A(M(z))$ for all $z\in \mathbb{C}_+$, we will usually
assume that
$\ker(A_{1,1}+A_{1,2}M(z))=\{0\}$ for all
$z\in \mathbb{C}_+$. Moreover,
in a slight abuse of notation, we shall abbreviate
$M_A(M(z))$ by
\begin{equation} \lb{6.19}
M_A(z)=(A_{2,1}+A_{2,2}M(z))(A_{1,1}+
A_{1,2}M(z))^{-1},\quad A\in\mathcal{A}_{2n},\quad z\in
\mathbb{C}_+
\end{equation}
from now on.

We start with a series of results concerning linear
fractional
transformations of the type \eqref{6.19}.

\begin{lemma}\label{l6.2} Suppose $A=(A_{p,q})_{1\leq p,q
\leq 2}\in
\mathcal{A}_{2n}$ and $M\in M_n(\mathbb{C})$ with
$\Im (M)\geq 0$. Then
\begin{equation}\lb{6.20}
\ker (A_{1,1}+A_{1,2}M)\subseteq \ker (\Im
(M)).
\end{equation}
In particular,
\begin{equation}\lb{6.21}
\Im (M)>0 \text{ implies } \ker
(A_{1,1}+A_{1,2}M)=\{0\}.
\end{equation}
\end{lemma}

\begin{proof} Suppose the existence of an $\ul x_0\in
\mathbb{C}^n\backslash\{0\}$ such that
\begin{equation}\lb{6.22}
(A_{1,1}+A_{1,2}M)\ul x_0=0.
\end{equation}
Then
\begin{align}\lb{6.23}
&(\ul x_0,\Im(M)\ul x_0)_{\mathbb{C}^n}
=(2i)^{-1}\big((\ul x_0,M\ul x_0)_{\mathbb{C}^n}-
(M\ul x_0,\ul x_0)_{\mathbb{C}^n}\big) \no \\
&=(2i)^{-1}\big((\ul x_0,(A_{1,1}^*A_{2,2}
-A^*_{2,1}A_{1,2})M\ul x_0)_{\mathbb{C}^n}-((A^*_{1,1}A_{2,2}
-A_{2,1}^*A_{1,2})M\ul x_0,
\ul x_0)_{\mathbb{C}^n}\big) \no \\
&=(2i)^{-1}\big((A_{1,1}
\ul x_0,A_{2,2}M\ul x_0)_{\mathbb{C}^n}
+(\ul x_0,A^*_{2,1}A_{1,1}\ul x_0)_{\mathbb{C}^n}
-(A_{2,2}M\ul x_0,A_{1,1}\ul x_0)_{\mathbb{C}^n} \no \\
& \hspace*{1.6cm}
-(\ul x_0,A^*_{1,1}A_{2,1}\ul x_0)_{\mathbb{C}^n}\big) \no \\
&=(2i)^{-1}\big(-(A_{1,2}M\ul
x_0,A_{2,2}M\ul x_0)_{\mathbb{C}^n}
+(A_{2,2}M\ul x_0,A_{1,2}M\ul x_0)_{\mathbb{C}^n}\big) \no \\
&=(2i)^{-1}(M\ul x_0,(A^*_{2,2}A_{1,2}-
A^*_{1,2}A_{2,2})M\ul x_0)_{\mathbb{C}^n}=0,
\end{align}
where we repeatedly used \eqref{6.11} and \eqref{6.22}.
Since $\Im
(M)\geq 0$ by hypothesis, \eqref{6.23} yields
$\Im (M)\ul x_0=0$
and hence \eqref{6.20}.
\end{proof}

\begin{lemma}\label{l6.3}
Suppose $A=(A_{p,q})_{1\leq p,q \leq
2}\in \mathcal{A}_{2n}$, $M\in M_n(\mathbb{C})$, and
$\ker (A_{1,1}+A_{1,2}M)=\{0\}$. Define
\begin{equation}\lb{6.24}
M_A=(A_{2,1}+A_{2,2}M)(A_{1,1}+
A_{1,2}M)^{-1}.
\end{equation}
Then

\noindent (i). \begin{equation}\lb{6.25}
\Im
(M_A)=((A_{1,1}+A_{1,2}M)^{-1})^*\Im
(M)(A_{1,1}+A_{1,2}M)^{-1}.
\end{equation}

\noindent (ii). \begin{align}&
(A^*_{2,2}-A^*_{1,2}M_A)(A_{1,1}+A_{1,2}M)=I_n,\lb{6.26}
\\ &\ker
(A^*_{2,2}-A^*_{1,2}M_A)=\{0\}.\lb{6.27}\end{align}

\noindent (iii). \begin{align}\lb{6.28}&
M=-(A^*_{2,1}-A^*_{1,1}M_A)(A^*_{2,2}-A^*_{1,2}M_A)^{-1},
\\ &\Im
(M)=((A^*_{2,2}-A^*_{1,2}M_A)^{-1})^*\Im
(M_A)(A^*_{2,2}-A^*_{1,2}M_A)^{-1}.\lb{6.29}
\end{align}
\end{lemma}

\begin{proof} (i) is a straightforward consequence of
\eqref{6.24} and
\eqref{6.11}. \eqref{6.26} is a simple consequence of
\eqref{6.24}
and \eqref{6.27} follows from \eqref{6.26}. (iii) is
readily derived
from \eqref{6.11} and \eqref{6.24}.
\end{proof}

Applying Lemmas~\ref{l6.2} and \ref{l6.3} to $M_A(z)$ in
\eqref{6.19}
then yields the following result.

\begin{theorem}\label{t6.4}
Assume $A=(A_{p,q })_{1\leq p,q \leq
2}\in \mathcal{A}_{2n}$, let $M(z)$ be an $n\times n$
Herglotz matrix,
and suppose $\ker (A_{1,1}+A_{1,2}M(z))=\{0\}$ for all $z\in
\mathbb{C}_+$. Define $M_A(z)$, $z\in
\mathbb{C}_+$ as in \eqref{6.19}. Then

\noindent (i). $M_A(z)$ is an $n\times n$ Herglotz matrix and
\begin{align} \lb{6.30}
\Im(M_A(z))&=((A_{1,1}+A_{1,2}M(z))^{-1})^*\Im (M(z))(A_{1,1}
+A_{1,2}M(z))^{-1}\geq 0, \\
& \hspace*{8.5cm} z\in\mathbb{C}_+. \no
\end{align}

\noindent (ii). For all $z\in \mathbb{C}_+$,
\begin{align}\lb{6.31}
&(A^*_{2,2}-A^*_{1,2}M_A(z))(A_{1,1}+
A_{1,2}M(z))=I_n,\\ \lb{6.32}
&\ker (A^*_{2,2}-A_{1,2}^*M_A(z))=\{0\},\\ \lb{6.33}
&M(z)=-(A^*_{2,1}-A_{1,1}^*M_A(z))(A^*_{2,2}-
A^*_{1,2}M_A(z))^{-1},\\
\lb{6.34}&\Im (M(z))=((A^*_{2,2}-A^*_{1,2}M_A(z))^{-1})^*\Im
(M_A(z))(A^*_{2,2}-A^*_{1,2}M_A(z))^{-1}.
\end{align}
\end{theorem}

\begin{proof} Assertions \eqref{6.30}--\eqref{6.34} are
clear from
Lemmas~\ref{l6.2} and \ref{l6.3}. Since
$M(z)^*=M(\overline z)$ clearly
implies $M_A(z)^*=M_A(\overline z)$ by \eqref{6.11},
$M_A(z)$ is an $n\times n$
Herglotz matrix.\end{proof}

We note in connection with Lemma~\ref{l5.2a} and
Theorem~\ref{t6.4} that
\begin{align}
&\Im (M(z_0))>0 \text{ for some } z_0\in \mathbb{C}_+ \no \\
&\text{implies }\ker (A_{1,1}+A_{1,2}M(z))=\{0\}
\text{ for all } z\in
\mathbb{C}_+ \lb{6.35}
\end{align}
by \eqref{6.21}.

\begin{remark}\label{r6.5} The condition $A\in
\mathcal{A}_{2n}$ in the
definition \eqref{6.19} of $M_A(z)$ for $M_A(z)$ to be
a Herglotz
matrix (assuming $M(z)$ to be a Herglotz matrix) is
not a necessary one.
As discussed in detail by Krein and Shmulyan \cite{KS74},
the condition $A\in
\mathcal{A}_{2n}$ in Theorem~\ref{t6.4} can be replaced
by the pair of
conditions
\begin{equation} \lb{6.36}
iA^*J_{2n}A\geq icJ_{2n},\quad iAJ_{2n}A^*\geq icJ_{2n}
\end{equation}
for some $c>0$. In a sense, by using the
condition $A\in \mathcal{A}_{2n}$, we chose equality in
\eqref{6.36}
(and $c=1$). From the point of view of applications
of matrix Herglotz functions to spectral theory of
matrix Schr\"odinger
and Jacobi operators and more generally, even-order
Hamiltonian systems,
with various boundary conditions involved, our
restrictive hypothesis
\eqref{6.9} is sufficiently general to cover all such
cases. Pertinent
references to spectral theory for even-order Hamiltonian
systems are, for
instance, \cite{At64}, Ch. 9, \cite{Cl90},
\cite{HS93}--\cite{HS86}, \cite{KS88}--\cite{Kr89b},
\cite{Kr95}, Chs. 7, 8, \cite{MHR77},
\cite{Ni72}, \cite{Sa92}--\cite{Sa94b} and
the literature cited
therein.\end{remark}

Finally, we turn to $\mathcal{E}(\Omega_{A,ac})$ and
$\mathcal{E}_r(\Omega_{A,ac})$ the equivalence classes of
$S_{\Omega_{A,ac}}$ and $S_{\Omega_{A,ac,r}}$,
$1\leq r\leq n$ associated
with $M_A(z)$, $A\in \mathcal{A}_{2n}$ (cf.~(3.1) and
the paragraph
following Theorem~\ref{t6.1}). We recall that
$\mathcal{E}(\Omega_{ac})$ and $\mathcal{E}_r(\Omega_{ac})$
are the
corresponding equivalence classes of
$S_{\Omega_{ac}}$ and $S_{\Omega_{ac},r}$,
$1\leq r\leq n$ associated
with $M(z)$ (cf.~\eqref{6.1}, \eqref{6.2}). We also
recall
(cf. \eqref{6.19})
\begin{equation} \lb{6.36a}
M_A(z)=(A_{2,1}+A_{2,2}M(z))(A_{1,1}+
A_{1,2}M(z))^{-1},\quad A\in\mathcal{A}_{2n},\,\, z\in
\mathbb{C}_+
\end{equation}
and its general version
\begin{align}
&M_A(z)=((AB^{-1})_{2,1}+(AB^{-1})_{2,2}M_B(z))((AB^{-1})_{1,1}
+(AB^{-1})_{1,2}M_B(z))^{-1}, \no \\
& \hspace*{8cm} A,B\in\calA_{2n}, \,\, z\in\bbC_+. \lb{6.54}
\end{align}

Our principal result on
the absolutely continuous part of $\Omega_A$, the matrix
analog of
Theorem~\ref{t3.2}\,(i)--(iii), then reads as follows.

\begin{theorem}\label{t6.6}
Let $M(z)$ and $M_A(z)$,
$A=(A_{p,q})_{1\leq p,q \leq 2}\in \mathcal{A}_{2n}$ be
Herglotz
matrices related by \eqref{6.36a} assuming
$\ker(A_{1,1}+A_{1,2}M(z))=\{0\}$ for all $z\in
\mathbb{C}_+$. Let
$\Omega $ and $\Omega_A$ be the measures associated
with $M(z)$ and
$M_A(z)$, respectively. Then

\noindent (i). For all $A\in \mathcal{A}_{2n}$,
\begin{align}\lb{6.37}\mathcal{E}_r(\Omega_{A,ac})&=
\mathcal{E}_r(\Omega _{ac}),\quad 1\leq r\leq n,\\
\lb{6.38}\mathcal{E}(\Omega_{A,ac})&
=\mathcal{E}(\Omega_{ac}),
\end{align} that is, $\mathcal{E}_r(\Omega_{A,ac})$,
$1\leq r\leq n$ and
$\mathcal{E}(\Omega _{A,ac})$ are independent of $A\in
\mathcal{A}_{2n}$ (and hence denoted by
$\mathcal{E}_{ac,r}$, $1\leq
r\leq n$ and $\mathcal{E}_{ac}$ below) and
$\Omega_{A,ac}\sim \Omega_{ac}$
for all $A\in\calA_{2n}$.


\noindent (ii). Suppose $\Omega_B$ is a discrete point
measure,
$\Omega_B=\Omega_{B,d}$, for some $B\in \mathcal{A}_{2n}$.
Then
$\Omega_A=\Omega_{A,d}$ is a discrete point measure for
all $A\in
\mathcal{A}_{2n}$.


\noindent (iii). Define
\begin{equation}\lb{6.39}
S=\{\lambda \in
\mathbb{R} \, | \, \text{there is no } A\in
\mathcal{A}_{2n} \text{ for which } \Im (M_A(\lambda
+i0)) \text{ exists and equals } 0\}.
\end{equation}
Then $S\in \mathcal{E}_{ac}$.
\end{theorem}

\begin{proof} (i). Define
\begin{equation}
\widehat{S}_{A,r}=S_{\Omega_{A,ac,r}}\cap
\{\lambda \in \mathbb{R} \, | \, M(\lambda +i0)
\text{ exists finitely}\}.
\end{equation}
Then Theorem~\ref{t5.3}(ii) yields
\begin{equation}\lb{6.41}
|S_{\Omega_{A,ac,r}}\backslash
\widehat{S}_{A,r}|=0,
\end{equation}
where $|\cdot|$ abbreviates Lebesgue
measure on $\mathbb{R}$. Since by \eqref{6.31},
\begin{equation}\lb{6.42}
(A_{1,1}+A_{1,2}M(\lambda
+i0))^{-1}=A^*_{2,2}-A_{1,2}^*M_A(\lambda +i0)
\text{ exists
for } \lambda \in \widehat{S}_{A,r},
\end{equation}
$(A_{1,1}+A_{1,2}M(\lambda +i0))^{-1}:\mathbb{C}^n\to
\mathbb{C}^n$ is a
bijection for $\lambda \in \widehat{S}_{A,r}$ and
\eqref{6.30} yields
\begin{align}\lb{6.43}
&\Im (M_A( \lambda +i0)) \\
&=((A_{1,1}+A_{1,2}M(\lambda +i0))^{-1})^*\Im (M(\lambda
+i0))(A_{1,1} +A_{1,2}M(\lambda +i0))^{-1}, \no \\
& \hspace*{9.5cm} \lambda \in\widehat{S}_{A,r} \no
\end{align}
and hence
\begin{equation}\lb{6.44}
\text{rank}(\Im (M_A(\lambda +i0)))=\rank
(\Im (M(\lambda +i0))),\quad \lambda \in
\widehat{S}_{A,r}.
\end{equation}
Thus,
\begin{equation}\lb{6.45}
\widehat{S}_{A,r}\subseteq
S_{\Omega_{ac},r}.
\end{equation}
Similarly, defining
\begin{equation}
\widehat{S}_r=S_{\Omega_{ac},r}\cap\{\lambda
\in \mathbb{R} \, | \, M_A(\lambda +i0)
\text{ exists finitely\}},
\end{equation}
then
\begin{equation}\lb{6.47}
|S_{\Omega_{ac},r}\backslash
\widehat{S}_r|=0.
\end{equation}
By \eqref{6.31} we conclude the existence of
\begin{equation}\lb{6.48}
(A^*_{2,2}-A^*_{1,2}M_A(\lambda
+i0))^{-1}=(A_{1,1}+A_{1,2}M(\lambda+i0)),\quad \lambda \in
\widehat{S}_r
\end{equation}
and hence $(A^*_{2,2}-A^*_{1,2}M_A(\lambda
+i0))^{-1}:\mathbb{C}^n\to \mathbb{C}^n$ is a bijection
for $\lambda\in
\widehat{S}_r$. Thus \eqref{6.34} yields
\begin{align}\lb{6.49}
&\Im (M(\lambda+i0)) \\
&=((A^*_{2,2}-A^*_{1,2}M_A(\lambda +i0))^{-1})^*
\Im (M_A(\lambda
+i0))(A^*_{2,2}-A^*_{1,2}M_A(\lambda +i0))^{-1}, \no \\
& \hspace*{10.5cm} \lambda \in\widehat{S}_r \no
\end{align}
and consequently,
\begin{equation}\lb{6.50}
\text{rank}(\Im
(M(\lambda+i0)))=\text{rank} (\Im (M_A(\lambda +i0))),
\quad \lambda \in
\widehat{S}_r.
\end{equation}
Thus,
\begin{equation}\lb{6.51}
\widehat{S}_r\subseteq
S_{\Omega_{A,ac},r}.
\end{equation}
By \eqref{6.41}, \eqref{6.45},
\eqref{6.47}, and \eqref{6.51},
\begin{equation}\lb{6.52}
|S_{\Omega_{A,ac},r}\triangle
S_{\Omega_{ac},r}|=0.
\end{equation}
Since $\Omega_{A,ac}$ and
$\Omega_{ac}$ are absolutely continuous with respect to
Lebesgue measure
on $\mathbb{R}$,
\eqref{6.52} yields
\begin{equation}\lb{6.53}
\Omega_{A,ac}(S_{\Omega_{A,ac},r}\triangle
S_{\Omega_{ac},r})=\Omega_{ac}(S_{\Omega_{A,ac},r}\triangle
S_{\Omega_{ac},r})=0
\end{equation}
proving \eqref{6.37}. Equality
\eqref{6.38} is then obvious from \eqref{6.37} and \eqref{6.2}.

(ii). Part (ii) follows from \eqref{6.54} (cf.~\eqref{6.18})
and
the fact that $\Omega_A=\Omega_{A,d}$ if and only if
$M_A(z)$ is meromorphic on $\mathbb{C}$.

(iii). We follow the proof of Corollary 1 in \cite{GP87}
in the scalar case
$n=1$. Since by hypothesis $\Im(M(\lambda +i0))>0$ for all
$\lambda\in S_{\Omega_{ac}}$ one concludes from \eqref{6.20}
that
$\ker(A_{1,1}+A_{1,2}M(\lambda +i0))
=\ker(\Im(M(\lambda +i0)))=\{0\}$ and
hence $\Im(M_A(\lambda +i0))>0$ for a.e.
$\lambda\in S_{\Omega_{ac}}$ and
all $A\in\calA_{2n}$ by \eqref{6.30}. Thus, one computes
\begin{align}
|S_{\Omega_{ac}} \backslash S|  =  |\{\lambda\in
S_{\Omega_{ac}}
\, | \, \text{there is an } A_{\lambda}\in\calA_{2n}
\text{ s.t. }
\Im(M_{A_{\lambda}}(\lambda +i0))=0 \}| =0.  \lb{6.55}
\end{align}
Similarly,
\begin{align}
|S \backslash S_{\Omega_{ac}}| & = |\{\lambda\in\bbR \, | \,
\text{there is no }
A\in\calA_{2n} \text{ s.t. } \Im(M_A(\lambda +i0))=0
\text{ exists } \no \\
& \hspace*{.59cm} \text{ and either } M(\lambda +i0)
\text{ does not exist, or } M(\lambda +i0)_{j,j}
\text{ exists } \no \\
& \hspace*{.59cm} \text{ and equals } \infty
\text{ for some } 1\leq j \leq n \}| \no \\
& \leq |\{\lambda\in\bbR \, | \, \text{either } M(\lambda +i0)
\text{ does not exist, or } M(\lambda +i0)_{j,j}
\text{ exists } \no \\
& \hspace{.59cm}
\text{ and equals } \infty
\text{ for some } 1\leq j \leq n \}| =0 \lb{6.56}
\end{align}
by Theorem~\ref{t5.3}\,(ii). Thus $|S_{\Omega_{ac}}
\triangle S|=0$. Since
$\Omega_{ac} \ll |\cdot|$, one infers
$\Omega_{ac}(S_{\Omega_{ac}} \triangle S)$ $=0$ and hence
$S\sim S_{\Omega_{ac}}$, or equivalently, $S\in\calE_{ac}$.
\end{proof}

\begin{remark}\label{r6.7} One might ask whether the
first part of
Theorem~\ref{t3.2}\,(iv) extends to the matrix-valued
situation. However,
the simple counter example
\begin{equation}
M_1(z)=\left(\begin{array}{cc}m(z) & 0\\ 0 &
-m(z)^{-1}\end{array}\right), \quad M_2(z)=-M_1(z)^{-1}=
\left(\begin{array}{cc}-m(z)^{-1} & 0\\ 0 &
m(z)\end{array}\right), \lb{6.57}
\end{equation}
with $m(z)$ a scalar Herglotz function with representation
\eqref{2.14}
and $\omega_{pp}\neq 0$ or $\omega_{sc}\neq 0$, immediately
destroys such
hopes since the measures $\Omega_1$ and $\Omega_2$
corresponding to
$M_1(z)$ and $M_2(z)$ are clearly equivalent.
\end{remark}

\begin{remark}\label{r6.8} Theorem~\ref{t6.6}(i) is quite
familiar in the
context of finite-rank perturbations of the resolvent of a
self-adjoint
operator in a (complex, separable) Hilbert space. For
instance, the
absolutely continuous ({\it ac}) parts of self-adjoint
extensions of a densely
defined symmetric operator with deficiency indices $(n,n)$
are all
unitarily equivalent. In particular, their absolutely
continuous spectra
and the multiplicity functions (associated with the
{\it ac} spectra)
coincide, which is essentially \eqref{6.37} and
\eqref{6.38}.
However, even-order Hamiltonian systems do not necessarily
have such an
underlying Hilbert space formulation (cf., e.g.,
\cite{HS93}--\cite{HS86},
\cite{Kr89a}, \cite{Kr89b}, \cite{Ni72} and the
literature cited therein)
and in these cases Theorem~\ref{t6.6}(i) appears to be
an ideal tool for
identifying {\it ac} spectra associated with $\Omega_A$.
\end{remark}

As in the scalar case, the relationship between
$\Im(\ln(M_A(z)))$
(respectively, $\Xi_A(\lambda)$) and $\Im(\ln(M(z)))$
(respectively, $\Xi(\lambda)$),
analogous to \eqref{6.36a}, in general, is quite
involved. The special
case $A=J_{2n}$, that is,
\begin{equation}
M_A(z)=-M(z)^{-1}, \lb{6.58}
\end{equation}
however, is particularly simple and leads to the
analog of \eqref{3.19},
\begin{equation}
\Xi_{J_{2n}}(\lambda)=I_n-\Xi(\lambda)
\text{ for a.e. } \lambda\in\bbR.
\lb{6.59}
\end{equation}

The analog of Lemma~\ref{l3.3} in the present matrix-valued
context
appears to be more involved.

\section{Applications of Matrix-Valued Herglotz
Functions} \lb{s7}

In this section we extend some of the applications of
scalar Herglotz
functions in Section~\ref{s4} to the matrix-valued
context. In
particular, we will study self-adjoint finite-rank
perturbations of
self-adjoint operators, Friedrichs and Krein extensions
of densely
defined symmetric operators bounded from below with
finite deficiency
indices, and a class of Hamiltonian systems on a half-line.
Throughout this section we closely follow the setup in
Section~\ref{s4}.
In particular, we omit proofs whenever they parallel
the corresponding
scalar situation and focus on those arguments which
require new
elements when compared to Section~\ref{s4}.

Before we enter a discussion of these three cases, we
briefly digress
into the definition of $L^2$-spaces with underlying
matrix-valued
measures (see, e.g., \cite{DS88}, Sect.~XIII.5,
\cite{Na68}, Ch.~VI).
Suppose $\Omega=(\Omega_{j,k})_{1\leq j,k \leq n}$
generates a
matrix-valued measure on an interval
$\Lambda\subseteq\bbR$ as in
\eqref{5.4}--\eqref{5.7} with $\omega^{tr}
=\sum_{j=1}^n \Omega_{j,j}$, the corresponding scalar
trace measure. Let
$\ul{\hat f} =(\hat f_1,\dots ,\hat f_n)^t\in C_0(\Lambda)^n$,
where ``t''
abbreviates transpose and $C_0(\Lambda)$ denotes the set of
complex-valued continuous functions of compact support
contained in
$\Lambda$. On $C_0(\Lambda)^n$ we define the inner product,
\begin{equation}
(\ul{\hat f}, \ul{\hat g})_0=\sum_{j,k=1}^n
\int_\Lambda d\Omega_{j,k}
(\lambda)
\hat f_j(\lambda)\hat g_k(\lambda), \quad \ul{\hat f},
\ul{\hat g} \in C_0(\Lambda)^n.
\lb{7.1}
\end{equation}
The Hilbert space $L^2(\Lambda;d\Omega)$ is then defined
as the
completion of $C_0(\Lambda)^n$ with respect to the norm
$\|\cdot\|_0$
induced by \eqref{7.1}. A perhaps more useful, though
equivalent,
characterization of $L^2(\Lambda;d\Omega)$ can be obtained
as follows.
Introduce the density matrix $\rho$,
\begin{equation}
\rho(\lambda)=\big(\rho_{j,k}(\lambda)\big)_{1\leq j,k
\leq n}, \quad
\rho_{j,k}(\lambda) =\frac{d\Omega_{j,k}}
{d\omega^{tr}}(\lambda), \quad
j,k=1,\dots,n. \lb{7.2}
\end{equation}
Consider all complex-valued $\hat f_j:\Lambda\to\bbC$
such that
$\sum_{j,k=1}^n \overline{\hat f_j(\lambda)}\rho_{j,k}
(\lambda) \hat
f_k(\lambda)
\geq0$ is $\omega^{tr}$ integrable over $\Lambda$ and define
$\widehat \calH (\Lambda)$ as the set of equivalence classes
$\ul{\hat f}=
(\hat f_1,\dots ,\hat f_n)^t$ modulo $\Omega$-null
functions. (Here
$\ul{\hat g} =
(\hat g_1,\dots ,\hat g_n)^t$ is defined to be an
$\Omega$-null function
if $\int_\Lambda d\omega^{tr}(\lambda)
(\sum_{j,k=1}^n\overline{g_j(\lambda)}
\rho_{j,k}(\lambda)g_k(\lambda))=0$.) This space is
complete with respect
to the norm induced by the scalar product
\begin{equation}
(\ul{\hat f},\ul{\hat g})_{\widehat \calH (\Lambda)}=
\int_\Lambda d\omega^{tr}(\lambda)
\bigg(\sum_{j,k=1}^n \overline{f_j(\lambda)}
\rho_{j,k}(\lambda)
g_k(\lambda)\bigg),
\quad \ul{\hat f}, \ul{\hat g} \in \widehat \calH (\Lambda)
\lb{7.3}
\end{equation}
and coincides with $L^2(\Lambda;d\Omega)$,
\begin{equation}
\widehat \calH (\Lambda)=L^2(\Lambda;d\Omega). \lb{7.4}
\end{equation}

Now we turn to self-adjoint finite-rank perturbations of
self-adjoint
operators. Let $\calH$ be a separable complex Hilbert
space with
scalar product $(\cdot,\cdot)_{\calH}$, $H_0$ a self-adjoint
operator in
$\calH$ (which may or may not be bounded), and
$\{f_1,\dots,f_n\}\subset
\calH$ an orthogonal generating basis for $H_0$ (i.e.,
$(f_j,f_k)_\calH =\delta_{j,k}$, $j,k=1,\dots,n$ and
$\calH=\overline{\text{linspan}
\{(H_0-z)^{-1}f_j\in\calH\,|\,j=1,\dots,n, \, z\in
\bbC\backslash\bbR\}}$ =
$\overline{\text{linspan}\{E_0(\lambda)f_j\in\calH\,|
\, j=1,}$
\linebreak $\overline{\dots,n, \,
\lambda\in\bbR\}}$, $E_0(\cdot)$ the family of orthogonal
spectral projections of $H_0$).
Introducing the self-adjoint diagonal matrix
\begin{equation}
\alpha=\big(\alpha_j\delta_{j,k}\big)_{1\leq j,k\leq n},
\quad \alpha_j\in\bbR, \,\, j=1,\dots,n, \lb{7.5}
\end{equation}
we use the notation
\begin{equation}
\ul c=(c_1,\dots,c_n)^t\in\bbC^n, \quad \alpha\ul
c=(\alpha_1c_1, \dots,
\alpha_nc_n)^t, \text{ etc.} \lb{7.6}
\end{equation}
Moreover, we define
\begin{align}
&K:\bbC^n\to\calH, \quad \ul c\to \sum_{j=1}^nc_jf_j,
\lb{7.7}\\
&K^*:\calH\to\bbC^n, \quad f\to ((f_1,f)_\calH,
\dots,(f_n,f)_\calH )^t
\lb{7.8}
\end{align}
and note that
\begin{equation}
K\alpha K^*=\sum_{j=1}^n\alpha_j(f_j,\cdot)_\calH f_j.
\lb{7.9}
\end{equation}
After these preliminaries we can define the self-adjoint
finite-rank
perturbation of $H_0$ by
\begin{equation}
H_\alpha=H_0+K\alpha K^*=H_0+\sum_{j=1}^n\alpha_j(f_j,
\cdot)_\calH f_j,
\lb{7.10}
\end{equation}
with $\calD (H_\alpha)=\calD (H_0)$, $\alpha_j\in\bbR$,
$j=1,\dots,n$.
Denoting by $E_\alpha(\lambda)$, $\lambda\in\bbR$ the family of
orthogonal spectral projections of $H_\alpha$ one introduces
\begin{align}
&\Omega_\alpha(\lambda)=\big(\Omega_{\alpha,j,k}(\lambda)
\big)_{1\leq j,k\leq n}, \quad d\Omega_{\alpha,j,k}(\lambda)=
(f_j,dE_\alpha(\lambda)f_k)_\calH, \no\\
& \int_\bbR
d\Omega_{\alpha,j,k}(\lambda) =(f_j,f_k)_\calH =\delta_{j,k},
\quad
j,k=1,\dots,n. \lb{7.11}
\end{align}
By the canonical representation of self-adjoint operators
with finite
spectral multiplicity (cf., e.g., \cite{Na68}, Sect.~20),
$H_\alpha$ in
$\calH$ is unitarily equivalent to
$\widehat H_\alpha$ in $\hatt \calH_\alpha=L^2(\bbR;
d\Omega_\alpha)$, where
\begin{align}
& (\widehat H_\alpha \ul{\hat g})(\lambda)=
\lambda \ul{\hat g}(\lambda), \quad
\ul{\hat g}\in\calD (\widehat H_\alpha)=
L^2(\bbR;(1+\lambda^2)d\Omega_\alpha),
\lb{7.12} \\
& H_\alpha=U_\alpha \widehat H_\alpha U_\alpha^{-1}, \quad
\calH=U_\alpha L^2(\bbR;d\Omega_\alpha), \lb{7.13}
\end{align}
with $U_\alpha$ unitary,
\begin{align}
&U_\alpha:\widehat \calH_\alpha=L^2(\bbR;
d\Omega_\alpha)\to\calH,
\lb{7.14} \\
&\ul{\hat g}\to (U_\alpha \ul{\hat g}) = \slim_{N\to\infty}
\sum_{j=1}^n \int_{-N}^N d(E_\alpha(\lambda)f_j)\hat
g_j(\lambda), \quad
\ul{\hat g}=(\hat g_1,\dots,\hat g_n)^t\in
L^2(\bbR;d\Omega). \no
\end{align}
Moreover,
\begin{equation}
f_j=U_\alpha \ul{\hat f}_j, \quad \ul {\hat f}_j(\lambda)=
(\delta_{j,1},\dots,\delta_{j,n})^t, \quad
\lambda\in\bbR. \lb{7.15}
\end{equation}
The family of spectral projections $\widehat
E_\alpha(\lambda)$,
$\lambda\in\bbR$ of $\widehat H_\alpha$ is then
given by
\begin{equation}
(\widehat E_\alpha (\lambda) \ul{\hat g})(\mu)=
\theta(\lambda - \mu)
\ul{\hat g} (\mu) \text{ for $\Omega_\alpha$--a.e. }
\mu\in\bbR,
\,\, \ul{\hat g}\in L^2(\bbR;d\Omega_\alpha). \lb{7.16}
\end{equation}
Introducing the matrix-valued Herglotz function
\begin{align}
M_{\alpha}(z)&=\big((\ul e_j,K^*(H_\alpha-z)^{-1}
K\ul e_k)_{\bbC^n}
\big)_{1\leq j,k\leq n}
=\big((f_j,(H_{\alpha}-z)^{-1}f_k)_{\calH}\big)_{1\leq
j,k\leq n} \no \\
&=\int_{\bbR}
\frac{d\Omega_{\alpha}}{\lambda-z}, \quad z\in\bbC_+,
\lb{7.17}
\end{align}
with
\begin{equation}
\ul e_j=(\delta_{j,1},\dots,\delta_{j,n})^t \in\bbC^n,
\quad j=1,\dots,n, \lb{7.17a}
\end{equation}
one verifies
\begin{align}
&M_{\beta}(z)=M_{\alpha}(z)(I_n+(\beta-\alpha)
M_{\alpha}(z))^{-1},
\lb{7.18} \\
& \alpha=\big(\alpha_j\delta_{j,k}\big)_{1\leq
j,k\leq n}, \,
\beta=\big(\beta_j\delta_{j,k}\big)_{1\leq j,k\leq n},
\,\, \alpha_j,
\beta_j \in\bbR, \,\, j=1,\dots,n. \no
\end{align}
A comparison of \eqref{7.18} and \eqref{6.36a} suggests
the introduction of
\begin{align}
&A(\alpha,\beta)=\left(\begin{array}{cc} I_n &
\beta-\alpha \\ 0 &
I_n \end{array}\right) \in\calA_{2n},  \lb{7.19} \\
& \alpha=\big(\alpha_j\delta_{j,k}\big)_{1\leq
j,k\leq n}, \,
\beta=\big(\beta_j\delta_{j,k}\big)_{1\leq j,k\leq n},
\,\, \alpha_j,
\beta_j \in\bbR, \,\, j=1,\dots,n. \no
\end{align}
In particular, Theorem~\ref{t6.6} applies (with
$A_{1,1}(\alpha,\beta)=
A_{2,2}(\alpha,\beta)=I_n$, $A_{1,2}(\alpha,\beta)
=\beta-\alpha$,
$A_{2,1}(\alpha,\beta)=0$).

If $\{f_1,\dots,f_n\}$ is not a generating basis for
$H_0$, then
$\calH$ decomposes into
$\calH=\calH^n\oplus\calH^{n,\bot}$, with
$\calH^n=\ol{\text{linspan}\{(H_0-z)^{-1}f_j\,|
\,z\in\bbC,\,j=1,\dots,n\}}$ separable and $\calH^n$,
$\calH^{n,\bot}$ reducing subspaces for all $H_\alpha$.
The part
$H_\alpha^n$ of $H_\alpha$ in $\calH^n$ then plays
the role
analogous to $H_\alpha^1$ in the context of
\eqref{4.8aa}--\eqref{4.8e}.

Introducing the following set of Herglotz matrices
\begin{align}
\calN_1^{n \times n}=\{&M:\bbC_+\to M_n(\bbC)
\text{ Herglotz}\,|\,
M(z)=\smallint_\bbR d\Omega(\lambda)(\lambda-z)^{-1}, \no \\
& \text{\,for all } \ul c\in\bbC^n, \,\,
\smallint_\bbR
(\ul c,d\Omega(\lambda)\ul c)_{\bbC^n}<\infty \}, \lb{7.20}
\end{align}
we now turn to a realization theorem for Herglotz matrices
of the type
\eqref{7.17} and state the analog of Theorem~\ref{t4.1a}.

\begin{theorem} \mbox{\rm } \lb{t7.1}

\noindent (i). Any $M\in\calN_1^{n \times n}$ with
associated
measure $\Omega$ can be
realized in the form
\begin{align}
M(z)&=\big((\ul e_j,K^*(H-z)^{-1}K\ul e_k)_{\bbC^n}
\big)_{1\leq j,k\leq n}
\no \\
&=\big((f_j,(H-z)^{-1}f_k)_{\calH}
\big)_{1\leq j,k\leq n}, \quad z\in\bbC_+, \lb{7.21}
\end{align}
where $H$ denotes a self-adjoint operator in some
separable complex
Hilbert space $\calH$, $\{f_1,\dots,f_n\}\subset\calH$,
$(f_j,f_k)_\calH
=\delta_{j,k}$, $j,k=1,\dots,n$ and
\begin{equation}
\int_{\bbR}d\Omega(\lambda)=\big(\|f_j\|^2_{\calH}
\delta_{j,k}\big)_{1\leq j,k\leq n}. \lb{7.22}
\end{equation}

\noindent (ii). Suppose $M_\ell\in\calN_1^{n \times n}$
with
corresponding measures
$\Omega_\ell$, $\ell=1,2$, and $M_1\neq M_2$. Then $M_1$
and $M_2$ can be
realized as
\begin{align}
&M_\ell(z)=\big((\ul e_j,K^*(H_\ell-z)^{-1}K
\ul e_k)_{\bbC^n}\big)_{1\leq j,k\leq n} \no \\
&=\big((f_j,(H_\ell
-z)^{-1}f_k)_{\calH}\big)_{1\leq j,k\leq n},
\quad \ell=1,2, \,\,
z\in\bbC_+, \lb{7.23}
\end{align}
where $H_\ell$, $\ell=1,2$ are self-adjoint
finite-rank perturbations
of one and
the same self-adjoint operator $H_0$ in some complex
Hilbert space
$\calH$ (which may be chosen separable) with
$\{f_1,\dots,f_n\}\subset\calH$, $(f_j,f_k)_\calH
=\delta_{j,k}$, $j,k=1,\dots,n$, that is,
\begin{equation}
H_\ell=H_0+K\alpha_\ell K^*=H_0+\sum_{j=1}^n \alpha_{\ell,j}
(f_j,\cdot)_\calH f_j \lb{7.24}
\end{equation}
for some $\alpha_\ell=\big(\alpha_{\ell,j}
\delta_{j,k}\big)_{1\leq j,k\leq n}$,
$\alpha_{\ell,j}\in\bbR$, $j=1,\dots,n$, $\ell=1,2$, if
and only if
the following conditions hold:
\begin{equation}
\int_{\bbR}d\Omega_1 (\lambda)=\int_{\bbR}d\Omega_2 (\lambda)
=\big(\|f_j\|^2_{\calH}\delta_{j,k}\big)_{1\leq j,k\leq n},
\lb{7.25}
\end{equation}
and for all $z\in\bbC_+$,
\begin{equation}
M_2(z)=M_1(z)\big(I_n+
\big(\|f_j\|^{-2}_{\calH}\delta_{j,k}\big)_{1\leq j,k\leq n}
(\alpha_2 -\alpha_1)M_1(z)\big)^{-1}. \lb{7.26}
\end{equation}
\end{theorem}

Since the proof parallels that of Theorem~\ref{t4.1a}
step by step
we omit further details.

It is possible to extend this formalism to more general
classes of
(possibly unbounded) symmetric finite-rank (form)
perturbations of
$|H_0|$, see, for instance, \cite{AK97} and the references
therein.\\

Next we turn to a characterization of Friedrichs and Krein
extensions
of densely defined operators bounded from below with
deficiency indices $(n,n)$ (see also \cite{AT82}, \cite{DM91},
\cite{DM95},
\cite{DMT88}, \cite{Do65}, \cite{Kr47a}, \cite{Sk79},
\cite{Ts80}--\cite{Ts92}).

We start by describing a canonical representation of
densely defined
closed symmetric operators with deficiency indices $(n,n)$
as discussed
in \cite{Na68} (in close analogy to the scalar case treated
in detail by
Donoghue \cite{Do65}). Let $\calH$ be a separable
complex Hilbert space, $H$ a closed densely defined symmetric
operator
with domain $\calD (H)$ and deficiency indices $(n,n)$. Let
\begin{equation}
U_\alpha:\ker(H^*-i)\to\ker(H^*+i) \lb{7.27}
\end{equation}
be a linear isometric isomorpism
and parametrize all self-adjoint extensions $H_\alpha$ of
$H$
according to von Neumann's formula by
\begin{align}
&H_{\alpha}(g+(1 +U_\alpha )u_+)=Hg+i(1 -U_\alpha) u_+,
\lb{7.28} \\
&\calD (H_{\alpha})=\{(g+(1 +U_\alpha) u_+)\in\calD(H^*)\,|\,
g\in\calD (H),\, u_+\in\ker(H^*-i) \}. \no
\end{align}
In order to resemble the notation employed in
Section~\ref{s4}, we
think of $2\alpha$ as a self-adjoint
matrix representing $U_\alpha=e^{2i\alpha}\in U(n)$ with
respect to
fixed bases in
$\ker(H^*\mp i)$. (Here $U(n)$ denotes the set of unitary
$n \times n$
matrices with entries in $\bbC$.) Next, we assume that
$\{u_{+,j} \}_{1\leq j \leq n}$ is a generating
basis for $H_\alpha$ for some (and hence for all)
$e^{2i\alpha}\in U(n)$.
Let
$E_{\alpha}(\cdot)$ be the corresponding family of
orthogonal spectral
projections of $H_{\alpha}$ and define
\begin{align}
&d\Upsilon_\alpha (\lambda)=\big(d\Upsilon_{\alpha,j,k}
(\lambda)\big)_{1\leq j,k\leq n}, \quad
d\Upsilon_{\alpha,j,k}(\lambda)=
(u_{+,j},dE_{\alpha}(\lambda) u_{+,k})_{\calH}, \no \\
&\int_{\bbR}d\Upsilon_{\alpha,j,k}(\lambda)=(u_{+,j},
u_{+,k})_{\calH}
=\delta_{j,k},, \quad j,k=1,\dots,n, \,\,e^{2i\alpha}
\in U(n). \lb{7.31}
\end{align}
Then $H_{\alpha}$ is unitarily equivalent to
multiplication by $\lambda$
in $L^2(\bbR;d\Upsilon_{\alpha})$ and $u_{+,j}$ can be
mapped into the
vector $(\delta_{j,1},\dots,\delta_{j,n})^t$. However,
it is more
convenient to define
\begin{equation}
d\Omega_{\alpha}(\lambda)=(1+\lambda^2)
d\Upsilon_{\alpha}(\lambda),
\lb{7.32}
\end{equation}
such that
\begin{equation}
\int_{\bbR}\frac{d\Omega_{\alpha}(\lambda)}
{1+\lambda^2}=I_n, \quad
\int_{\bbR} (\ul c,d\Omega_{\alpha}(\lambda)
\ul c)_{\bbC^n}=\infty
\text{ for all } \ul c\in\bbC^n\backslash\{0\}, \,\,
e^{2i\alpha}\in U(n) \lb{7.33}
\end{equation}
(by \eqref{7.31} and the fact that $u_{+,j}
\notin\calD (H_{\alpha})$).
Thus,
$H_{\alpha}$ is unitarily equivalent to
${\hatt H}_{\alpha}$ in
$\hatt \calH_\alpha=L^2(\bbR;d\Omega_{\alpha})$, where
\begin{align}
 ({\hatt H}_{\alpha}\ul{\hat g})(\lambda)&=\lambda
\ul{\hat g}(\lambda),
\quad  \ul{\hat g}\in\calD ({\hatt H}_{\alpha})=
L^2(\bbR;(1+\lambda^2)d\Omega_{\alpha}), \lb{7.34} \\
& H_{\alpha}=U_{\alpha}{\hatt H}_{\alpha}
U_{\alpha}^{-1}, \quad
\calH =U_{\alpha}L^2(\bbR;d\Omega_{\alpha}), \lb{7.35}
\end{align}
with $U_{\alpha}$ unitary,
\begin{align}
& U_{\alpha}: \hatt \calH_{\alpha}=L^2(\bbR;
d\Omega_{\alpha}) \to
\calH, \no \\
& \ul{\hat g} \to U_{\alpha}\ul{\hat g} =\slim_{N\to\infty}
 \sum_{j=1}^n \int_{-N}^N
 d(E_{\alpha}(\lambda)u_{+,j})(\lambda-i) \hat g_j(\lambda),
\lb{7.36} \\
& \hspace*{3.5cm} \ul{\hat g}=(\hat g_1,\dots,\hat g_n)^t
\in L^2(\bbR;d\Omega_\alpha). \no
\end{align}
Moreover,
\begin{equation}
u_{+,j}=U_{\alpha}\ul{\hat u}_{+,j}, \quad
\ul{\hat u}_{+,j}(\lambda)=
(\lambda - i)^{-1} \ul e_j, \,\, j=1,\dots,n \lb{7.37}
\end{equation}
and
\begin{align}
&({\hatt H}(\alpha)\ul{\hat g})(\lambda)=
\lambda\ul{\hat g}(\lambda),
\lb{7.38} \\
&\ul{\hat g}\in\calD ({\hatt H}(\alpha))=
\{\ul{\hat h}\in\calD ({\hatt H}_{\alpha})\,|\,
\smallint_{\bbR}(\ul{e}_j,d\Omega_{\alpha}
(\lambda)\ul{\hat h}(\lambda))_{\bbC^n}=0,\,
j=1,\dots,n\}, \no
\end{align}
where $\ul e_j$ has been defined in \eqref{7.17a} and
\begin{equation}
H=U_{\alpha}{\hatt H}(\alpha)U_{\alpha}^{-1}. \lb{7.39}
\end{equation}
Thus $\hatt H(\alpha)$ in $L^2(\bbR;d\omega_{\alpha})$
is a canonical
representation for a densely defined closed symmetric
operator $H$ with deficiency indices $(n,n)$ in a
separable Hilbert space
$\calH$ and a generating basis $\{u_{+,j}
\in\ker(H^*-i)\}_{1\leq j\leq n}$.
We shall prove in Theorem~\ref{t7.2} below that
$\hatt H(\alpha)$
in $L^2(\bbR;d\Omega_\alpha)$ is actually a model
for all such operators.
Moreover, since
\begin{align}
& ((H-\overline z)g,U_{\alpha}(\cdot -z)^{-1}
\ul e_j)_{\calH}=
\int_{\bbR} (\lambda -z)((\ul{U_\alpha^{-1}g})(\lambda),
d\Omega_{\alpha}(\lambda)\ul e_j)_{\bbC^n}(\lambda -z)^{-1}
=0,  \no \\
& \hspace*{7.9cm} g\in\calD(H), \,\, z\in\bbC
\backslash\bbR, \lb{7.40}
\end{align}
by \eqref{7.38}, one infers that
$U_{\alpha}(\cdot -z)^{-1}\ul e_j\in\calD(H^{*})$.
Since $\calD(H)$
is dense in $\calH$, one concludes
\begin{equation}
\ker(\hatt H(\alpha)^{*} -z)=\{c_j(\cdot -z)^{-1}
\ul e_j\,|\,
c_j\in\bbC, \, j=1,\dots,n \}, \quad z\in
\bbC\backslash\bbR, \lb{7.41}
\end{equation}
where
\begin{equation}
H^{*}=U_{\alpha}\hatt H (\alpha)^{*}U_{\alpha}^{-1}. \lb{7.42}
\end{equation}

If $\{u_{+,j}\in\ker(H^*-i)\}_{1\leq j\leq n}$ is not a
generating basis
for $H_{\alpha}$ then, in close analogy to Section~\ref{s4},
$\calH$ (not necessarily assumed to be separable at
this point)
decomposes into two orthogonal subspaces $\calH^0$ and
$\calH^{0,\bot}$,
\begin{equation}
\calH=\calH^0\oplus\calH^{0,\bot}, \lb{7.43}
\end{equation}
with $\calH^0$ separable, each of which is
a reducing subspace for all $H_{\alpha}$, $e^{2i\alpha}
\in U(n)$
and
\begin{align}
& \calH^0=\overline{\text{linspan}\{(H_\alpha-z)^{-1}u_{+,j}
\in\calH
\,|\,z\in\bbC\backslash\bbR,\,j=1,\dots,n \}} \no \\
& \hspace*{.81cm} \text{ is independent of }
\alpha\in U(n), \lb{7.44} \\
& (H_{\alpha} -z)^{-1}=(H_{\beta} -z)^{-1} \text{ on }
\calH^{0,\bot}
\text{ for all } e^{2i\alpha},e^{2i\beta}\in U(n),
\,\, z\in\bbC_+. \lb{7.45}
\end{align}
In particular, the part $H^{0,\bot}$ of $H$ in
$\calH^{0,\bot}$ is
then self-adjoint,
\begin{align}
& H=H^0\oplus H^{0,\bot}, \quad H_{\alpha}=
H_{\alpha}^0\oplus H^{0,\bot},
\quad e^{2i\alpha}\in U(n), \lb{7.46} \\
& \ran(H^{0,\bot} -z) = \calH^{0,\bot}, \quad z
\in\bbC\backslash\bbR,
 \\
& u_+ = u_+^0 \oplus 0, \lb{7.47}
\end{align}
with $H^0$ a densely defined closed symmetric operator
in $\calH^0$
and deficiency indices $(n,n)$. One then computes
\begin{align}
& z(u_{+,j},u_{+,k})_{\calH} +
(1+z^2)(u_{+,j},(H_{\alpha} -z)^{-1}u_{+,k})_{\calH}
\no  \\
& =z(u_{+,j}^{0},u_{+,k}^{0})_{\calH^0} +
(1+z^2)(u_{+,j}^0,(H_{\alpha}^0 -z)^{-1}u_{+,k}^0)_{\calH^0},
\lb{7.48}\\
& \hspace*{4.3cm} j,k=1,\dots,n, \,\, e^{2i\alpha}\in
U(n)  \no
\end{align}
and hence $\alpha$-dependent spectral properties of
$H_{\alpha}$
in $\calH$ are again effectively
reduced to those of $H_{\alpha}^0$ in $\calH^0$, where
$H_{\alpha}^0$
are self-adjoint operators with a generating basis
$\{u_{+,j}^0\in\ker((H^0)^* -i)\}_{1\leq j\leq n}$. We
shall call the
densely defined closed symmetric operator $H$ with
deficiency indices
$(n,n)$ {\it prime} if $\calH^{0,\bot}=\{0\}$ in \eqref{7.43}.

Next we show the model character of
$(\hatt \calH_\alpha, \hatt H(\alpha), \hatt H_\alpha)$
following the
approach outlined in Theorem~\ref{t4.1b}.

\begin{theorem} \lb{t7.2}
Let $H$ be a densely defined closed prime symmetric operator
with deficiency
indices $(n,n)$ and normalized deficiency vectors
$u_{\pm,j}\in\ker(H^* \mp i)$, $\|u_{\pm,j}\|_\calH$ $=1$,
$j=1,\dots,n$
in some separable
complex Hilbert space $\calH$. Let
$H_\alpha$ be a self-adjoint extension of $H$ with
generating orthonormal
basis $\{u_{+,j}\in\ker(H^*-i)\,|\, j=1,\dots,n \}$. Then
the pair
$(H,H_\alpha)$ in $\calH$ is unitarily equivalent to the pair
$(\hatt H(\alpha), \hatt H_\alpha)$ in $\hatt \calH$
defined in
\eqref{7.38} and \eqref{7.34} with unitary operator
$U_\alpha$
defined in
\eqref{7.36} (cf. \eqref{7.39} and \eqref{7.35}).
Conversely, given a
matrix-valued measure $d\widetilde \Omega$ satisfying
\begin{equation}
\int_{\bbR}\frac{d\widetilde\Omega(\lambda)}
{1+\lambda^2}=I_n, \quad
\int_{\bbR} (\ul c,d\widetilde
\Omega(\lambda)\ul c)_{\bbC^n}=\infty
\text{ for all } \ul c\in\bbC^n\backslash\{0\},\lb{7.49}
\end{equation}
define the
self-adjoint operator $\widetilde H$ of multiplication
by $\lambda$
in $\widetilde \calH=L^2(\bbR;d\widetilde \Omega)$,
\begin{equation}
({\widetilde H} \ul {\widetilde g})(\lambda)=\lambda
\ul {\widetilde g}(\lambda), \quad
\ul {\widetilde g}\in\calD(\widetilde H)
= L^2(\bbR;(1+\lambda^2)d\widetilde \Omega) \lb{7.50}
\end{equation}
and the linear operator $H$ in $\widetilde \calH$,
\begin{equation}
\calD(H)=\{\ul {\widetilde g}\in\calD(\widetilde H)\,|\,
\smallint_{\bbR}
(\ul e_j,d\widetilde \Omega(\lambda)
\ul {\widetilde g}(\lambda))_{\bbC^n}=0, \,
j=1,\dots,n \}, \quad H=
\widetilde H \big|_{\calD(H)}. \lb{7.51}
\end{equation}
Then $H$ is a densely defined closed symmetric operator in
$\widetilde \calH$ with
deficiency indices $(n,n)$ and deficiency subspaces
\begin{equation}
\ker(H^{*} \mp i)=\{c_j(\lambda \mp i)^{-1}\ul e_j
\,|\, c_j\in\bbC, \,
j=1,\dots,n \}. \lb{7.52}
\end{equation}
\end{theorem}
\begin{proof}
Except for a few modifications one can follow the
corresponding
proof
of Theorem~\ref{t4.1b} step by step. In particular, the
first part
goes through with the obvious changes indicated in
\eqref{7.31}--\eqref{7.48}. Hence we briefly turn to the
proof of the
second part of Theorem~\ref{t7.2}. Given
\eqref{7.2}--\eqref{7.4},
the scale of Hilbert spaces is still defined by
 $\widetilde \calH_{2r}=L^2(\bbR;(1+\lambda^2)^rd\widetilde
\Omega)$,
$r\in\bbR$, $\widetilde \calH_0=\widetilde \calH$ and
one considers
again the unitary operator $R$,
\begin{equation}
 R:\widetilde \calH_2\to\widetilde \calH_{-2}, \quad
{\ul {\widetilde f}}\to (1+\lambda^2){\ul {\widetilde f}}.
\lb{7.53}
\end{equation}
We note that $\bbC^n\subset\widetilde \calH_{-2}$. Again
$\calD(H)$ is
well-defined, and as a restriction of the self-adjoint
operator
$\widetilde H$, $H$ is clearly
symmetric. By \eqref{7.53} and \eqref{7.51} one infers
\begin{equation}
\calD(\widetilde H)=\widetilde \calH_2=
\calD(H)\oplus_{\tilde \calH_2}R^{-1}\bbC^n, \lb{7.54}
\end{equation}
which allows one to prove that $\calD(H)$ is dense in
$\widetilde \calH$ as in the proof of Theorem~\ref{t4.1b}.
That
$H$ is a closed operator is also proved as in Section~\ref{s4}.
Since $\widetilde H$ is self-adjoint, $\ran(\widetilde H -z)=
\widetilde \calH$
for all $z\in\bbC\backslash\bbR$, and
$(\widetilde H \pm i):\widetilde \calH_2\to\widetilde
\calH$ is unitary.
Together with \eqref{7.54} this yields
\begin{align}
\widetilde \calH&=(\widetilde H\pm i)\widetilde \calH_2 =
(\widetilde H\pm i)(\calD(H)\oplus_{\tilde \calH_2}R^{-1}\bbC)
=(H\pm i)\calD(H)\oplus_{\tilde \calH}
\big \{\tfrac{\lambda\pm i}{1+\lambda^2}c
\,\big |\, c\in\bbC \big \} \no \\
& =\ran(H\pm i)\oplus_{\tilde \calH}
\{c(\lambda \mp i)^{-1} \,|\,
c\in\bbC \} \lb{7.56}
\end{align}
and hence \eqref{7.52}.
\end{proof}

Introducing the Herglotz function
\begin{align}
M_{\alpha}(z)&=\int_{\bbR}d\Omega_{\alpha}(\lambda)
((\lambda-z)^{-1}-
\lambda(1+\lambda^2)^{-1})  \lb{7.57} \\
&=zI_n+(1+z^2)\big((u_{+,j},(H_{\alpha}-z)^{-1}
u_{+,k})_{\calH}
\big)_{1\leq j,k\leq n} \lb{7.58}
\end{align}
(the last equality being a simple consequence of
$\int_{\bbR}d\Omega_{\alpha}(\lambda)
(1+\lambda^2)^{-1}=I_n$) one verifies,
\begin{align}
M_{\beta}(z)&=(-e^{-i\beta}(\sin(\beta)\cos(\alpha)-
\cos(\beta)\sin(\alpha))e^{i\alpha}  \no \\
&\quad \,\, +e^{-i\beta}(\cos(\beta)\cos(\alpha)
+\sin(\beta)\sin(\alpha))
e^{i\alpha} M_{\alpha}(z)) \no \\
&\quad \times (e^{-i\beta}(\cos(\beta)\cos(\alpha)
+\sin(\beta)\sin(\alpha))e^{i\alpha}  \no \\
&\quad \quad \, +e^{-i\beta} (\sin(\beta)\cos(\alpha)-
\cos(\beta)\sin(\alpha))e^{i\alpha}M_{\alpha}(z))^{-1},
\lb{7.59} \\
& \hspace*{4.5cm} \exp(2i\alpha),\exp(2i\beta)\in U(n).
\no
\end{align}
Since \eqref{7.59} does not seem to be a well-known
result, we will
provide its derivation, following \cite{GMT97}, in
Appendix~\ref{B}.
A comparison of \eqref{7.59} and \eqref{6.36a} suggests
invoking
\begin{align}
&A(\alpha,\beta)= \big(A(\alpha,\beta)_{j,k}
\big)_{1\leq j,k\leq n}
\in\calA_{2n}, \lb{7.60} \\
&A(\alpha,\beta)_{1,1}=e^{-i\beta}(\cos(\beta)
\cos(\alpha)
+\sin(\beta)\sin(\alpha))e^{i\alpha}, \no \\
& A(\alpha,\beta)_{1,2}=e^{-i\beta} (\sin(\beta)
\cos(\alpha)-
\cos(\beta)\sin(\alpha))e^{i\alpha}, \no \\
&A(\alpha,\beta)_{2,1}=e^{-i\beta}(\cos(\beta)
\sin(\alpha)-
\sin(\beta)\cos(\alpha))e^{i\alpha}, \no \\
&A(\alpha,\beta)_{2,2}=e^{-i\beta}(\cos(\beta)
\cos(\alpha)
+\sin(\beta)\sin(\alpha))e^{i\alpha},  \no \\
& \hspace*{3.8cm} \exp(2i\alpha),\exp(2i\beta)\in U(n).
\no
\end{align}
Moreover, Theorem~\ref{t6.6} applies (with
$A_{1,1}(\alpha,\beta)=
A_{2,2}(\alpha,\beta)=e^{-i\beta}(\cos(\beta)
\cos(\alpha)$
$+\sin(\beta)\sin(\alpha))e^{i\alpha}$,
$A_{1,2}(\alpha,\beta)=-A_{2,1}(\alpha,\beta)
=e^{-i\beta} (\sin(\beta)\cos(\alpha)-
\cos(\beta)\sin(\alpha))$ $e^{i\alpha}$). Since by
definition,
$M_\alpha(i)=iI_n$, Lemma~\ref{l5.2a} yields
$\Im(M_\alpha(z))>0$ for
all $z\in\bbC_+$. In fact, Lemma~\ref{l3a} yields an
explicit lower
bound for $\Im(z)\Im(M_\alpha(z))$. \\

Next, assuming that $H$ is nonnegative, $H\geq 0$, we
again intend to
characterize
the Friedrichs and Krein extensions, $H_F$ and $H_K$, of
$H$. In order
to apply Krein's results \cite{Kr47a} (see also
\cite{AT82}, \cite{Ts80},
\cite{Ts81}) in a slightly different
form (see, e.g., \cite{Sk79}, Sect.~4 for an efficient
summary of
Krein's results most relevant in our context) we start
with the analog of
Theorem~\ref{t4.1c}.

\begin{theorem} \mbox{\rm } \lb{t7.3}

\noindent (i). $H_{\alpha}=H_F$ for some
$e^{2i\alpha}\in U(n)$ if and
only if for all $R>0$,
$\int_R^{\infty}d\|E_{\alpha}(\lambda)u_+\|_{\calH}^2
\lambda$ $=\infty$
for all $0\neq u_+\in\ker(H^*-i)$, or
equivalently, if and only if  $\int_R^{\infty}
(\ul c,d\Omega_{\alpha}(\lambda)\ul c)_{\bbC^n}
\lambda^{-1}$ $=\infty$
for all $\ul c\in\bbC^n\backslash\{0\}$.

\noindent (ii). $H_{\beta}=H_K$ for some
$e^{2i\beta}\in U(n)$ if and
only if for all $R>0$,
$\int^R_{0} d\|E_{\beta}(\lambda)u_+\|_{\calH}^2$
$\lambda^{-1}=\infty$ for all
$0\neq u_+\in\ker(H^*-i)$, or
equivalently, if and only if $\int^R_{0}
(\ul c,d\Omega_{\beta}(\lambda)\ul c)_{\bbC^n}$
$\lambda^{-1}=\infty$
for all $\ul c\in\bbC^n\backslash\{0\}$.

\noindent (iii). $H_{\gamma}=H_F=H_K$ for some
$e^{2i\gamma}\in U(n)$ if
and only if
$\int_R^{\infty} d\|E_{\gamma}(\lambda)u_+\|_{\calH}^2
\lambda
=\infty=$ $\int^R_{0}
d\|E_{\gamma}(\lambda)u_+\|_{\calH}^2 \lambda^{-1}$
for all $R>0$ and all $0\neq u_+\in\ker(H^*-i)$, or
equivalently, if and only if  $\int_R^{\infty}
(\ul c,d\Omega_{\gamma}(\lambda)\ul c)_{\bbC^n}
\lambda^{-1} =
\int^R_{0} (\ul c,d\Omega_{\gamma}(\lambda)\ul c)_{\bbC^n}
\lambda^{-1}
=\infty$ for all $\ul c\in\bbC^n\backslash\{0\}$.
\end{theorem}
\begin{proof}
As in Section~\ref{s4}, in order to reduce the above
statements (i)--(iii)
to those in Krein
\cite{Kr47a} (as summarized in Skau \cite{Sk79}), it
suffices
to notice  that $(\mu +1)^{-1}-(\mu -i)^{-1}
=\Oh(\mu^{-2})$ as
$\mu\uparrow\infty$ and $(\mu +1)^{-1}-(\mu -i)^{-1}
=\Oh(1)$ as
$\mu\downarrow 0$.
\end{proof}

Of course \eqref{4.21}--\eqref{4.22a} remain valid in
the present case
of deficiency indices $(n,n)$.

Applying Theorem~\ref{t7.3} to $H_{\alpha}$ then yields
the analog of
Theorem~\ref{t4.1}\,(i)--(iii).

\begin{theorem} \mbox{\rm (\cite{DM91}, \cite{DM95},
\cite{DMT88},
\cite{KO78}). } \lb{t7.4}

\noindent (i). $H_{\alpha}=H_F$ if and only
if $\lim_{\lambda\downarrow
-\infty} (\ul c,M_{\alpha}(\lambda)\ul c)_{\bbC^n}=-\infty$
for all $\ul c\in\bbC^n\backslash\{0\}$.

\noindent (ii). $H_{\beta}=H_K$ if and only if
$\lim_{\lambda\uparrow
0}(\ul c,M_{\beta}(\lambda)\ul c)_{\bbC^n}=\infty$
for all
$\ul c\in\bbC^n\backslash\{0\}$.

\noindent (iii). $H_{\gamma}=H_F=H_K$ if and only if for all
$\ul c\in\bbC^n\backslash\{0\}$,
$\lim_{\lambda\downarrow
-\infty}(\ul c,M_{\gamma}(\lambda)\ul c)_{\bbC^n}=-\infty$
and
$\lim_{\lambda\uparrow
0}(\ul c,M_{\gamma}(\lambda)\ul c)_{\bbC^n}=\infty$.

\end{theorem}

Since the proof parallels the corresponding one in
Section~\ref{s4}
step by step we omit further details.

If $H_\alpha$ and $H_\beta$ are two distinct self-adjoint
extensions
of the symmetric operator $H$ with deficiency indices
$(n,n)$, $n\geq 2$
considered in Theorem~\ref{t7.2}, then $\calD (H_\alpha)$
and
$\calD (H_\beta)$ may have a nontrivial intersection,
that is,
\begin{equation}
\calD (H_\alpha) \cap \calD (H_\beta) \underset{\not=}
{\supset} \calD (H),
\quad e^{2i\alpha},e^{2i\beta} \in U(n), \,\,
U_\alpha\neq U_\beta.
\lb{7.56a}
\end{equation}
Next, we characterize the case where the domain of a
nonnegative
self-adjoint extension
$\widetilde H$ of $H$ has only trivial intersection with
that of $H_F$ or
$H_K$. These results go beyond those in \cite{Kr47a} and
appear to be new.

\begin{theorem} \mbox{\rm } \lb{t7.5} Suppose
$\widetilde H\geq 0$ is a nonnegative self-adjoint extension
of a
densely defined nonnegative closed operator $H\geq 0$ with
deficiency
indices $(n,n)$. We denote by $\widetilde E(\lambda)$ the
family of
spectral projections of $\widetilde H$ and by
$d\widetilde \Omega(\lambda)$ the measure defined in
\eqref{7.31},
\eqref{7.32}. Then

\noindent (i). $\calD (\widetilde H)\cap\calD (H_F)
=\calD (H)$ if and
only if for all $R>0$,
$\int_R^{\infty} d\|\widetilde E(\lambda)u_+\|_{\calH}^2
\lambda
<\infty$ for all $u_+\in\ker(H^*-i)$, or
equivalently, if and only if  $\int_R^{\infty}
(\ul c,d\widetilde \Omega(\lambda)
\ul c)_{\bbC^n} \lambda^{-1} <\infty$
for all $\ul c\in\bbC^n$.

\noindent (ii). $\calD (\widetilde H)\cap\calD (H_K)
=\calD (H)$
if and only if for all $R>0$,
$\int^R_{0} d\|\widetilde E(\lambda)u_+\|_{\calH}^2
\lambda^{-1} <\infty$
for all $u_+\in\ker(H^*-i)$, or
equivalently, if and only if $\int^R_{0}
(\ul c,d\widetilde \Omega(\lambda)\ul c)_{\bbC^n}
\lambda^{-1} <\infty$
for all $\ul c\in\bbC^n$.

\noindent (iii). $\calD (\widetilde H)\cap\calD (H_F)
=\calD (H)=
\calD (\widetilde H)\cap\calD (H_K)$ if and only if for
all $R>0$,
$\int_R^{\infty} d\|\widetilde E(\lambda)u_+\|_{\calH}^2
\lambda $
$+\int^R_{0} d\|\widetilde E(\lambda)u_+\|_{\calH}^2
\lambda^{-1} <\infty$
for all $u_+\in\ker(H^*-i)$, or
equivalently, if and only if  $\int_0^{\infty}
(\ul c,d\widetilde \Omega(\lambda)\ul c)_{\bbC^n}
\lambda^{-1} <\infty$
for all $\ul c\in\bbC^n$.
\end{theorem}
\begin{proof}
It is sufficient to prove item~(i) since the remaining proofs
offer no new
details. We use the terminology introduced in
Appendix~\ref{B} and
identify $A_1=\wti H$, $A_2=H_F$, $P_{1,2}=\wti P_F$,
${\calU}_1=\wti {\calU}$, ${\calU}_2={\calU}_F$, etc.
First we
suppose that $\calD(\wti H)\cap\calD(H_F)=\calD(H)$. Using
von Neumann's
parametrization of $\wti H$ and $H_F$ in terms of the linear
isometric
isomorphisms $\wti {\calU}$ and ${\calU}_F$ from
$\ker(H^*-i)$ onto
$\ker(H^*+i)$, this assumption is equivalent to $-1$ not
being an eigenvalue
of $U_F$ (the matrix representation of ${\calU}_F$ in the
orthogonal bases of
$\ker(H^*\mp i)$ as discussed in Appendix~\ref{B}). By
\eqref{e19}, this is
equivalent to the existence of the inverse of $\wti P_F (i)$.
In order to prove
that $\int_R^{\infty} (\ul c,d\widetilde \Omega(\lambda)
\ul c)_{\bbC^n} \lambda^{-1} <\infty$ for all
$\ul c\in\bbC^n$, it suffices
to prove that the Herglotz matrix $\widetilde M(z)$
associated
with the
measure $d\widetilde \Omega(\lambda)$ corresponding to
$\wti H$ has a limit
as $z\to\ -\infty$. Using \eqref{e22}, one computes
\begin{equation}
\wti M(z)-\Re(\wti P_F (i)^{-1})=(2iI_n-\wti P_F(-i)^{-1})
(\wti P_F (-i)^{-1}
-iI_n+M_F(z))^{-1}\wti P_F (-i)^{-1}. \lb{7.60aa}
\end{equation}
Here $M_F(z)$ denotes the Herglotz matrix associated to
$H_F$ and we used
the fact
\begin{equation}
\Re(\wti P_F (i)^{-1})=iI_n+\wti P_F (i)^{-1}=-iI_n+\wti
P_F (-i)^{-1}.
\lb{7.60a}
\end{equation}
Next, recalling Theorem~\ref{t7.4}\,(i), we will invoke
that
\begin{equation}
\lim_{\lambda\downarrow -\infty}
(\ul c,M_F(\lambda)\ul c)_{\bbC^n} =
-\infty \text{ for all } \ul c\in\bbC^n\backslash\{0\}.
\lb{7.60b}
\end{equation}
Since $(\ul c,M_F(\lambda)\ul c)_{\bbC^n}^{-1}$ converges
monotonically to
zero pointwise
for any $\ul c\in\{\ul d\in\bbC^n\,|\,
\|\ul d\|_{\bbC^n}=1\}$,
the compact
unit
sphere in $\bbC^n$, Dini's theorem yields in fact uniform
convergence to zero.
Consequently,
\begin{equation}
\wti P_F (-i)^{-1}-iI_n+M_F(\lambda)\leq \gamma(\lambda)I_n,
\lb{7.60c}
\end{equation}
for $\lambda$ sufficiently negative and some
$\gamma(\lambda)$ with $\gamma
(\lambda)\downarrow -\infty$ as $\lambda\downarrow -\infty$.
In particular,
\begin{equation}
(\wti P_F (-i)^{-1}-iI_n+M_F(\lambda))^{-1}\to 0
\text{ as }
\lambda\downarrow -\infty. \lb{7.60d}
\end{equation}
\eqref{7.60aa} and \eqref{7.60d} then prove
\begin{equation}
\lim_{\lambda\downarrow -\infty} \wti M (\lambda)
=\Re(\wti P_F (i)^{-1}).
\lb{7.60e}
\end{equation}
Conversely, we suppose $\int_R^{\infty}
(\ul c,d\widetilde \Omega(\lambda)
\ul c)_{\bbC^n} \lambda^{-1} <\infty$ for all
$\ul c\in\bbC^n$ or
equivalently,
$\lim_{\lambda\downarrow -\infty} \wti M(\lambda)=
\wti M(-\infty)$ exists.
Similarly to \eqref{e28}, one derives
\begin{equation}
M_F(z)-\wti M(z)=(iI_n+\wti M(z))(I_n+i\wti P_F (i)-
\wti P_F (i)
\wti M(z))^{-1}\wti P_F (i)
(-iI_n+\wti M(z)) \lb{7.60f}
\end{equation}
and hence
\begin{align}
&((-iI_n+\wti M(\lambda))^{-1}\ul c,(M_F(\lambda)-
\wti M(\lambda))
(-iI_n+\wti M(\lambda))^{-1}\ul c)_{\bbC^n} \no \\
&=(\ul c,(I_n+i\wti P_F (i)-\wti P_F (i)
\wti M(\lambda))^{-1}
\wti P_F (i)\ul c)_{\bbC^n},
 \quad \lambda < 0, \,\, \ul c\in\bbC^n\backslash\{0\}.
\lb{7.60g}
\end{align}
By \eqref{7.60b} and the existence of $\wti M(-\infty)$,
the left-hand side
of \eqref{7.60g} tends to $-\infty$ as
$\lambda\downarrow -\infty$.
Consequently,
$\ker(\wti P_F (i))=\{0\}$, that is, $\wti P_F (i)$ is
invertible. By
\eqref{e19}, this is
equivalent to $-1$ not being an eigenvalue of $U_F$
implying
$\calD(\wti H)\cap\calD(H_F)=\calD(H)$.
\end{proof}

Theorem~\ref{t7.5} then yields the following result.

\begin{theorem} \mbox{\rm } \lb{t7.6} Suppose
$\widetilde H\geq 0$ is a nonnegative self-adjoint
extension of a
densely defined nonnegative closed operator $H\geq 0$
with deficiency
indices $(n,n)$. We denote by $\widetilde M(z)$ the
corresponding
Herglotz matrix associated with the measure
$d\widetilde \Omega(\lambda)$ defined in \eqref{7.31},
\eqref{7.32}.
Then

\noindent (i). $\calD (\widetilde H)\cap\calD (H_F)
=\calD (H)$ if and only
if  $\lim_{\lambda\downarrow
-\infty} |(\ul c,\widetilde M(\lambda)
\ul c)_{\bbC^n}|<\infty$ for all
$\ul c\in\bbC^n$.

\noindent (ii). $\calD (\widetilde H)\cap\calD (H_K)
=\calD (H)$ if and
only if  $\lim_{\lambda\uparrow
0} |(\ul c,\widetilde M(\lambda)\ul c)_{\bbC^n}|
<\infty$ for all
$\ul c\in\bbC^n$.

\noindent (iii). $\calD (\widetilde H)\cap\calD (H_F)
=\calD (H)=
\calD (\widetilde H)\cap\calD (H_K)$ if and only if
for all
$\ul c\in\bbC^n$,
$\lim_{\lambda\downarrow
-\infty} |(\ul c,\widetilde M(\lambda)\ul c)_{\bbC^n}|+
\lim_{\lambda\uparrow
0} |(\ul c,\widetilde M(\lambda)\ul c)_{\bbC^n}|<\infty$.
\end{theorem}

An analog of Theorem~\ref{t4.1}\,(iv) can now be obtained
as follows.

\begin{theorem} \mbox{\rm } \lb{t7.6a} Let
$\widetilde H\geq 0$ be a nonnegative self-adjoint
extension of a
densely defined nonnegative closed operator $H\geq 0$
with deficiency
indices $(n,n)$. We denote by $\widetilde M(z)$ the
corresponding
Herglotz matrix associated with the measure
$d\widetilde \Omega(\lambda)$ defined in \eqref{7.31},
\eqref{7.32} and
identify $A_1=\wti H$, $A_2=H_F \text{ or } H_K$,
$P_{1,2}=\wti P_F
\text{ or } \wti P_K$, ${\calU}_1=\wti \calU$,
${\calU}_2={\calU}_F
\text{ or } {\calU}_K$, etc., in Appendix~\ref{B}. Then

\noindent (i). If $\calD (\widetilde H)\cap\calD (H_F)
=\calD (H)$ then
\begin{equation}
 \lim_{\lambda\downarrow -\infty} \widetilde M(\lambda)=
\Re(\wti P_F (i)^{-1})=-\int_\bbR d\wti \Omega (\lambda)
\lambda
(1+\lambda^2)^{-1}. \lb{7.60h}
\end{equation}

\noindent (ii). If $\calD (\widetilde H)\cap\calD (H_K)
=\calD (H)$ then
\begin{equation}
 \lim_{\lambda\uparrow 0} \widetilde M(\lambda)=
\Re(\wti P_K (i)^{-1})=\int_\bbR d\wti \Omega (\lambda)
(\lambda^{-1}
-\lambda (1+\lambda^2)^{-1}). \lb{7.60i}
\end{equation}

\noindent (iii). If $\calD (\widetilde H)\cap\calD (H_F)
=\calD (H)=
\calD (\widetilde H)\cap\calD (H_K)$ then
\begin{equation}
\int_\bbR d\wti \Omega (\lambda) \lambda^{-1}=\Re(\wti
P_K (i)^{-1}
-\wti P_F (i)^{-1}). \lb{7.60j}
\end{equation}
\end{theorem}
\begin{proof}
Item (i) is clear from \eqref{7.57} and \eqref{7.60e}.
Similarly,
(ii) follows from \eqref{7.57} and
\begin{equation}
\lim_{\lambda\uparrow 0} \wti M(\lambda)=\Re(\wti
P_K (i)^{-1}. \lb{7.60k}
\end{equation}
Finally, (iii) is obvious by taking the difference of
\eqref{7.60i} and
\eqref{7.60h}.
\end{proof}

Next, we turn to a realization theorem for Herglotz
functions of the type
\eqref{7.58}. It will be convenient to introduce the
following sets of Herglotz matrices,
\begin{align}
&\calN_0^{n \times n} =\{M:\bbC_+\to M_n (\bbC)
\text{ Herglotz}\,|\,
M(z)=\smallint_{\bbR} d\Omega(\lambda)((\lambda -z)^{-1}
-\lambda(1+\lambda^2)^{-1}), \no \\
& \hspace*{1.15cm} \text{ for all }
\ul c\in\bbC^n\backslash\{0\}, \,\,\smallint_{\bbR}
(\ul c,d\Omega(\lambda)\ul c)_{\bbC^n}=\infty, \,\,
\smallint_{\bbR} (\ul c,d\Omega(\lambda)\ul c)_{\bbC^n}
(1+\lambda^2)^{-1} <\infty \}, \lb{7.61} \\
&\calN_{0,F}^{n \times n} =\{M\in\calN_0^{n \times n}\,|\,
\supp(\omega^{tr})\subseteq[0,\infty),
\text{ for all } \ul c\in\bbC^n\backslash\{0\}, \no \\
& \hspace*{1.6cm} \smallint_R^{\infty}
(\ul c,d\Omega(\lambda)
\ul c)_{\bbC^n}\lambda^{-1}
=\infty \text{ for some }R>0 \}, \lb{7.62} \\
&\calN_{0,K}^{n \times n} =\{M\in\calN_0^{n \times n}\,|\,
\supp(\omega^{tr})\subseteq[0,\infty), \text{ for all }
\ul c\in\bbC^n\backslash\{0\}, \no \\
& \hspace*{1.6cm}\smallint_0^{R} (\ul c,d\Omega(\lambda)
\ul c)_{\bbC^n}\lambda^{-1}
=\infty \text{ for some }R>0 \}, \lb{7.63} \\
&\calN_{0,F,K}^{n \times n} =\{M\in\calN_0^{n \times n}\,|\,
\supp(\omega^{tr})\subseteq[0,\infty), \text{ for all }
\ul c\in\bbC^n\backslash\{0\}, \no \\
& \hspace*{1.6cm} \smallint_R^{\infty}
(\ul c,d\Omega(\lambda)
\ul c)_{\bbC^n}\lambda^{-1}
=\smallint_0^{R} (\ul c,d\Omega(\lambda)
\ul c)_{\bbC^n}\lambda^{-1}
=\infty \text{ for some }R>0 \} \\
&\hspace*{1.1cm}=\calN_{0,F}^{n \times n}
\cap\calN_{0,K}^{n \times n}.
\lb{7.64} \\
&\calN_{0,F^{\bot}}^{n \times n}
=\{M\in\calN_0^{n \times n}\,|\,
\supp(\omega^{tr})\subseteq[0,\infty),
\text{ for all } \ul c\in\bbC^n, \no \\
& \hspace*{1.6cm} \smallint_R^{\infty}
(\ul c,d\Omega(\lambda)
\ul c)_{\bbC^n}\lambda^{-1}
<\infty \text{ for some }R>0 \}, \lb{7.65} \\
&\calN_{0,K^{\bot}}^{n \times n}
=\{M\in\calN_0^{n \times n}\,|\,
\supp(\omega^{tr})\subseteq[0,\infty), \text{ for all }
\ul c\in\bbC^n, \no \\
& \hspace*{1.6cm}\smallint_0^{R} (\ul c,d\Omega(\lambda)
\ul c)_{\bbC^n}\lambda^{-1}
<\infty \text{ for some }R>0 \}, \lb{7.66} \\
&\calN_{0,F^{\bot},K^{\bot}}^{n \times n}
=\{M\in\calN_0^{n \times n}\,|\,
\supp(\omega^{tr})\subseteq[0,\infty),
\text{ for all }
\ul c\in\bbC^n, \no \\
&\hspace*{2.1cm}\smallint_0^{\infty}
(\ul c,d\Omega(\lambda)
\ul c)_{\bbC^n}\lambda^{-1}
<\infty \}=
\calN_{0,F^{\bot}}^{n \times n}
\cap\calN_{0,K^{\bot}}^{n \times n}.
\lb{7.67}
\end{align}

The sets in \eqref{7.62}--\eqref{7.66} are of course
independent of $R>0$. The analog of Theorem~\ref{t4.3}
then reads as
follows.

\begin{theorem}  \mbox{\rm } \lb{t7.7}

\noindent (i). Any \,
$\widetilde M\in\calN_0^{n \times n}$ can be
realized in the form
\begin{align}
\widetilde M(z) =
\,\,&z\big(\|u_{+,j}\|^2_{\tilde\calH}\delta_{j,k}
\big)_{1\leq j,k\leq n}
\lb{7.68} \\
&+(1+z^2)\big((u_{+,j},
(\widetilde H -z)^{-1}u_{+,k})_{\tilde\calH}
\big)_{1\leq j,k\leq n}, \quad z\in\bbC_+, \no
\end{align}
where $\widetilde H$ denotes the self-adjoint extension of
some densely
defined closed symmetric operator $H$ with deficiency
indices $(n,n)$
and deficiency subspace $\{u_{+,j}
\in\ker(H^{*}-i)\}_{1\leq j,k\leq n}$ in
some separable complex Hilbert space $\widetilde \calH$.

\noindent (ii). Any \,
$\widetilde M_{F(\text{resp.~}K)}\in
\calN_{0,F(\text{resp.~}K)}^{n \times n}$ can be
realized in the form
\begin{align}
\widetilde M_{F(\text{resp.~}K)}(z) =
\,\,&z\big(\|u_{+,j}\|^2_{\tilde\calH}\delta_{j,k}
\big)_{1\leq j,k\leq n}
\lb{7.69} \\
&+(1+z^2)\big((u_+,(\widetilde H_{F(\text{resp.~}K)}
-z)^{-1}
u_+)_{\tilde\calH}\big)_{1\leq j,k\leq n},
\quad z\in\bbC_+, \no
\end{align}
where $\widetilde H_{F(\text{resp.~}K)}\geq 0$ denotes
the Friedrichs
(respectively, Krein) extension of some densely
defined closed operator $H\geq 0$ with deficiency
indices $(n,n)$
and deficiency subspace
$\{u_{+,j}\in\ker(H^{*}-i)\}_{1\leq j,k\leq n}$ in
some separable complex Hilbert space $\widetilde \calH$.

\noindent (iii). Any \, $\widetilde M_{F,K}\in
\calN_{0,F,K}^{n \times n}$ can be realized in the form
\begin{align}
\widetilde M_{F,K}(z) =
\,\,&z\big(\|u_{+,j}\|^2_{\tilde\calH}
\delta_{j,k}\big)_{1\leq j,k\leq n}
\lb{7.70} \\
&+(1+z^2)\big((u_+,
(\widetilde H_{F,K} -z)^{-1}u_+)_{\tilde\calH}
\big)_{1\leq j,k\leq n}, \quad z\in\bbC_+, \no
\end{align}
where $\widetilde H_{F,K}\geq 0$ denotes the unique
nonnegative
self-adjoint extension of some densely
defined closed operator $H\geq 0$ with deficiency
indices $(n,n)$
and deficiency subspace
$\{u_{+,j}\in\ker(H^{*}-i)\}_{1\leq j,k\leq n}$ in
some separable complex Hilbert space $\widetilde \calH$.

\noindent (iv). Any \, $\widetilde M_{F^{\bot}
(\text{resp.~}K^{\bot})}
\in
\calN_{0,F^{\bot}(\text{resp.~}K^{\bot})}^{n \times n}$
can be realized
in the form
\begin{align}
\widetilde M_{F^{\bot}(\text{resp.~}K^{\bot})}(z) =
\,\,&z\big(\|u_{+,j}\|^2_{\tilde\calH}
\delta_{j,k}\big)_{1\leq j,k\leq n}
\lb{7.71} \\
&+(1+z^2)\big((u_+,(\widetilde H_{F^{\bot}
(\text{resp.~}K^{\bot})}
-z)^{-1}
u_+)_{\tilde\calH}\big)_{1\leq j,k\leq n},
\,\, z\in\bbC_+, \no
\end{align}
where $\widetilde H_{F^{\bot}(\text{resp.~}K^{\bot})}
\geq 0$ denotes
a nonnegative self-adjoint extension of some densely
defined closed operator $H\geq 0$ with deficiency
indices $(n,n)$
and deficiency subspace $\{u_{+,j}
\in\ker(H^{*}-i)\}_{1\leq j,k\leq n}$
in
some separable complex Hilbert space $\widetilde
\calH$ such that
$\calD (\widetilde H_{F^{\bot}}\cap\calD
(H_F)=\calD (H)$ (respectively,
$\calD (\widetilde H_{K^{\bot}}\cap\calD (H_K)=
\calD (H)$).

\noindent (v). Any \, $\widetilde M_{F^{\bot},
K^{\bot}}\in
\calN_{0,F^{\bot},K^{\bot}}^{n \times n}$ can be
realized in the form
\begin{align}
\widetilde M_{F^{\bot},K^{\bot}}(z) =
\,\,&z\big(\|u_{+,j}\|^2_{\tilde\calH}
\delta_{j,k}\big)_{1\leq j,k\leq n}
\lb{7.72} \\
&+(1+z^2)\big((u_+,(\widetilde H_{F^{\bot},K^{\bot}}
-z)^{-1}u_+)_{\tilde\calH}
\big)_{1\leq j,k\leq n}, \quad z\in\bbC_+, \no
\end{align}
where $\widetilde H_{F^{\bot},K^{\bot}}\geq 0$ denotes
a nonnegative
self-adjoint extension of some densely
defined closed operator $H\geq 0$ with deficiency
indices $(n,n)$
and deficiency subspace $\{u_{+,j}
\in\ker(H^{*}-i)\}_{1\leq j,k\leq n}$
in
some separable complex Hilbert space $\widetilde \calH$
such that
$\calD (\widetilde H_{F^{\bot}}\cap\calD (H_F)=\calD (H)=
 \calD (\widetilde H_{K^{\bot}}\cap\calD (H_K)$.

\noindent In each case (i)--(v) one has
\begin{equation}
\int_{\bbR}d\widetilde \Omega(\lambda)
(1+\lambda^2)^{-1}=
\big(\|u_{+,j}\|^2_{\tilde\calH}
\delta_{j,k}\big)_{1\leq j,k\leq n},
\lb{7.73}
\end{equation}
where $\widetilde \Omega$ denotes the measure in the
Herglotz representation of $\widetilde M(z)$. Moreover,
$H$ may be
chosen prime and $\widetilde \calH$ separable.
\end{theorem}

\begin{proof}
We use the notation established in Theorem~\ref{t7.2}.
Define
\begin{equation}
u_{+,j}(\lambda)=(\lambda-i)^{-1}\ul e_j, \quad j=1,
\dots ,n, \lb{7.73a}
\end{equation}
then
\begin{align}
\|u_{+,j}\|^2_{\ti \calH}\delta_{j,k}=(u_{+,j},
u_{+,k})_{\ti \calH}
&=\int_\bbR (\ul e_j,d\widetilde \Omega(\lambda)
\ul e_k)_{\bbC^n}
(1+\lambda^2)^{-1} \no \\
&=\int_\bbR d\widetilde \Omega_{j,k}(\lambda)
(1+\lambda^2)^{-1}
\lb{7.73b}
\end{align}
and
\begin{align}
&z\|u_{+,j}\|^2_{\ti \calH}\delta_{j,k} + (1+z^2)
(u_{+,j},
(\widetilde H-z)^{-1}u_{+,k})_{\ti \calH} \no \\
&=\int_\bbR (\ul e_j,d\widetilde
\Omega(\lambda)\ul e_k)_{\bbC^n}
(z(1+\lambda^2)^{-1}+(1+z^2)(\lambda-z)^{-1}
(1+\lambda^2)^{-1}) \no \\
&=\int_\bbR d\widetilde \Omega_{j,k}(\lambda)
((\lambda-z)^{-1}
-\lambda(1+\lambda^2)^{-1}) = \widetilde
M_{j,k}(z) \lb{7.73c}
\end{align}
proves \eqref{7.68} and hence part (i). Parts (ii)
and (iii) then follow
in the same manner from Theorems~\ref{t7.2} and
\ref{t7.3}. Similarly,
parts (iv) and (v) follow from Theorems~\ref{t7.2}
and \ref{t7.5}.
\end{proof}

We also formulate the analog of Theorem~\ref{t4.5}.

\begin{theorem} \lb{t7.8}
Suppose $M_\ell\in\calN_0^{n \times n}$ with
corresponding measures
$\Omega_\ell$ in
the Herglotz representation of $M_\ell$, $\ell=1,2$,
and $M_1\neq M_2$.
Then $M_1$ and $M_2$ can be realized as
\begin{align}
M_\ell(z) = \,\,&z\big(\|u_{+,j}\|_{\calH}^2\delta_{j,k}
\big)_{1\leq j,k\leq n} \lb{7.74} \\
&+(1+z^2)\big((u_{+,j},(H_\ell -z)^{-1}u_{+,k})_{\calH}
\big)_{1\leq j,k\leq n}, \quad \ell=1,2, \,\, z\in\bbC_+,
\no
\end{align}
where $H_\ell$, $\ell=1,2$ are distinct self-adjoint
extensions of one
and the same
densely defined closed symmetric operator $H$ (which may
be chosen prime)
with deficiency indices
$(n,n)$ and deficiency subspace
$\{u_{+,j}\}_{1\leq j,k\leq n}\in\ker(H^{*}-i)$ in some
complex
Hilbert space $\calH$ (which may be chosen separable) if
and only if
the following conditions hold:
\begin{equation}
\int_{\bbR} d\Omega_1(\lambda)(1+\lambda^2)^{-1}=
\int_{\bbR} d\Omega_2(\lambda)(1+\lambda^2)^{-1}=
\big(\|u_{+,j}\|^2_{\calH}\delta_{j,k}\big)_{1\leq
j,k\leq n},
\lb{7.75}
\end{equation}
and for all $z\in\bbC_+$,
\begin{align}
&M_{2}(z)=\big(-\big(\|u_{+,j}\|_{\calH}^2\delta_{j,k}
\big)_{1\leq j,k\leq n}e^{-i\alpha_2}(\sin(\alpha_2)
\cos(\alpha_1)-
\cos(\alpha_2)\sin(\alpha_1))e^{i\alpha_1}  \no \\
& +\big(\|u_{+,j}\|_{\calH}^2\delta_{j,k}
\big)_{1\leq j,k\leq n}e^{-i\alpha_2}(\cos(\alpha_2)
\cos(\alpha_1) \no \\
&+\sin(\alpha_2)\sin(\alpha_1))
e^{i\alpha_1}\big(\|u_{+,j}\|_{\calH}^{-2}\delta_{j,k}
\big)_{1\leq j,k\leq n} M_{1}(z)\big) \no \\
& \times \big(e^{-i\alpha_2}(\cos(\alpha_2)\cos(\alpha_1)
+\sin(\alpha_2)\sin(\alpha_1))e^{i\alpha_1}  \no \\
& +e^{-i\alpha_2} (\sin(\alpha_2)\cos(\alpha_1)-
\cos(\alpha_2)\sin(\alpha_1))e^{i\alpha_1}
\big(\|u_{+,j}\|_{\calH}^{-2}\delta_{j,k}
\big)_{1\leq j,k\leq n}M_{1}(z)\big)^{-1}, \no \\
& \hspace*{5.5cm} \text{ for some }\alpha_\ell
=\alpha^*_\ell
\in M_n(\bbC), \,\, \ell=1,2. \lb{7.76}
\end{align}
\end{theorem}

\begin{proof}
Assuming \eqref{7.74}, \eqref{7.75} follows from
\begin{equation}
M_\ell(i)=i\big(\|u_{+,j}\|^2_\calH\delta_{j,k}
\big)_{1\leq j,k\leq n}
=i\int_\bbR d\Omega(\lambda)(1+\lambda^2)^{-1}, \quad
\ell=1,2, \lb{7.77}
\end{equation}
and \eqref{7.76} is clear from \eqref{7.59} upon
identifying $\alpha_1=
\alpha$, $\alpha_2=\beta$, $(\|u_{+,j}\|^{-2}_\calH$ $
\delta_{j,k})_{1\leq j,k\leq n}M_1(z)=M_\alpha(z)$, and
$(\|u_{+,j}\|^{-2}_\calH
\delta_{j,k})_{1\leq j,k\leq n}M_2(z)
=M_\beta(z)$. (Without loss of generality we may assume
that
$\{u_{+,j}\}_{1\leq j\leq n}$ is a generating basis for
$H_\ell$,
$\ell=1,2$, since otherwise we may apply the reduction
\eqref{7.48}.)
Conversely, assume \eqref{7.75} and \eqref{7.76}. By
Theorem~\ref{t7.7}\,(i), we may realize $M_1(z)$ as
\begin{align}
&\big(\|u_{+,j}\|^{-2}_\calH \delta_{j,k}\big)_{1\leq
j,k\leq n}
M_1(z)=zI_n \no \\
&+(1+z^2)\big(\|u_{+,j}\|^{-2}_\calH
\delta_{j,k}\big)_{1\leq j,k\leq n}
\big((u_{+,j},(H_1-z)^{-1}u_{+,k})\big)_{1\leq j,k\leq n}.
\lb{7.78}
\end{align}
By \eqref{7.48} we may assume that $\{u_{+,j}\}_{1\leq
j\leq n}$ is
a generating basis for $H_1$ and identify $H_1$ with
$H_\alpha$
defined in \eqref{7.28}. If $H_\beta$ is another
self-adjoint
extension of $H$ defined as in \eqref{7.28}, distinct
from $H_\alpha$,
introduce
\begin{equation}
M_\beta(z)=zI_n+(1+z^2)\big(\|u_{+,j}
\|_{\calH}^{-2}\delta_{j,k}
\big)_{1\leq j,k\leq n}\big((u_{+,j},
(H_\beta-z)^{-1}u_{+,k})
\big)_{1\leq j,k\leq n}. \lb{7.79}
\end{equation}
By \eqref{7.59} one obtains ($M_\alpha(z)
=(\|u_{+,j}\|^{-2}_\calH
\delta_{j,k})_{1\leq j,k\leq n}M_1(z)$),
\begin{align}
&M_\beta(z)=\big(-e^{-i\beta}(\sin(\beta)\cos(\alpha)-
\cos(\beta)\sin(\alpha))e^{i\alpha}  \no \\
& \quad +e^{-i\beta}(\cos(\beta)\cos(\alpha)+
\sin(\beta)\sin(\alpha))
e^{i\alpha}\big(\|u_{+,j}\|_{\calH}^{-2}\delta_{j,k}
\big)_{1\leq j,k\leq n} M_{1}(z)\big) \no \\
& \times \big(e^{-i\beta}(\cos(\beta)\cos(\alpha)
+\sin(\beta)\sin(\alpha))e^{i\alpha}  \lb{7.80} \\
& \quad +e^{-i\beta} (\sin(\beta)\cos(\alpha)-
\cos(\beta)\sin(\alpha))e^{i\alpha}
\big(\|u_{+,j}\|_{\calH}^{-2}\delta_{j,k}
\big)_{1\leq j,k\leq n}M_{1}(z)\big)^{-1}. \no
\end{align}
A comparison of \eqref{7.76} and \eqref{7.80} then
yields
$(\|u_{+,j}\|^{2}_\calH
\delta_{j,k})_{1\leq j,k\leq n}M_\beta(z)=M_2(z)$,
$\alpha=\alpha_1$,
and $\beta=\alpha_2$, completing the proof.
\end{proof}

Clearly Remark~\ref{rm4.7} applies in the present
matrix-valued context.

For different types of realization theorems in the
context of
conservative systems, see, for instance,
\cite{BGK79}--\cite{BT98}, \cite{Ts92}, and \cite{TS77}.

\vspace{1mm}

Finally we briefly turn to Hamiltonian systems on a half-line
following Hinton and Shaw \cite{HS81}--\cite{HS86} (see
also \cite{HS93},
\cite{Kr89a}, \cite{Kr89b}). These systems describe
matrix-valued
Schr\"odinger and Dirac-type differential and difference
operators (see,
e.g., \cite{Ga68}, \cite{GHL98}, \cite{GR98}, \cite{KR74},
\cite{Kr95},
and \cite{Wa74}). Let
$A, B \in L^1 ([0,R])^{2n \times 2n}$ for all $R>0$,
$A(x)=A(x)^*$,
$B(x)=B(x)^*$ for a.e.~$x>0$. Moreover, suppose that
for some
$1\leq r\leq n$, $A(x)=
\left(\begin{smallmatrix}W(x)& 0 \\ 0 &
0\end{smallmatrix}\right)$,
$W\in L^1 ([0,R])^{r \times r}$ for all $R>0$,
$W(x)>0$ for
a.e.~$x>0$. For $\ul \psi (z,\cdot)\in AC([0,R])^{2n}$
for all $R>0$,
$z\in\bbC$, consider the formally symmetric Hamiltonian
system
\begin{equation}
J_{2n}\ul \psi' (z,x)=(zA(x)+B(x))\ul \psi (z,x),
\quad x>0 \lb{7.81}
\end{equation}
and suppose Atkinson's definiteness condition
\begin{equation}
\int_a^b dx (\ul {\hat \psi}(z,x),A(x)\ul {\hat
\psi}(z,x))_{\bbC^{2n}}>0
\text{ for all } 0\leq a < b < \infty \lb{7.82}
\end{equation}
whenever $\ul {\hat \psi}$ satisfies $0\neq
\ul {\hat \psi}(z,x)
\in AC([0,R])^{2n}$ for all $R>0$ and
\eqref{7.81}. Introduce for $0\leq c < d \leq \infty$,
\begin{equation}
L^2_A ((c,d))=\{\ul f:(c,d)\to\bbC^{2n}
\text{ measurable}\,|\,
\int_c^d dx (\ul f (x),A(x)\ul f(x))_{\bbC^{2n}}
<\infty\} \lb{7.83}
\end{equation}
and for $x_0>0$,
\begin{align}
&N(z,0)=\{\ul f\in L^2_A ((0,x_0))\,|\, J_{2n}\ul
f'=(zA+B)\ul f
\text{ a.e.~on }(0,x_0)\}, \lb{7.84} \\
&N(z,\infty)=\{\ul f\in L^2_A ((x_0,\infty))\,|\,
J_{2n}\ul f'=(zA+B)\ul f
\text{ a.e.~on }(x_0,\infty)\}. \lb{7.85}
\end{align}
Then \eqref{7.81} is defined to be in the limit point
(respectively, limit
circle) case at
$e\in\{c,d\}$ if $\dim_\bbC(N(z_{\pm},e))=n$ (respectively,
$\dim_\bbC(N(z_{\pm},e))=2n$) for some (and hence for all)
$z_+\in\bbC_+$
and $z_-\in\bbC_-$. There are of course also intermediate
cases between
the limit point and limit circle case but we omit such
considerations for
simplicity. For more details in this connection see,
for instance,
\cite{KR74} and \cite{Or76}.

Next, consider $\alpha_p,\beta_p\in M_n(\bbC)$, $p=1,2$,
satisfying
$\rank(\alpha)=\rank(\beta)=n$, where
$\alpha=(\alpha_1,\alpha_2)$,
$\beta=(\beta_1,\beta_2)$ are $2n \times n$ matrices
over $\bbC$, and
\begin{equation}
\alpha_1\alpha_1^*+\alpha_2\alpha_2^*=I_n
=\beta_1\beta_1^*+\beta_2\beta_2^*,
\quad \alpha_1\alpha_2^*=\alpha_2\alpha_1^*, \quad
\beta_1\beta_2^*=\beta_2\beta_1^*. \lb{7.86}
\end{equation}
Let $\Psi_\alpha (z,x)\in M_{2n}$, $z\in\bbC$ be a
fundamental system
of solutions of \eqref{7.81} satisfying
\begin{equation}
\Psi_\alpha (z,0)=\left(\begin{array}{cc} \alpha_1^*
& -\alpha_2^*\\ \alpha_2^* & \alpha_1^*\end{array}\right),
\quad z\in\bbC \lb{7.87}
\end{equation}
and partition $\Psi_\alpha (z,x)$ into $n \times n$
blocks,
\begin{equation}
\Psi_\alpha (z,x)=\left(\begin{array}{cc}
\theta_{\alpha,1}(z,x)
& \phi_{\alpha,1}(z,x)\\ \theta_{\alpha,2}(z,x)
& \phi_{\alpha,2}(z,x)\end{array}\right). \lb{7.88}
\end{equation}
Then $\Psi_\alpha (z,x)$ is entire in $z\in\bbC$
and one defines
\begin{equation}
M_{\alpha,\beta,R}(z)=-(\beta_1\phi_{\alpha,1}(z,R)+
\beta_2 \phi_{\alpha,2}(z,R))^{-1}(\beta_1
\theta_{\alpha,1}(z,R)
+\beta_2\theta_{\alpha,2}(z,R)). \lb{7.89}
\end{equation}
$M_{\alpha,\beta,R}(z)$ is the Weyl-Titchmarsh matrix
corresponding to
the boundary value problem
\begin{align}
&J_{2n}\ul \psi' (z,x)=(zA(x)+B(x))\ul \psi (z,x),
\quad 0\leq x\leq R,
\lb{7.90} \\
&\alpha\ul \psi (z,0)=0, \,\, \beta\ul \psi (z,R)=0.
\no
\end{align}
As shown in detail by Hinton and Shaw \cite{HS81},
\cite{HS83}, \cite{HS84},
\begin{equation}
\lim_{R\uparrow\infty}M_{\alpha,\beta,R}(z)=M_\alpha(z),
\quad
z\in\bbC\backslash\bbR \lb{7.91}
\end{equation}
exists and is independent of $\beta$ if and only if
\eqref{7.81} is
in the limit point case at $\infty$. In the limit circle
case at
$\infty$, uniqueness and independence of $\beta$ is lost
and we denote
by $\hatt M_\alpha (z)$ a parametrization of all
possible limit points
of $M_{\alpha,\beta,R}(z)$ as $R\uparrow\infty$.
$M_\alpha (z)$
(respectively, $\hatt M_\alpha (z)$) are matrix-valued
Herglotz functions
with representations
\begin{equation}
M_\alpha(z)=C_\alpha+\int_\bbR d\Omega_\alpha(\lambda)
((\lambda-z)^{-1}
-\lambda(1+\lambda^2)^{-1}), \quad C_\alpha=C_\alpha^*
\lb{7.92}
\end{equation}
and one verifies
\begin{equation}
(\theta_{\alpha,p}(z,\cdot)+\phi_{\alpha,p}
(z,\cdot)M_\alpha(z))
\in L^2_A((0,\infty)), \quad p=1,2, \,\,
z\in\bbC\backslash\bbR.
\lb{7.93}
\end{equation}
Moreover,
\begin{equation}
M_\alpha(z)=(-\alpha_2+\alpha_1M(z))(\alpha_1+
\alpha_2M(z))^{-1},
\lb{7.94}
\end{equation}
where $M(z)=M_{(I_n,0)}(z)$. Analogous relationships
hold for
$\hatt M_\alpha (z)$ in the limit circle case at $\infty$.
A comparison
of \eqref{7.94} and \eqref{6.36a} suggests the
introduction of
\begin{equation}
A(\alpha)=\left(\begin{array}{cc} \alpha_1
& \alpha_2 \\ -\alpha_2 & \alpha_1 \end{array}\right)
\in{\calA}_{2n}.
\lb{7.95}
\end{equation}
In particular, Theorem~\ref{t6.6} applies (with
$A_{1,1}(\alpha)=
A_{2,2}(\alpha)=\alpha_1$, $A_{1,2}(\alpha)=-A_{2,1}
(\alpha)
=\alpha_2$).

\medskip
{\bf Acknowledgments.} We would like to thank Steve
Clark,
Henk de Snoo, Seppo Hassi, Don Hinton, and Barry Simon
for valuable
correspondence and discussions and, especially, Dirk
Buschmann and
Konstantin Makarov for a
critical reading of large portions of this manuscript.
Moreover, we
are indebted to Dirk Buschmann for communicating
Lemma~\ref{l5.2a} to us.

\appendix

\section{Examples of Scalar Herglotz Functions} \lb{A}
\renewcommand{\theequation}{A.\arabic{equation}}
\renewcommand{\thetheorem}{A.\arabic{theorem}}
\setcounter{theorem}{0}
\setcounter{equation}{0}

For convenience of the reader we collect some standard
examples
of scalar Herglotz functions and their explicit
representations
(cf.~\cite{AD56}, \cite{Bh97}, Ch.~V, \cite{Do74}, Ch.~II,
\cite{EK82}, Ch.~2).

In the following we denote Lebesgue measure on $\bbR$ by
$d\lambda$
and a pure point measure supported at $x\in\bbR$ with
mass one by
$\mu_{\{x\}}$,
\begin{equation}
\supp(\mu_{\{x\}}) =\{x\}, \quad \mu_{\{x\}}(\{x\})=1.
\lb{A.1}
\end{equation}

We start with very simple examples and progressively
discuss more
sophisticated ones.

\begin{equation}
c+id=c+d\pi^{-1}\int_{\bbR} d\lambda \,
((\lambda -z)^{-1}-\lambda(1+\lambda^2)^{-1}),
\quad c\in\bbR, \, d\geq 0. \lb{A.2}
\end{equation}
\begin{align}
\ln(id)=&\ln(d) +(i\pi /2)=\ln(d) +2^{-1}\int_{\bbR}
d\lambda \,
((\lambda -z)^{-1}-\lambda(1+\lambda^2)^{-1}), \lb{A.3} \\
&\hspace*{8.5cm} d\geq 0. \no
\end{align}
\begin{equation}
c+dz, \quad c\in\bbR, \, d\geq 0. \lb{A.4}
\end{equation}
\begin{equation}
-z^{-1}=\int_{\bbR} d\mu_{\{0\}} (\lambda) \,
(\lambda-z)^{-1}. \lb{A.5}
\end{equation}
\begin{align}
\ln(z)&=\int_{-\infty}^{0} d\lambda \,
((\lambda -z)^{-1}-\lambda(1+\lambda^2)^{-1}), \quad
 \lb{A.6} \\
\ln(-z^{-1})&=\int_{0}^{\infty} d\lambda \,
((\lambda -z)^{-1}-\lambda(1+\lambda^2)^{-1}), \lb{A.7}
\end{align}
where $\ln(\cdot)$ denotes the principal value of the
logarithm (i.e., with cut along $(-\infty,0]$ and
$\ln(\lambda)>0$ for $\lambda>0$).
\begin{align}
z^r&=\exp(r\ln(z)) \no \\
&=\cos(r\pi /2)+\pi^{-1}\sin(r\pi)\int_{-\infty}^{0}
d\lambda \, |\lambda|^r
((\lambda -z)^{-1}-\lambda(1+\lambda^2)^{-1}), \lb{A.8} \\
& \hspace*{8.5cm} 0<r<1, \no \\
-z^{-r}&=-\exp(-r\ln(z)) \no \\
&=-\cos(r\pi /2)+\pi^{-1}\sin(r\pi)\int_{-\infty}^{0}
d\lambda \, |\lambda|^{-r}
((\lambda -z)^{-1}-\lambda(1+\lambda^2)^{-1}), \lb{A.9} \\
& \hspace*{9.1cm} 0<r<1. \no
\end{align}
\begin{align}
\tan(z)&=\sum_{n\in\bbZ}(((n+\tfrac12) \pi-z)^{-1}-
(n+\tfrac12)\pi
(1+(n+\tfrac12 )^2\pi^2)^{-1}) \no \\
&=\int_{\bbR} d\omega(\lambda) \,
((\lambda -z)^{-1}-\lambda(1+\lambda^2)^{-1}), \lb{A.10} \\
\omega &=\sum_{n\in\bbZ} \mu_{\{(n+\tfrac12)\pi\}},
\lb{A.10a} \\
-\cot(z)&=\sum_{n\in\bbZ}((n \pi-z)^{-1}-n\pi
(1+n^2\pi^2)^{-1}) \no \\
&=\int_{\bbR} d\omega(\lambda) \,
((\lambda -z)^{-1}-\lambda(1+\lambda^2)^{-1}), \lb{A.11} \\
\omega &=\sum_{n\in\bbZ} \mu_{\{n\pi\}}. \lb{A.11a}
\end{align}
The psi or digamma function,
\begin{align}
\psi(z)=\Gamma'(z)/\Gamma(z)&=C+\sum_{n\in\bbN_0}((-n-z)^{-1}
+n(1+n^2)^{-1}) \no \\
&=\int_{\bbR} d\omega(\lambda) \,
((\lambda -z)^{-1}-\lambda(1+\lambda^2)^{-1}),
\lb{A.12} \\
& \hspace*{-2.8cm} \omega =\sum_{n\in\bbN_0}\mu_{\{-n\}},
\quad
C=-\gamma + \sum_{n\in\bbN_0}((n+1)^{-1}-n(1+n^2)^{-1}).
\lb{A.12a}
\end{align}
Here $\Gamma(z)$ denotes the gamma function,
$\gamma=-\psi(1)=.572\dots$
Euler's constant (cf.~\cite{AS72}, Ch.~6), and
$\bbN_0=\bbN\cup\{0\}$.
\begin{align}
\frac{z-\lambda_2}{z-\lambda_1} &=1+
(\lambda_2-\lambda_1)\int_{\bbR}d\mu_{\{\lambda_1\}}
(\lambda)
\, (\lambda-z)^{-1}, \quad \lambda_1<\lambda_2,
\lb{A.13} \\
\ln \bigg(\frac{z-\lambda_2}{z-\lambda_1} \bigg)
&=\int_{\lambda_1}^{\lambda_2} d\lambda
\, (\lambda-z)^{-1}, \quad \lambda_1<\lambda_2. \lb{A.14}
\end{align}

Next we turn to Weyl-Titchmarsh m-functions
$m_{a(\alpha)}(z)$ associated with the
operator $H_\alpha$ in $L^2([0,\infty);dx)$ defined by
\begin{align}
&(H_\alpha g)(x)=-g''(x), \quad x>0, \no \\
&\calD (H_\alpha)=
\{g\in L^2([0,\infty);dx)\,|\,g,g'\in AC([0,R])
\text{ for all $R>0$ }; \lb{A.14a} \\
& \hspace{1.6cm} -g''\in L^2([0,\infty);dx); \,\,
 \sin(\alpha)g'(0_+)+\cos(\alpha)g(0_+)=0 \}, \quad
\alpha\in (0,\pi). \no
\end{align}
These are
Herglotz functions of the type
\begin{equation}
m_{a(\alpha)}(z)=\frac{-\sin(\alpha)+\cos(\alpha)iz^{1/2}}
{\cos(\alpha)+\sin(\alpha)iz^{1/2}}=\cot(\alpha)
+\int_{\bbR} d\omega_{a(\alpha)}(\lambda) \,
(\lambda-z)^{-1}, \lb{A.15}
\end{equation}
\begin{equation}
\omega_{a(\alpha)}(\lambda)=\begin{cases} 0, \,
\lambda <0,
\, (\pi /2) \leq \alpha < \pi  \\
-2\frac{\cot(\alpha)}{\sin^2(\alpha)}, \,
-\infty<\lambda<-\cot^2(\alpha), \, 0<\alpha<(\pi /2) \\
0, \,  -\cot^2(\alpha) < \lambda <0, \, 0<\alpha<(\pi /2) \\
\frac{2}{\pi}\lambda^{1/2}, \, \lambda\geq 0, \,
\alpha=(\pi /2) \\
\frac{2}{\pi \sin^2(\alpha)}\big(\lambda^{1/2}
-\cot(\alpha)\arctan \big(\frac{\lambda^{1/2}}
{\cot(\alpha)}\big)\big),
\, \lambda\geq 0, \,
\alpha\in(0,\pi)\backslash\{\pi /2\}, \end{cases}
\lb{A.15a}
\end{equation}
where (cf.~\eqref{2.11} and \eqref{3.10})
\begin{equation}
a(\alpha)=\left(\begin{array}{cc} \cos(\alpha)
& \sin(\alpha)\\ -\sin(\alpha) & \cos(\alpha)
\end{array}\right)
\in \calA_2. \lb{A.16}
\end{equation}
Similarly,
\begin{equation}
m_{a(0)}=iz^{1/2}=-2^{-1/2} +\pi^{-1}\int_0^\infty
d\lambda \,\lambda^{1/2}
((\lambda -z)^{-1}-\lambda(1+\lambda^2)^{-1}) \lb{A.17}
\end{equation}
corresponds to the remaining self-adjoint (Friedrichs)
boundary condition
$\alpha=0$, that is, to $g(0_+)=0$ in \eqref{A.14a}.

Finally we describe a class of Herglotz functions
fundamental in
Floquet theory of periodic Schr\"odinger operators
on $\bbR$.
Consider a sequence $\{\lambda_n\}_{n\in\bbN_0}
\subset\bbR$,
\begin{equation}
0=\lambda_0<\lambda_1\leq\lambda_2<\lambda_3
\leq\lambda_4<\cdots
\lb{A.20}
\end{equation}
such that asymptotically
\begin{equation}
\lambda_{2n}, \, \lambda_{2n-1} \underset{n\to\infty}{=}
(n\pi)^2+\Oh(1). \lb{A.21}
\end{equation}
Define an entire function $\Delta(z)$ such that
\begin{align}
\Delta(z)-1&=\frac{(\lambda_0-z^2)}{2}\prod_{n\in\bbN}
\frac{(\lambda_{4n-1}-z^2)(\lambda_{4n}-z^2)}{(2n\pi)^4},
\lb{A.22} \\
\Delta(z)+1&=2\prod_{n\in\bbN_0}
\frac{(\lambda_{4n+1}-z^2)(\lambda_{4n+2}-z^2)}
{((2n+1)\pi)^4} \lb{A.23}
\end{align}
and hence,
\begin{equation}
\Delta(z)^2-1=(\lambda_0-z^2)\prod_{n\in\bbN}
\frac{(\lambda_{2n-1}-z^2)(\lambda_{2n}-z^2)}{(n\pi)^4}.
\lb{A.24}
\end{equation}
Moreover, define
\begin{equation}
\theta(z)=-\int_0^z d\zeta \, \frac{\Delta'(\zeta)}
{(1-\Delta(\zeta)^2)^{1/2}}, \quad z\in\bbC_+,
\lb{A.25}
\end{equation}
where the square root branch in \eqref{A.25} is
chosen to be
positive on the interval $(0,\lambda_1^{1/2})$.
Then
\begin{equation}
\cos(\theta(z))=\Delta(z) \lb{A.26}
\end{equation}
and, as shown in
\cite{Ma86}, Sect.~3.4, $\theta$ is a Herglotz
function with a
representation of the type
\begin{equation}
\theta(z)=c+ z + \pi^{-1} \int_{\bbR} d\lambda \,
\Im(\theta(\lambda) ((\lambda-z)^{-1}-
\lambda(1+\lambda^2)^{-1})
\lb{A.27}
\end{equation}
for some $c\in\bbR$. In the case where the sequence
$\{\lambda_n\}_{n\in\bbN_0}$
represents the periodic and antiperiodic eigenvalues
associated
with a Schr\"odinger operator $H=-\tfrac{d^2}{dx^2}+q$,
with $q\in L^1_{\loc}(\bbR)$ real-valued and of
period one,
$\Delta(z)$ represents the corresponding Floquet
discriminant
and $\theta(z)$ the Floquet (Bloch)
momentum associated with $H$. In this case one verifies
(see, e.g., \cite{JM82}, \cite{Ko85})
\begin{equation}
\theta(z)=\frac{i}{2} \int_0^1 dx \, G(z^2,x,x)^{-1},
\quad z\in\bbC_+, \lb{A.28}
\end{equation}
with $G(\zeta,x,y)=(H-\zeta)^{-1}(x,y)$ the Green's
function of $H$.

Analogous observations apply to one-dimensional Dirac-type
operators.

\section{Krein's Formula and Linear Fractional
Transformations} \lb{B}
\renewcommand{\theequation}{B.\arabic{equation}}
\renewcommand{\thetheorem}{B.\arabic{theorem}}
\setcounter{theorem}{0}
\setcounter{equation}{0}

The main purpose of this appendix is to provide a proof
of \eqref{7.59} (cf.~Theorem~\ref{t5})
following its derivation in \cite{GMT97}. Our method of
proof is based
on Krein's formula,
which describes the resolvent difference of two self-adjoint
extensions
$A_1$ and $A_2$
of a densely defined closed symmetric linear operator
$A$ with
deficiency indices $(n,n)$, $n\in \bbN$. (Reference
\cite{GMT97} treats
this topic in the general case where
$n\in\bbN\cup\{\infty\}$.
Here we
specialize to the case $n<\infty$.) Since the latter
formula is
interesting in its own right we start with the basic
setup following
\cite{AG93}.

Let $\calH$ be a complex separable Hilbert space,
$\dot A:\calD(\dot A)\to\calH, \,\, \ol{\calD(\dot A)}
=\calH$ a
densely defined closed symmetric linear operator in
$\calH$ with finite
and equal deficiency indices
 ${\rm def}(\dot A)=(r,r),$
$r\in \bbN$. Let $A_\ell, \, \ell=1,2$,  be two distinct
self-adjoint
extensions
of $\dot A$  and denote by $A$ the maximal common part
of $A_1$ and $A_2$,
that is, $A$ is the largest closed extension of
$\dot A$  with
$\calD(A)=\calD(A_1)\cap \calD(A_2).$ Let
$0\le p\le r-1$ be the
maximal number of elements in $\calD(A)=\calD(A_1)
\cap \calD(A_2)$
which are linearly independent modulo $\calD(\dot A)$.
Then $A$ has
 deficiency indices ${\rm def} (A)=(n,n)$, $n=r-p$.
Next, denote by
 $\ker (A^*-z), \, z\in \bbC\backslash\bbR$ the deficiency
subspaces
of $A$ and define
\begin{equation}\label{e1}
U_{1,z,z_0}=I+(z-z_0)(A_1-z)^{-1}=(A_1-z_0)(A_1-z)^{-1},
\,\,z,z_0\in
\rho (A_1),
\end{equation}
where $I$ denotes the identity operator in $\calH$ and
$\rho(T)$
abbreviates the resolvent set of $T$. One verifies
\begin{equation}\label{e2}
U_{1,z_0,z_1}U_{1,z_1,z_2}=U_{1,z_0,z_2},\,\, z_0, z_1,
z_2\in \rho (A_1)
\end{equation}
and
\begin{equation}\label{e3}
U_{1,z,z_0}\ker(A^*-z_0)=\ker(A^*-z).
\end{equation}
Let $\{ u_j(i) \}_{1\le j \le n}$ be an orthonormal basis
for
$ \ker(A^*-i)$ and define
\begin{equation}\label{e4}
u_{1,j}(z)=U_{1,z,i}\,u_{j}(i), \,\, 1\le j\le n, \,\,
z\in \rho(A_1).
\end{equation}
Then $\{ u_{1,j}(z) \}_{1\le j \le n}$ is a basis for
$\ker(A^*-z)$,
 $z\in \rho(A_1)$ and since
$ U_{1,-i,i}=(A_1-i)(A_1+i)^{-1}$  is the unitary
Cayley transform of $A_1$,
$\{  u_{1,j}(-i) \}_{1\le j \le n}$ is in fact an
orthonormal basis for
 $\ker(A^*+i)$.

The basic result on Krein's formula, as presented by
Akhiezer and
 Glazman \cite{AG93}, Sect. 84, then reads as follows.

\begin{theorem} \mbox{\rm (Krein's formula, \cite{AG93},
Sect. 84)} \lb{t1} \\
There exists a $P_{1,2}(z)=(P_{1,2,j,k}(z))_{1\le j,k
\le n}\in M_n(\bbC)$,
$z\in \rho(A_2)\cap\rho(A_1)$, such that
\begin{align}
&\det (P_{1,2}(z))\ne 0, \,\,z\in \rho(A_2)\cap\rho(A_1),
\label{e5} \\
&P_{1,2}(z)^{-1}=P_{1,2}(z_0)^{-1}-(z-z_0)(u_{1,j}
(\bar z),u_{1,k}( z_0)),
\,\,z,z_0\in \rho(A_1), \label{e6} \\
&\Im \, (P_{1,2}(i)^{-1})=-I_n, \label{e7} \\
&(A_2-z)^{-1}=(A_1-z)^{-1}+
\sum_{j,k=1}^nP_{1,2,j,k}(z)(u_{1,k} (\bar z), \cdot \,
) u_{1,j} (z),
\quad z\in \rho(A_2)\cap\rho(A_1). \label{e8}
\end{align}
\end{theorem}
We note that $P_{1,2}(z)^{-1}$ extends by continuity
from $z\in\rho(A_2)\cap\rho(A_1)$ to all of $\rho(A_1)$
since the right-hand
side of
(\ref{e6}) is continuous for $z\in \rho(A_1)$. The
normalization condition
(\ref{e7}) is not mentioned  in \cite{AG93} but
it trivially follows
from (\ref{e6}) and the fact
\begin{equation}
(u_j(i), u_k(i))=\delta_{j,k}, \quad 1\le j,k \le n \lb{e8a}
\end{equation}
(where $\delta_{j,k}$ denotes Kronecker's symbol)
and from
\begin{equation}
P_{1,2}^*(z)=P_{1,2}(\bar z), \quad z\in\rho(A_1)
\cap\rho(A_2). \lb{e8b}
\end{equation}
Taking $z=\overline{z_0}$ in (\ref{e6}) shows that
$-P_{1,2}(z)^{-1}$
and hence
 $P_{1,2}(z)$ is a matrix-valued Herglotz function,
that is,
\begin{equation}\label{e8c}
\Im\, (P_{1,2}(z))>0, \,\, z\in \bbC_+.
\end{equation}
Strict positive
definiteness in (\ref{e8c}) follows from the fact that
$\{u_{1,k}(z)\}_{1\le k \le n}$ are linearly independent
for $z\in \bbC_+$
and hence
$((u_{1,j}(z), u_{1,k}(z))_{1\le j,k \le n} >0.$

Krein's
formula has been used in a great variety of
problems in mathematical physics as can be seen from the
extensive number of
  references provided, for instance, in \cite{AGHKH88}.
(A complete
bibliography on Krein's formula is impossible in this
context.)

Next we describe the connection between $P_{1,2}(z)$ and
von Neumann's
parametrization of self-adjoint extensions of $A$. Due to
(\ref{e6}), $P_{1,2}(z)^{-1}$ is determined for all
$z\in \rho(A_1)$
in terms of
$P_{1,2}(i)^{-1}, \,\,(A_1-z)^{-1}$ and
$\{ u_j(i)\}_{1 \le j \le n} $,
\begin{align}
&P_{1,2}(z)^{-1}=P_{1,2}(i)^{-1}-(z-i)I_n-(1+z^2)(u_j(i),
(A_1-z)^{-1}u_k(i))_{1 \le j,k \le n}), \no \\
& \hspace*{9.5cm} z\in \rho(A_1). \label{e9}
\end{align}
Hence it suffices to focus on
\begin{equation}
P_{1,2}(i)^{-1}=\Re\,(P_{1,2}(i)^{-1})-iI_n. \lb{e9a}
\end{equation}
Let
\begin{equation}\label{e10}
\calU_\ell:\ker (A^*-i) \to \ker (A^*+i),
\quad \ell=1,2,
\end{equation}
be the linear isometric isomorphisms that parameterize
$A_\ell$ according
to von Neumann's formula
\begin{align}
&A_\ell(f+(I+\calU_\ell)u_+)=Af+i(I-\calU_\ell)u_+,
\lb{e11} \\
&\calD(A_\ell)=\{ (g+(I+\calU_\ell)u_+)\in \calD(A^*) \,
\vert \, g\in \calD(A), \, u_+\in \ker (A^*-i) \},
\,\,\ell=1,2. \no
\end{align}
Denote by $U_\ell=(U_{\ell,j,k})_{1\le j,k \le n}\in
M_n(\bbC),\,\,  \ell=1,2$
the unitary matrix representation of $\calU_\ell$ with
respect to the bases
$\{ u_j(i) \}_{1 \le j \le n}$ and
$\{ u_{1,j}(-i) \}_{1 \le j \le n}$
of $\ker (A^*-i)$ and $\ker (A^*+i)$ respectively, that is,
\begin{equation}\label{e12}
\calU_\ell u_j(i)=\sum_{k=1}^n U_{\ell,k,j}u_{1,k}(-i),
\,\, 1\le j \le n, \,\,
\ell=1,2.
\end{equation}

\begin{lemma} \mbox{\rm } \lb{l2}

\noindent (i). $U_1=-I_n.$

\noindent (ii). $(u_{1,j}(-i) + \sum_{k=1}^n
\ol{U_{\ell,j,k }}u_k(i) )
\in \calD(A_\ell),
\,\, 1 \le j \le n, \,\,\ell=1,2, $ and
\begin{equation}\label{e13}
(A_\ell-i)(u_{1,j}(-i) +\sum_{k=1}^n
\ol{U_{\ell,j,k}}u_k(i))=
-2i u_{1,j}(-i), \,\,1
\le j \le n, \,\,\ell=1,2.
\end{equation}
\end{lemma}
\begin{proof} (i). Since  $u_{1,j}(-i)-u_{j}(i)=-2i
(A_1+i)^{-1}u_j(i)
\in \calD(A_1), \,\,1 \le j \le n$
by (\ref{e1}) and (\ref{e4}), one infers
\begin{equation}
\label{e14}
u_{1,j}(-i)-u_j(i)=c_j(I+\calU_1)u_{+,j}, \,\,1 \le j
\le n
\end{equation}
for some $u_{+,j}\in \ker (A^*-i)$ and $c_j\in \bbC$.
Since $\calD(A^*)$
 decomposes into the direct sum
 $\calD(A^*)=\calD(A) \dot + \ker (A^*-i) \dot +
\ker (A^*+i), $
one infers
\begin{equation}\label{e15}
c_ju_{+,j}=-u_j(i), \quad
c_j\calU_1 u_{+,j}=-\calU_1 u_j(i)=u_{1,j}(-i),
\end{equation}
and hence $U_1=-I_n.$

\noindent (ii). Using (\ref{e12}) one computes
\begin{equation}
\calU_\ell \sum_{k=1}^n\ol{U_{\ell,j,k}}u_k(i)
=\sum_{k=1}^n\sum_{m=1}^n
\ol{U_{\ell,j,k}} U_{\ell,m,k}u_{1,m}(-i)=u_{1,j}(-i),
\,\, 1 \le j \le n,
\lb{e16}
\end{equation}
utilizing unitarity of $U_\ell, \,\ell=1,2$. Hence,
\begin{equation}\label{e17}
u_{1,j}(-i)+\sum_{k=1}^n
\ol{U_{\ell,j,k}} u_k(i) =
(I+\calU_\ell) (\sum_{k=1}^n
\ol{U_{\ell,j,k}} u_k(i)) \in \calD(A_\ell),
\end{equation}
and thus (\ref{e13}) follows from
\begin{align}
(A_\ell-i)(u_{1,j}(-i)+\sum_{k=1}^n
\ol{U_{\ell,j,k}} u_k(i)) &=
(A^*-i)(u_{1,j}(-i)+\sum_{k=1}^n
\ol{U_{\ell,j,k}} u_k(i)) \no \\
&=-2iu_{1,j}(-i).  \label{e18}
\end{align}
\end{proof}

This yields the desired connection between $P_{1,2}(i)$
and $U_\ell, \,\,\ell=1,2.$

\begin{corollary} \mbox{\rm } \lb{c3}
\begin{equation}\label{e19}
P_{1,2}(i)=\frac{i}{2} (I_n+U_2^{-1})=\frac{i}{2}(U_2^{-1}
-U_1^{-1}).
\end{equation}
\end{corollary}
\begin{proof} By (\ref{e13}), \eqref{e8}, and (\ref{e8a}),
\begin{equation}\label{e20}
((A_2-i)^{-1}-(A_1-i)^{-1})u_{1,j}(-i)=\frac{i}{2}\sum_{k=1}^n
(\delta_{j,k}
+\ol{U_{2,j,k}})u_k(i)=\sum_{k=1}^n P_{1,2,k,j}(i)u_k(i).
\end{equation}
Unitarity of $U_2$ and linear independence of the $u_k(i)$
then complete
the proof of (\ref{e19}).
\end{proof}

Finally, we turn to our main goal, the Weyl-Titchmarsh
$M$-matrices
$M_1(z)$ and $M_2(z)$ associated with $A_1$
and $A_2$. Define (cf. \cite{Do65}, \cite{GMT97},
\cite{Ma92})

\begin{equation}\label{e21}
M_\ell(z)=zI_n +(1+z^2) \big((u_j(i), (A_\ell-z)^{-1}
u_k(i))_{1\le j,k \le n}\big), \quad z\in \rho(A_\ell),
\,\, \ell=1,2.
\end{equation}

$M_\ell(z)$ as defined in (\ref{e21}) are known to be
matrix-valued
Herglotz functions. More precisely, one can prove

\begin{lemma} \mbox{\cite{GMT97}} \lb{l3a}
Assume $A_1$ to be a self-adjoint extension of $A$. Then
the Weyl-Titchmarsh
matrix $M_1 (z)$ is analytic for $z\in \bbC\backslash
\bbR$ and
\begin{equation}\label{e21a}
\Im(z) \Im(M_1 (z)) \geq (\max (1,|z|^2)+|\Re(z)|)^{-1},
\quad z\in \bbC\backslash \bbR.
\end{equation}
In particular, $M_1 (z)$ is an $n \times n$ matrix-valued
Herglotz
function.
\end{lemma}
\begin{proof}
 Using (\ref{e21}), an explicit computation yields
\begin{align}
& \Im (z) \Im (M_1 (z))  \label{e21b} \\
&=\big((u_j(i), (I+A_1^2)^{1/2}
((A_1-\Re (z))^2 +(\Im(z))^2)^{-1}(I+A_1^2)^{1/2}u_k(i))
\big)_{1\leq j,k\leq n}. \no
\end{align}
Next we note that for $z\in\bbC\backslash\bbR$,
\begin{equation}
\frac{1+\lambda^2}{(\lambda-\Re(z))^2+(\Im(z))^2}\geq
\frac{1}{\max (1,|z|^2)+|\Re(z)|}, \quad \lambda\in\bbR.
\lb{e21c}
\end{equation}
Since by the Rayleigh-Ritz technique, projection onto a
subspace contained
in the domain of a self-adjoint operator bounded from
below can only
raise the lowerbound of the spectrum (cf.
\cite{RS78}, Sect.~XIII.1), \eqref{e21b} and \eqref{e21c}
prove
\eqref{e21a}.
\end{proof}

Combining \eqref{e9}, \eqref{e9a}, and \eqref{e21} for
$\ell=1$ yields
\begin{equation}
P_{1,2}(z)^{-1}=\Re\,((P_{1,2}(i)^{-1})-M_1(z). \lb{e21aa}
\end{equation}

One infers the following result relating $M_1(z)$ and
$M_2(z)$.

\begin{theorem} \mbox{\rm \cite{GMT97}} \lb{t4}
\begin{align}
M_2(z)&=(P_{1,2}(i)+(I_n+iP_{1,2}(i))M_1(z))
((I_n+iP_{1,2}(i))-
P_{1,2}(i)M_1(z))^{-1} \label{e22} \\
&=e^{-i\alpha_2}(\cos (\alpha_2)+\sin (\alpha_2)M_1(z) )
(\sin (\alpha_2)-\cos (\alpha_2)M_1(z) )^{-1}e^{i\alpha_2},
\label{e23}
\end{align}
where $\alpha_2\in M_n(\bbC)$ denotes a self-adjoint
matrix related to
$U_2$ by $U_2=e^{2i\alpha_2}.$
\end{theorem}
\begin{proof}
Using
\begin{align}
&(u_j(i), u_k(i))=\delta_{j,k}, \label{e24} \\
&(u_j(i), u_{1,k}(z))=\delta_{j,k}+(z-i)(u_j(i),
(A_1-z)^{-1} u_k(i)),
\label{e25} \\
&(u_{1,j}(\bar z), u_k(i))=\delta_{j,k}+(z+i)(u_j(i),
(A_1-z)^{-1} u_k(i)),
\label{e26}
\end{align}
and Krein's formula (\ref{e8}), one infers
\begin{align}
M_{2,j,k}(z)&=z\delta_{j,k}+(1+z^2)(u_j(i),
(A_2-z)^{-1} u_k(i)) \no \\
&=z\delta_{j,k}+(1+z^2)(u_j(i),(A_1-z)^{-1} u_k(i)) \no \\
& \quad +\sum_{s,t=1}^n (u_j(i),u_{1,s}(z))
P_{1,2,s,t}(z)(u_{1,t}(\bar z),
u_k(i))(1+z^2) \no \\
&=M_{1,j,k}(z)+
\sum_{s,t=1}^n((z+i)\delta_{j,s}+(1+z^2)(u_j(i),
(A_1-z)^{-1}u_s(i)) \no \\
& \quad \times(P_{1,2}(i)^{-1}-(z-i)I_n-(1+z^2)((u_p(i),
(A_1-z)^{-1}u_q(i))_{1\le p,
q \le n}))^{-1}_{s,t} \no \\
& \quad \times((z-i)\delta_{t,k} +(1+z^2)(u_t(i),
(A_1-z)^{-1}u_k(i))). \label{e27}
\end{align}
Hence,
\begin{equation}\label{e28}
M_2(z)=M_1(z)+(iI_n+M_1(z))(P_{1,2}(i)^{-1}
+iI_n-M_1(z))^{-1}(-iI_n+M_1(z)),
\end{equation}
which easily reduces to (\ref{e22}). Equation (\ref{e23})
then follows from
 (\ref{e19}) and the elementary  trigonometric identity
\begin{equation}
\label{e29}
\Re\,\,(P_{1,2}(i)^{-1} ) = \tan(\alpha_2), \,\,\,
U_2=e^{2i\alpha_2}.
\end{equation}
\end{proof}

Equation \eqref{e23} is connected with the pair
$(U_2,U_1)=
(e^{2i\alpha_2},-I_n)$. If one is interested in a
general pair of
self-adjoint extensions $(A_\alpha,A_\beta)$ of $A$,
associated with
$(U_\alpha,U_\beta)$, one proceeds as follows:

\begin{theorem} \mbox{\cite{GMT97}} \lb{t5}
Let $U_\alpha=e^{2i\alpha}, U_\beta=e^{2i\beta}\in
M_n(\bbC)$ be the
matrix representations of the operators $\calU_\alpha,
\calU_\beta$
associated with two self-adjoint
extensions $A_\alpha, A_\beta$ of $A$ with respect to
the basis
$\{u_j(i)\}_{1\leq j\leq n}$ of $\ker(A^*-i)$ and any (not
necessarily
orthogonal) basis $\{v_j\}_{1\leq j\leq n}$ of $\ker(A^*+i)$.
Then
\begin{align}
M_{\beta}(z)&=(-e^{-i\beta}(\sin(\beta)\cos(\alpha)-
\cos(\beta)\sin(\alpha))e^{i\alpha}  \no \\
&\quad \quad +e^{-i\beta}(\cos(\beta)\cos(\alpha)
+\sin(\beta)\sin(\alpha))
e^{i\alpha} M_{\alpha}(z)) \no \\
&\quad \times (e^{-i\beta}(\cos(\beta)\cos(\alpha)
+\sin(\beta)\sin(\alpha))e^{i\alpha}  \no \\
&\quad \quad \,\,\,\, +e^{-i\beta} (\sin(\beta)\cos(\alpha)-
\cos(\beta)\sin(\alpha))e^{i\alpha}M_{\alpha}(z))^{-1}.
\lb{e30}
\end{align}
\end{theorem}
\begin{proof}
We start by proving \eqref{e30} with respect to the
orthogonal bases
$\{u_j(i)\}_{1\leq j\leq n}$ and $\{u_{1,j}(-i)\}_{1\leq
j\leq n}$
applying Theorem~\ref{t4}.
Assuming that the pairs $(A_\alpha,A_1)$ and
$(A_\beta,A_1)$
are
relatively prime, one infers from \eqref{e23},
\begin{align}
M_\alpha(z)&=e^{-i\alpha}(\cos (\alpha)+
\sin (\alpha)M_1(z) )
(\sin (\alpha)-\cos (\alpha)M_1(z) )^{-1}e^{i\alpha},
\label{e31} \\
M_\beta(z)&=e^{-i\beta}(\cos (\beta)+\sin (\beta)M_1(z) )
(\sin (\beta)-\cos (\beta)M_1(z) )^{-1}e^{i\beta},
\label{e32}
\end{align}
Computing $M_1(z)$ (corresponding to $A_1$ and
$U_1=-I_n$) from
\eqref{e31} yields
\begin{equation}
M_1(z)=-e^{i\alpha}(\cos(\alpha)-\sin(\alpha)
M_\alpha(z))
(\sin(\alpha)+\cos(\alpha)M_\alpha(z))^{-1}
e^{-i\alpha}. \lb{e33}
\end{equation}
Insertion of \eqref{e33} into \eqref{e32} then
proves \eqref{e30}.
Inspection of the eight trigonometric terms in
\eqref{e30} shows that
they are of the type $(c_1I_n+c_2U_\beta^{-1})
(c_3I_n+c_4U_\alpha)$ with
$c_m\in\{\pm1/4,\pm i/4 \}$, $1\leq m\leq 4$. That
is, they are matrix
representations of $F_{\alpha,\beta}:=(-c_1\calU_1^{-1}
+c_2\calU_\beta^{-1})
(-c_3\calU_1+c_4\calU_\alpha)$. But $F_{\alpha,\beta}$
map $\ker(A^*-i)$
into itself, and hence matrix representations of
$F_{\alpha,\beta}$ are
independent of the basis chosen in $\ker(A^*+i)$.
\end{proof}

The material of this appendix in the general case
where ${\rm def}(A)
=(n,n)$, $n\in\bbN\cup\{\infty\}$, is considered
in detail in \cite{GMT97}.

Since the boundary values $\lim_{\varepsilon
\downarrow 0}M_\alpha
(\lambda+i\varepsilon)$, $\lambda\in\bbR$, contain
spectral information
on the self-adjoint extension $A_\alpha$ of $A$,
relationships of the type
\eqref{e30} entail important connections between the
spectra of
$A_\alpha$ and $A_\beta$. In particular, the well-known
unitary equivalence
of the absolutely continuous parts $A_{\alpha,ac}$ and
$A_{\beta,ac}$ of
$A_\alpha$ and $A_\beta$ can be inferred from \eqref{e30}
as discussed in
detail in Section~\ref{s7}. \\

We conclude with a simple illustration.

\begin{example}
$\calH=L^2((0,\infty);dx)$,
\begin{align}
&A=-\frac{d^2}{dx^2}, \no \\
&\calD(A)=\{g\in L^2((0,\infty);dx) \, \vert \,
g,g'\in AC_{\loc}((0,\infty)), \,\, g(0_+)=g'(0_+)=0 \},
\no \\
&A^*=-\frac{d^2}{dx^2}, \no \\
&\calD (A^*)=\{g\in L^2((0,\infty);dx) \, \vert \,
g,g'\in AC_{\loc}((0,\infty)),\, g''\in
L^2((0,\infty);dx) \}, \no \\
&A_1=A_F=-\frac{d^2}{dx^2}, \quad
\calD(A_1)=\{g\in \calD(A^*) \,\, \vert\, g(0_+)=0 \},
\no \\
&A_2=-\frac{d^2}{dx^2}, \quad
\calD(A_2)=\{g\in \calD(A^*) \,\, \vert \,
g'(0_+)+2^{-1/2}(1-
\tan (\alpha_2))g(0_+)=0  \}, \no \\
& \hspace*{8.9cm} \alpha_2\in [0,\pi)
\backslash\{\pi/2\}, \no
\end{align}
where $A_F$ denotes the Friedrichs extension of $A$
(corresponding to
$\alpha_2=\pi/2$). One verifies,
\begin{align}
&\ker(A^*-z)=\{ce^{i\sqrt{z}x}, c \in \bbC \},
\,\, \Im \, (\sqrt{z}) >0,
\, z\in \bbC\setminus [0, \infty), \no \\
&{\rm def} (A)=(1,1), \quad u_1(i,x)=2^{1/4}
e^{i\sqrt{i}x}, \quad
u_{1,1}(-i,x)=2^{1/4} e^{i\sqrt{-i}x}, \no \\
&(A_2-z)^{-1}=(A_1-z)^{-1}
-(2^{-1/2}(1-\tan(\alpha_2)) +i\sqrt{z})^{-1}
(\ol{e^{i\sqrt{z} \cdot}}
, \, \cdot \, )e^{i\sqrt{z} \cdot}, \no \\
& \hspace*{7.2cm} z\in \rho(A_2), \,\, \Im \,
(\sqrt{z})>0, \no \\
&U_1=-1, \quad U_2=e^{2i\alpha_2}, \no \\
&P_{1,2}(z)=-(1-\tan(\alpha_2) +i\sqrt{2z})^{-1},
\,\, z\in \rho(A_2), \quad
P_{1,2}(i)^{-1}=\tan(\alpha_2)-i, \no \\
&M_1(z)=i\sqrt{2z} +1, \quad
M_2(z)=\frac{\cos(\alpha_2)+\sin(\alpha_2)(i\sqrt{2z} +1)}
{\sin(\alpha_2)-\cos(\alpha_2)(i\sqrt{2z} +1)}. \no
\end{align}
\end{example}

The Krein extension $A_2=A_K$ of $A$ corresponds to
$\tan (\alpha_2)=1$ and hence coincides with the
Neumann extension $A_N$ of $A$ (characterized by the
boundary condition $g'(0_+)=0$).


\end{document}